\documentclass{article}
\usepackage[utf8]{inputenc}

\usepackage{multirow}
\usepackage[letterpaper, margin=1in]{geometry} 

\usepackage{amssymb}
\usepackage{latexsym}
\usepackage{amsmath,amstext,amsthm,amssymb,bbm,amsbsy}
\usepackage{mathtools}
\usepackage{tabularx}
\usepackage{url}

\usepackage[letterpaper, margin=1in]{geometry} 
\usepackage[ruled,linesnumbered]{algorithm2e}
\usepackage{subcaption}

\usepackage{tikz}
\usetikzlibrary{shapes.geometric, arrows}
\usepackage{adjustbox}
\usepackage{siunitx}

\newcommand{\be}{\begin{equation}}
\newcommand{\ee}{\end{equation}}
\newcommand{\benn}{\begin{equation*}}
\newcommand{\eenn}{\end{equation*}}

\newtheorem{remark}{Remark}[section]

\DeclareMathOperator*{\argmin}{arg\,min}
\DeclareMathOperator*{\argmax}{arg\,max}

\def\R {\mathbb{R}}

\def\Rn {\R^n}

\def\cinh{\gamma}

\def\ind {I}  
\def\chara {\delta}  

\title{Primal-dual hybrid gradient algorithms for computing time-implicit Hamilton-Jacobi equations\thanks{T. Meng, S. Liu, and S. J. Osher are partially funded by AFOSR MURI FA9550-18-502, ONR N00014-18-1-2527 and N00014-20-1-2787.
W. Li’s work is partially supported by AFOSR YIP award No. FA9550-23-1-008, NSF DMS-2245097, and NSF RTG: 2038080.}}

\author{Tingwei Meng\footnotemark[2]\and Wenbo Hao\footnotemark[2]\ \footnotemark[4]\and  Siting Liu\footnotemark[2]\ \footnotemark[4]\and  Stanley J. Osher\footnotemark[2]\ \and Wuchen Li\footnotemark[3] \footnotemark[5]}

\begin{document}
\maketitle
\renewcommand{\thefootnote}{\fnsymbol{footnote}}
\footnotetext[2]{Department of Mathematics, UCLA, Los Angeles, CA 90025, USA.}
\footnotetext[3]{Department of Mathematics, University of South Carolina, Columbia, SC 29208, USA.}
\footnotetext[4]{W.H. and S.L. contributed equally to this work.}
\footnotetext[5]{Corresponding author (wuchen@mailbox.sc.edu).}
\renewcommand{\thefootnote}{\arabic{footnote}}

\begin{abstract}
Hamilton-Jacobi (HJ) partial differential equations (PDEs) have diverse applications spanning physics, optimal control, game theory, and imaging sciences. This research introduces a first-order optimization-based technique for HJ PDEs, which formulates the time-implicit update of HJ PDEs as saddle point problems. We remark that the saddle point formulation for HJ equations is aligned with the primal-dual formulation of optimal transport and potential mean-field games (MFGs). This connection enables us to extend MFG techniques and design numerical schemes for solving HJ PDEs. We employ the primal-dual hybrid gradient (PDHG) method to solve the saddle point problems, benefiting from the simple structures that enable fast computations in updates.
Remarkably, the method caters to a broader range of Hamiltonians, encompassing non-smooth and spatiotemporally dependent cases. The approach's effectiveness is verified through various numerical examples in both one-dimensional and two-dimensional examples, such as quadratic and $L^1$ Hamiltonians with spatial and time dependence.
\end{abstract}

\section{Introduction}
Hamilton-Jacobi (HJ) partial differential equations (PDEs) find applications in various fields such as physics \cite{Arnold1989Math, Caratheodory1965CalculusI,Caratheodory1967CalculusII,Courant1989Methods,landau1978course},
optimal control \cite{Bardi1997Optimal, Elliott1987Viscosity,fleming1976deterministic,fleming2006controlled,mceneaney2006max}, game theory \cite{BARRON1984213, Buckdahn2011Recent, Evans1984Differential, Ishii1988Representation}, imaging sciences \cite{darbon2015convex,darbon2019decomposition,Darbon2016Algorithms,darbon2022hamilton}, and machine learning~\cite{chaudhari2018deep}.
In existing literature, a variety of approaches have been explored to numerically address HJ PDEs. In lower dimensions, high-order grid-based techniques such as essentially nonoscillatory schemes (ENO)~\cite{Osher1991High}, weighted ENO scheme~\cite{Jiang2000Weighted}, and discontinuous Galerkin method~\cite{Hu1999Discontinuous} are commonly employed, while in higher dimensions, strategies have been proposed to manage the challenges arising from the curse of dimensionality.
These works include, but are not limited to, max-plus algebra methods \cite{mceneaney2006max,akian2006max,akian2008max, dower2015max,Fleming2000Max,gaubert2011curse,McEneaney2007COD,mceneaney2008curse,mceneaney2009convergence}, dynamic programming and reinforcement learning \cite{alla2019efficient,bertsekas2019reinforcement}, tensor decomposition techniques \cite{dolgov2019tensor,horowitz2014linear,todorov2009efficient}, sparse grids \cite{bokanowski2013adaptive,garcke2017suboptimal,kang2017mitigating}, model order reduction \cite{alla2017error,kunisch2004hjb}, polynomial approximation \cite{kalise2019robust,kalise2018polynomial}, optimization methods \cite{darbon2015convex,darbon2019decomposition,Darbon2016Algorithms,yegorov2017perspectives,chow2019algorithm,chen2021hopf,chen2021lax} and neural networks \cite{darbon2019overcoming,bachouch2018deep, Djeridane2006Neural,jiang2016using, Han2018Solving, hure2018deep, hure2019some, lambrianides2019new, Niarchos2006Neural, reisinger2019rectified,royo2016recursive, Sirignano2018DGM,darbon2021some,darbon2023neural}. 

In this study, we introduce an innovative optimization-based methodology for tackling time-implicit computation of HJ PDEs. Our approach formulates the HJ PDE as a min-max problem by introducing a Lagrange multiplier. We then employ the primal-dual hybrid gradient (PDHG) method \cite{chambolle2011first} for finding the saddle point, which corresponds to the HJ PDE solution. 
In particular, the saddle point formulation of the HJ equations connects to the primal-dual formulation of optimal transport and potential mean-field games. Mean-field games (MFGs), introduced in \cite{lasry2007mean, huang2006large} are a mathematical framework employed for modeling and analyzing the equilibrium state of strategic interactions within a large population. This framework has found widespread application in various domains~\cite{lachapelle2011mean,gomes2014mean,Gomes2015,carmona2015mean,cardaliaguet2018mean,kizilkale2019integral,han2019mean,lee2021controlling}. The MFG system can be described using a system of coupled PDEs: a Fokker-Planck equation evolves forward in time, and an HJ equation evolves backward in time.
In this study, we leverage this relationship to address HJ PDEs. Using this connection, techniques developed for solving MFGs can be seamlessly extended to tackle the computation of the solution to HJ PDEs with initial conditions. 
For instance, PDHG method has been applied to numerically solve potential MFGs that can be cast into a saddle point formulation~\cite{papadakis2014optimal,briceno2018proximal, briceno2019implementation}. We adapt such methodology to our proposed saddle point problem, which in turn enables us to effectively address the original HJ PDEs.
When compared to other grid-based methods, our current approach, while possessing a first-order accuracy level, attains numerical unconditional stability through the utilization of implicit time discretization. This feature enables us to adopt larger time steps.
Compared to other optimization-based methods, our technique boasts the capability to handle a broader range of Hamiltonian functions, including those that exhibit non-smooth behavior and dependence on both $x$ and $t$. Furthermore, our algorithm benefits from its straightforward saddle point formulation, which allows simple updates in each iteration. 
We emphasize that the updates in our method do not involve nonlinear inversion. The only non-trivial update is managed through the proximal point operator of the Hamiltonian $H$, which can be computed in parallel.

The method of transforming a PDE problem into a saddle point problem and subsequently utilizing PDHG for its solution is used in solving reaction-diffusion equations~\cite{liu2023firstorder,fu2023high,carrillo2023structure} as well as conservation laws~\cite{liu2023primal}. 
These works use PDHG together with integration by parts to solve initial value problems with or without constraints, obtaining simple implicit in time updates. 
It is well-known in the literature that one-dimensional HJ PDEs and conservation laws are equivalent to each other. 
However, applying the method used for solving conservation laws in~\cite{liu2023primal} directly to solve HJ PDEs is not straightforward. In~\cite{liu2023primal}, the key step for obtaining a saddle point problem involves employing integration by parts.
The method of integration by parts is applicable to conservation laws, as these equations entail gradients of flux functions, a feature that is absent in HJ PDEs.
Nevertheless, in this study, for convex HJ PDEs, we instead use duality via the Fenchel-Legendre transform in lieu of integration by parts, effectively avoiding nonlinear updates. The merit of this approach lies in the simplicity of the saddle point formulation. This simplicity facilitates updates within our method to have either explicit formulations or be conducive to parallel computation. 
It is pertinent to note that while the formulation in this work may resemble the adjoint method~\cite{evans2010adjoint}, the proposed numerical technique and saddle point formulation differ substantially.

We show several numerical examples in one dimension and two dimensions\footnote{Codes are available at \url{https://github.com/TingweiMeng/HJ_solver_PDHG}.}. These numerical examples show the ability of this method to handle certain Hamiltonians which may depend on $(x,t)$. In each iteration, the updates of the functions are independent on each point, which makes it possible to use parallel computing to accelerate the algorithm. Moreover, in the special case when the Hamiltonian $H(x,t,\cdot)$ is a separable and shifted $1$-homogeneous function for any $(x,t)$, the algorithm has a simpler form, and we obtain an explicit formula for the updates of the dual variables. 

The paper is organized as follows. In Section~\ref{sec:cont_formulation}, we show the saddle point problem related to the HJ PDE and the proposed algorithm in the function space. In Section~\ref{sec:1d}, we focus on the one-dimensional case and show both semi-discrete and fully-discrete formulations of the algorithm. In Section~\ref{sec:2d}, we show the algorithm in two-dimensional case. In Section~\ref{sec:numerics}, we provide several numerical results which demonstrate the ability of the algorithm to handle certain Hamiltonians which may depend on $(x,t)$. In Section~\ref{sec:summary}, we show the summary and future work. More details about the proposed method and different saddle point formulations with the corresponding algorithms are also shown in the appendix.

\section{The saddle point formulation of HJ PDEs}\label{sec:cont_formulation}
In this section, we give details on the derivation of the saddle point formulation of the HJ PDEs. Furthermore, we devise a primal-dual hybrid gradient algorithm that solves it.
\subsection{The saddle point formulation}
In this paper, we solve the following HJ PDE in the domain $\Omega \times [0,+\infty)$
\begin{equation}\label{eqt:cont_HJ_eqt}
\begin{dcases}
\frac{\partial \phi(x,t)}{\partial t} + H(x, t, \nabla_x \phi(x,t)) = \epsilon\Delta_x \phi(x,t), & x\in \Omega, t\in [0,T],\\
\phi(x,0) = g(x), & x\in \Omega,
\end{dcases}
\end{equation}
where $\Omega = \prod_{i=1}^n [a_i, b_i]$ is the spatial domain in $\R^n$ with periodic boundary condition, and $\epsilon\geq 0$ is the diffusion parameter. Note that when $\epsilon = 0$, equation~\eqref{eqt:cont_HJ_eqt} is the first-order HJ PDE. In general, we assume that the Hamiltonian $H(x,t,p)$ is convex with respect to $p$ for any $x\in \Omega$, $t\in [0,T]$.

To address these equations, we treat them as constraints within the following optimization problem
\begin{equation}
\min_{\phi \text{ satisfying~\eqref{eqt:cont_HJ_eqt}}} -c\int_{\Omega} \phi(x,T) dx,
\end{equation}
where $c>0$ is a hyperparameter. In our numerical experiments, we observe that the constant $c$ can influence the convergence of our proposed algorithm. Moreover, it is possible to formulate more intricate objective functions in this optimization problem. The question of selecting a suitable objective function may pose an interesting problem that requires further investigation.

Subsequently, we introduce the Lagrange multiplier $\rho$, leading to the following computation
\begin{equation}\label{eqt:saddle_cont_appendix}
{\scriptsize
\begin{split}
&\min_{\phi \text{ satisfying~\eqref{eqt:cont_HJ_eqt}}}  -c\int_{\Omega} \phi(x,T) dx\\
=& \min_{\substack{\phi\\ \phi(x,0)=g(x)}}\max_{\rho}  \int_0^T \int_\Omega \rho(x,t)\left(\frac{\partial \phi(x,t)}{\partial t} - \epsilon \Delta_x \phi(x,t)  + H\left(x, t,\nabla_x \phi(x,t)\right) \right)dxdt - c\int_\Omega \phi(x,T)dx\\
=& \min_{\substack{\phi\\ \phi(x,0)=g(x)}}\max_{\rho}  \int_0^T \int_\Omega \rho(x,t)\left(\frac{\partial \phi(x,t)}{\partial t}  - \epsilon \Delta_x \phi(x,t) +\max_{v\in\Rn} \left\{\langle v,\nabla_x\phi(x,t)\rangle - H^*\left(x, t,v\right)\right\} \right)dxdt - c\int_\Omega \phi(x,T)dx.
\end{split}
}
\end{equation}

If $\rho$ is non-negative, we can move the maximization with respect to $v$ outside of the integral and formulate the following saddle point problem.

\begin{equation} \label{eqt:saddle_cont}
\begin{split}
    \min_{\substack{\phi\\ \phi(x,0)=g(x)}}\max_{\substack{\rho,v\\ \rho \geq 0}} \, \mathcal{L}(\phi, \rho, v),
\end{split}
\end{equation}
where 
\begin{equation}\label{eqt:saddle_pt_objfn}
{\scriptsize
\begin{split}
\mathcal{L}(\phi,\rho,v) = \int_0^T \int_\Omega \rho(x,t)\left(\frac{\partial \phi(x,t)}{\partial t} + \langle v(x,t), \nabla_x\phi(x,t)\rangle - \epsilon \Delta_x \phi(x,t)  - H^*\left(x, t,v(x,t)\right) \right)dxdt - \int_\Omega c\phi(x,T)dx.
\end{split}}
\end{equation}
Here, $H^*(x,t,\cdot)$ is the Fenchel-Legendre transform of $H(x,t,\cdot)$ for any $x\in\Omega$ and $t\in [0,T]$, defined by $H^*(x,t,v) = \sup_{p\in\Rn}\{\langle v,p\rangle - H(x,t,p)\}$. For more discussion about this saddle point problem, see Appendix~\ref{appendix:cont}.

\begin{remark} \label{rem:relation_mfc}
The first-order optimality condition of~\eqref{eqt:saddle_cont} gives
\begin{equation}\label{eqt:2pde_cont}
\begin{dcases}
\frac{\partial \phi(x,t)}{\partial t} + H(x, t,\nabla_x \phi(x,t)) \leq \epsilon\Delta_x \phi(x,t), & x\in \Omega, t\in [0,T],\\
\frac{\partial \rho(x,t)}{\partial t} + \nabla_x \cdot (\rho(x,t)v(x,t)) + \epsilon\Delta_x \rho(x,t)=0, & x\in \Omega, t\in [0,T],\\
v(x,t) = \nabla_p H(x, t,\nabla_x \phi(x,t)),& x\in \Omega, t\in [0,T],\\
\phi(x,0) = g(x),\quad \rho(x,T) = c, & x\in \Omega,\\
\rho(x,t)\left(\frac{\partial \phi(x,t)}{\partial t} + H(x, t,\nabla_x \phi(x,t)) - \epsilon\Delta_x \phi(x,t)\right) = 0, & x\in \Omega, t\in [0,T].
\end{dcases}
\end{equation}
This coupled system of PDEs differs slightly from the HJ PDE~\eqref{eqt:cont_HJ_eqt}, with the distinction that the first row entails an inequality rather than an equality. As indicated by the final row, if $\rho(x,t)\neq 0$ for all $x,t$, the initial inequality in the first row transforms into an equality, thus yielding a solution to the HJ PDE~\eqref{eqt:cont_HJ_eqt}. 
The discrepancy between the first row and the HJ PDE~\eqref{eqt:cont_HJ_eqt} is attributed to the imposition of a constraint $\rho \geq 0$ in the saddle point problem~\eqref{eqt:saddle_cont}. For more details, see Appendix~\ref{appendix:cont}. 

The structure of this coupled system of PDEs closely resembles that of the coupled PDEs encountered in mean-field control problems. Essentially, when the density solution $\rho$ associated with the corresponding mean-field control problem maintains a positive value throughout the entire domain, the PDE system~\eqref{eqt:2pde_cont} gives a solution to the HJ PDE~\eqref{eqt:cont_HJ_eqt}. This encompasses scenarios where $\epsilon > 0$, as the diffusion term introduces the Brownian motion component within the underlying stochastic process, ensuring a nonzero density function $\rho$. In the case where $\epsilon$ is zero, although we cannot establish a theoretical assurance, we do provide numerical validation in Section~\ref{sec:numerics}.
\end{remark}

\subsection{PDHG algorithm}
Within existing literature, a widely recognized approach for solving saddle point problems is PDHG algorithm~\cite{chambolle2011first}. In this section, we provide a brief review of this method. It solves saddle point problems in the form of
\begin{align*}
    \min_{x\in X}\max_{y\in Y} \; \langle Kx, y\rangle + g(x) - f(y),
\end{align*} where $f$ and $g$ are convex functions, and $K$ represents a linear operator. 
This algorithm is an iterative method, where in each iteration, the primal variable $x$ and the dual variable $y$ are updated separately using the proximal point operators of $g$ and $f$.
With given suitable stepsizes $\sigma,\tau>0$, the $m$-th iteration update can be written as follows:
\begin{align*}
 \begin{dcases}
   x^m& = \argmin_{x\in X}\;  \langle x, K^T \bar{y}^{m-1} \rangle + g(x) + \dfrac{1}{2 \sigma} \| x-x^{m-1}\|^2_{X},\\
       y^m & = \argmax_{y\in Y}\;  \langle Kx^m, y \rangle  - f(y) - \dfrac{1}{2 \tau} \| y-y^{m-1}\|^2_{Y},\\
       \bar{y}^{m} & = 2 y^m - y^{m-1}.
\end{dcases}
\end{align*}
The above calculation requires proximal gradient descent (ascent) steps of the primal (dual) variables.
The computation for this method has also been expanded in scope by~\cite{valkonen2014primal} to encompass more generalized problems where the operator $K$ can be nonlinear. 

\subsection{PDHG algorithm for solving HJ PDEs}
In this section, we apply PDHG algorithm to solve the saddle point problem~\eqref{eqt:saddle_cont}.
Owing to the simplicity of~\eqref{eqt:saddle_pt_objfn}, the primal and dual updates either possess explicit formulas or can be computed in parallel using the proximal point operator of the function $H(x,t,\cdot)$.
This makes PDHG well-suited for solving the saddle point problem presented in~\eqref{eqt:saddle_cont}. 
Consequently, we employ this method, incorporating pre-conditioning on $\phi$, to solve this specific saddle point problem. The outlined algorithm is presented in Algorithm~\ref{alg:pdhg_cont}.
Further details about pre-conditioning and other techniques are summarized in Remark~\ref{rem:details_PDHG}. It is important to note that, although we have a simple formulation and updating rule in our method, in most cases, an explicit formula for the joint update of $\rho$ and $v$ is unavailable. As a consequence, each iteration in our algorithm involves multiple coordinate updating steps for $\rho$ and $v$ independently.

\begin{algorithm}[htbp]
\SetAlgoLined
\SetKwInOut{Input}{Inputs}
\SetKwInOut{Output}{Outputs}
\Input{Stepsize $\tau, \sigma>0$, error tolerance $\delta>0$, inner maximal iteration number $N_{inner}$ and outer maximal iteration number $N_{outer}$.}
\Output{the solution to the corresponding HJ PDE~\eqref{eqt:cont_HJ_eqt}.}
Initialize the functions by $\phi^0(x,t)=g(x)$ for all $x\in\Omega$ and $t\in[0,T]$, $\rho^0 \equiv c$, $v^0\equiv 0$.

 \For{$\ell = 0,1,\dots,N_{outer}-1$}{
 Update the function $\phi\colon \Omega\times [0,T]\to \R$ using
 \begin{equation}\label{eqt:phi_update_pdhg_cont}
 {\small
 \begin{split}
 \phi^{\ell+1} = \argmin_{\phi : \phi(x,0)=g(x)} \mathcal{L}(\phi, \rho^\ell, v^\ell) + \frac{1}{2\tau} \|\nabla \phi - \nabla \phi^\ell\|^2 = \phi^\ell + \tau (-\Delta)^{-1}(\partial_t\rho^\ell + \nabla_x\cdot (v^\ell\rho^\ell) + \epsilon \Delta_{x} \rho^\ell),
 \end{split}
 }
 \end{equation}
 where $\nabla$, $\Delta$ are differential operators with respect to both $x$ and $t$, and $(-\Delta)^{-1} f$ means the solution to $-\Delta u=f$ in $\Omega\times [0,T]$ with periodic spatial condition, Dirichlet initial condition $u(x,0)=0$, and Neumann terminal condition $\partial_t u(x,T)=0$.
 
 \If{$\|\partial_t \phi^{\ell+1}(x,t) + H(x, t,\nabla_x \phi^{\ell+1}(x,t)) - \epsilon\Delta_x \phi^{\ell+1}(x,t)\|_1 \leq \delta$}{
   Return $\phi^{\ell+1}$.
   }
 
 Set $\bar\phi^{\ell+1} = 2\phi^{\ell+1} - \phi^{\ell}$.
 
 Set $v^{\ell+1,0} = v^{\ell}$, $\rho^{\ell+1,0} = \rho^{\ell}$.

  \For{$m=0,1,\dots, N_{inner}-1$}{
 Update $v$ by
 \begin{equation}\label{eqt:v_update_pdhg_cont}
 {\small
 \begin{split}
     & v^{\ell+1,m+1} = \argmax_{v} \mathcal{L}(\bar\phi^{\ell+1}, \rho^{\ell+1,m}, v) - \frac{1}{2\sigma} \|\rho^{\ell+1,m} (v - v^{\ell+1,m})\|^2 \\
    &= \argmin_{v} \int_0^T\int_\Omega H^*(x,t,v(x,t)) + \frac{\rho^{\ell+1,m}(x,t)}{2\sigma} \left\|v(x,t) - v^{\ell+1,m}(x,t) - \sigma \frac{\nabla_x \bar\phi^{\ell+1}(x,t)}{\rho^{\ell+1,m}(x,t)} \right\|^2 dxdt.
    \end{split}
    }
 \end{equation}
 
 Update $\rho$ by 
 \begin{equation}\label{eqt:rho_update_pdhg_cont}
 {\small
 \begin{split}
     \rho^{\ell+1,m+1} &= \argmax_{\rho: \rho \geq 0} \mathcal{L}(\bar\phi^{\ell+1}, \rho, v^{\ell+1,m+1}) - \frac{1}{2\sigma} \|\rho - \rho^{\ell+1,m}\|^2
    = \max\{\mu^{\ell+1,m+1}, 0\},
    \end{split}
    }
 \end{equation}
 where $\mu^{\ell+1,m+1}$ is a function defined by $\mu^{\ell+1,m+1}(x,t) = \rho^{\ell+1,m}(x,t) + \sigma (\partial_t \bar\phi^{\ell+1}(x,t) +\langle v^{\ell+1,m+1}(x,t), \nabla_x\bar\phi^{\ell+1}(x,t)\rangle - H^*(x, t,v^{\ell+1,m+1}(x,t)) - \epsilon \Delta_x \bar\phi^{\ell+1}(x,t))$.
 }
 Set $v^{\ell+1} = v^{\ell+1,N_{inner}}$, $\rho^{\ell+1} = \rho^{\ell+1,N_{inner}}$.
 }
 Return $\phi^{N_{outer}}$.
 \caption{The proposed algorithm for solving~\eqref{eqt:saddle_cont}\label{alg:pdhg_cont}}
\end{algorithm}

\begin{remark} \label{rem:details_PDHG}
The value of the constant $c$ does not impact the solution of the HJ PDE~\eqref{eqt:cont_HJ_eqt}, yet we have observed its effects on the convergence during our numerical experiments.
While we opt for a constant value of $c$ in this paper, it's noteworthy that it could potentially take on any positive function. This choice serves as the terminal condition for $\rho$ in~\eqref{eqt:2pde_cont}. Further investigation is required to understand the influence of this terminal condition on the convergence of the proposed method and to determine the appropriate method for its selection.

There are several strategies to expedite the convergence of the algorithm. One approach is to modify the penalty terms in the updates of $\phi, \rho, v$. Specifically, for $\phi$, selecting a penalty term such as $\|\nabla \phi - \nabla \phi^\ell\|^2$ can enhance computational efficiency (this technique is called pre-conditioning), as demonstrated in previous studies~\cite{jacobs2019solving,liu2023primal}. In our exploration, we also investigated different penalty expressions for $\rho$ and $v$. Notably, a quadratic penalty demonstrated its effectiveness for $\rho$, while regarding $v$, we observed comparable performance between penalty terms such as $\|\rho^{\ell+1,m} (v - v^{\ell+1,m})\|^2$ and $\rho^{\ell+1,m}\| v - v^{\ell+1,m}\|^2$. These alternatives exhibited more favorable results compared to the quadratic penalty. In the subsequent sections of this paper, unless otherwise explicitly stated, we employ the notation $\|f\|$ to signify the $L^2$ norm when $f$ represents a function, and to indicate the $\ell^2$ norm if $f$ is a finite-dimensional vector.

Furthermore, we adopted a time interval partitioning strategy, employing the proposed algorithm separately within each subdivided interval. Although this necessitates solving multiple saddle point problems, the resulting smaller dimensions of each problem substantially expedite the algorithm's runtime, as detailed in prior research (Algorithm 3 in~\cite{liu2023primal}).
\end{remark}

\section{One-dimensional HJ PDEs}\label{sec:1d}
In this section, we focus on one-dimensional HJ PDEs, where the spatial domain $\Omega$ is defined as the interval $[a,b]$. We begin by introducing the semi-discrete method in Section~\ref{sec:semidisc_1d}, and then proceed to explain the fully-discrete method in Section~\ref{sec:fulldisc_1d}. 
Throughout the rest of this paper, we employ the notation $(a_i)_i$, $(a_{i,j})_{i,j}$, and $(a_{i,j,k})_{i,j,k}$ to represent a vector, matrix, and tensor, respectively, where the elements are denoted by $a_i$, $a_{i,j}$, and $a_{i,j,k}$.

\subsection{Semi-discrete formulation}\label{sec:semidisc_1d}

To solve the continuous HJ PDE~\eqref{eqt:cont_HJ_eqt}, our initial step involves discretizing the spatial domain $[a,b]$. We define $x_i$ as the $i$-th grid point on a uniformly spaced grid within the interval $[a,b]$. Specifically, $x_i$ is calculated as $a + \frac{(b-a)(i-1)}{n_x-1}$, where $n_x$ represents the total number of grid points. The semi-discrete approach for a general numerical Hamiltonian $\hat H$ is presented as follows:
\begin{equation}\label{eqt:semidisc_HJ_equation}
\begin{dcases}
\dot \phi_i(t) + \hat H\left(x_i, t, (D_x^+\phi)_i(t), (D_x^-\phi)_i(t)\right) = \epsilon (D_{xx} \phi)_i(t), & i=1,\dots, n_x,\, t\in [0,T],\\
\phi_i(0) = g(x_i), & i=1,\dots, n_x.
\end{dcases}
\end{equation}
According to the theory of first order monotone scheme, the numerical Hamiltonian $\hat H$ needs to be consistent (i.e., $\hat H(x,t,p,p) = H(x,t,p)$) and monotone (i.e., non-increasing with respect to $p^+$ and non-decreasing with respect to $p^-$). For more details, we refer readers to~\cite{shu2007high}. 
In this paper, we use $D_x^+$ and $D_x^-$ to represent the right and left finite difference approximations of the spatial derivative, respectively. Specifically, $(D_x^+\phi)_i(t)$ is computed as $\frac{\phi_{i+1}(t) - \phi_i(t)}{\Delta x}$, and $(D_x^-\phi)_i(t)$ is calculated as $\frac{\phi_{i}(t) - \phi_{i-1}(t)}{\Delta x}$. Additionally, $D_{xx}$ approximates the Laplace operator, given by $(D_{xx} \phi)_i(t) = \frac{\phi_{i+1}(t) - 2 \phi_{i}(t) + \phi_{i-1}(t)}{\Delta x^2}$.
Here, we assume $\hat H(x_i,t, p^+, p^-)$ is jointly convex with respect to $(p^+,p^-)$ for any $i=1,\dots, n_x$ and $t\in [0,T]$.

Similar to the continuous version, we incorporate these formulas into the constraints of an optimization problem and introduce the Lagrange multiplier $\rho$. Consequently, this semi-discrete equation can be solved using the following saddle point formulation
\begin{equation}\label{eqt:saddle_pt_general_semidisc}
{\small
\begin{split}
    \min_{\substack{\phi_i \forall i\\ \phi_i(0)=g(x_i)}}\max_{\substack{\rho_i, v_i^+, v_i^-\forall i\\ \rho_i\geq 0}}  \int_0^T \sum_{i=1}^{n_x}\rho_i(t) \Bigg(\dot \phi_i(t) + v^+_i(t)(D_x^+\phi)_i(t) + v^-_i(t)(D_x^-\phi)_i(t) 
    - \epsilon (D_{xx} \phi)_i(t) \\
    -\hat H^*(x_i,t, v_i^+(t), v_i^-(t)) \Bigg)dt - c\sum_{i=1}^{n_x}\phi_i(T),
\end{split}
}
\end{equation}
where $c>0$ is a hyper-parameter, and $\hat H^*(x_i, t, \cdot,\cdot)$ is the Fenchel-Legendre transform of $\hat H(x_i, t, \cdot, \cdot)$. The derivation of this saddle point formula is similar to the continuous case and is therefore omitted here. For more details, please refer to Appendix~\ref{sec:appendix_semidisc_1d}.

To solve this saddle point problem, we apply the PDHG method. Generally, deriving a direct updating formula for the combined variables $(\rho_i, v_i^+, v_i^-)$ is challenging.
As a result, an inner loop is utilized in which each iteration involves the sequential updates of $(v_i^+, v_i^-)$ and $\rho_i$.
Let's denote the objective function in~\eqref{eqt:saddle_pt_general_semidisc} as $\mathcal{L}_{semi}((\phi_i)_i, (\rho_i)_i, (v^+_i)_i, (v^-_i)_i)$. The proposed algorithm is summarized in Algorithm~\ref{alg:pdhg_semi_1d}.

\begin{remark}
In practice, many numerical Hamiltonians $\hat H$ are separable, i.e., $\hat H$ satisfies $\hat H(x, t, p^+,p^-) = \hat H_1(x, t, p^+) + \hat H_2(x, t, p^-)$, where $\hat H_1$ is non-increasing in $p^+$, and $\hat H_2$ is non-decreasing in $p^-$. In this scenario, the Fenchel-Legendre transform obeys the relationship $\hat H^*(x, t, v^+,v^-) = \hat H_1^*(x,t,v^+) + \hat H_2^*(x,t,v^-)$. As a result, the updates of $v^+$ and $v^-$ can be computed separately and in parallel.
\end{remark}

\subsection{Fully-discrete formulation}\label{sec:fulldisc_1d}
We proceed by introducing a fully-discrete method through time discretization. In this paper, we advocate implicit time discretization to enable the use of larger time steps and circumvent the restrictive Courant–Friedrichs–Lewy (CFL) time step restriction for explicit methods. The time interval $[0,T]$ is discretized into $n_t-1$ equi-spatial subintervals, each of length $\Delta t = \frac{T}{n_t-1}$. Given a function $f$ defined on $[a,b]\times [0,T]$, we denote $f_{i,k}$ as the function value at the grid point $(x_i, t_k)$.

We adopt the backward Euler scheme $(D_t^- \phi)_{i,k} = \frac{\phi_{i, k} - \phi_{i, k-1}}{\Delta t}$ to approximate the time derivative at $(x_i, t_k)$. Employing this discretization, the resulting numerical scheme is as follows:
\begin{equation}\label{eqt:disc_HJ_equation}
\begin{dcases}
(D_t^- \phi)_{i,k} + \hat H\left(x_i, t_k,(D_x^+\phi)_{i,k}, (D_x^-\phi)_{i,k}\right) = \epsilon (D_{xx} \phi)_{i,k}, & i=1,\dots, n_x,\, k=2,\dots, n_t,\\
\phi_{i,1} = g(x_i), & i=1,\dots, n_x,
\end{dcases}
\end{equation}
and the saddle point problem becomes
\begin{equation}\label{eqt:saddle_pt_general_fulldisc}
{\small
\begin{split}
\min_{\substack{\phi_{i,k}\forall i,k\\ \phi_{i,1}=g(x_i)}}\max_{\substack{\rho_{i,k},v^+_{i,k}, v^-_{i,k}\forall i,k\\ \rho_{i,k} \geq 0}}  \sum_{i=1}^{n_x}\sum_{k=1}^{n_t-1} \rho_{i,k}\Bigg((D_t^- \phi)_{i,k+1}+ v^+_{i,k}(D_x^+\phi)_{i,k+1} +v^-_{i,k}(D_x^-\phi)_{i,k+1} - \epsilon (D_{xx} \phi)_{i,k+1} \\
- \hat H^*\left(x_i, t_{k+1},v^+_{i,k}, v^-_{i,k}\right) \Bigg) - \frac{c}{\Delta t}\sum_{i=1}^{n_x}\phi_{i,n_t}.
\end{split}
}
\end{equation}

The set of equations in~\eqref{eqt:disc_HJ_equation} constitutes a total of $n_x\times (n_t-1)$ equations. As a result, the dual variables $\rho_{i,k}$, $v^+_{i,k}$, and $v^-_{i,k}$ possess indices ranging from $i=1$ to $n_x$, and $k$ ranges from $1$ to $n_t-1$.
We apply PDHG to solve this saddle point problem, and the details are summarized in Algorithm~\ref{alg:pdhg_full_1d} and Appendix~\ref{appendix:fulldisc_1d}.

The objective function in~\eqref{eqt:saddle_pt_general_fulldisc} is linear when considering either $\phi$ or $\rho$. This linearity enables us to have explicit formulas for updating $\phi$ and $\rho$. When iteratively updating $\phi$, we effectively address a discrete Poisson's equation within the temporal-spatial domain through the utilization of the Fourier transform, thus facilitating efficient computation. The only non-linear part in~\eqref{eqt:saddle_pt_general_fulldisc} pertains to updating $(v^+, v^-)$, involving solving the proximal point of $(v^+, v^-)\mapsto \sum_{i=1}^{n_x}\sum_{k=1}^{n_t-1} \hat H^*(x_i,t_{k+1},v_{i,k}^+, v_{i,k}^-)$. The element at the $(i,k)$ position of this proximal point corresponds to the proximal point of $\hat H^*(x_i,t_{k+1},\cdot, \cdot)$, facilitating parallel updates for $(v^+,v^-)$. Furthermore, when $\hat H$ takes the form $\hat H(x,t, p^+, p^-) = \hat H_1(x,t,p^+) + \hat H_2(x,t,p^-)$, the dual function simplifies to $\hat H^*(x,t,v^+,v^-) = \hat H_1^*(x,t,v^+) + \hat H_2^*(x,t,v^-)$, enabling further parallelization for $v^+$ and $v^-$.

In specific cases where the Hamiltonian $H$ has particular structures, it becomes possible to update the variables $(\rho, v^+, v^-)$ simultaneously, removing the need for the inner loop. An example of this occurs when $H(x, t, p)$ is convex, separable, and $1$-homogeneous with respect to $p$. For more detailed information, please refer to Appendix~\ref{appendix:L1_1d}.

\section{Two-dimensional HJ PDEs} \label{sec:2d}
In this section, we address the two-dimensional HJ PDEs and present the semi-discrete approach in Section~\ref{sec:semidisc_2d} and the fully-discrete approach in Section~\ref{sec:fulldisc_2d}. Due to the similarities between the two-dimensional and one-dimensional cases, certain details and explanations have been omitted.

\subsection{Semi-discrete formulation}\label{sec:semidisc_2d}

We apply discretization to the spatial domain $[a_1,b_1]\times [a_2,b_2]$ using $n_x$ grid points in the first dimension and $n_y$ grid points in the second dimension. The grid sizes in these dimensions are represented by $\Delta x = \frac{b_1-a_1}{n_x}$ and $\Delta y = \frac{b_2-a_2}{n_y}$. We denote $x_i$ as the $i$-th grid point in the first dimension and $y_j$ as the $j$-th grid point in the second dimension.
The semi-discrete formulation for a general numerical Hamiltonian is as follows:
\begin{equation}\label{eqt:semidisc_HJ_equation_2d}
\begin{dcases}
\dot \phi_{i,j}(t) + \hat H(x_{i},y_j, t,(D_x^+\phi)_{i,j}(t), (D_x^-\phi)_{i,j}(t), (D_y^+\phi)_{i,j}(t), (D_y^-\phi)_{i,j}(t)) \\
\quad\quad = \epsilon (D_{xx}\phi)_{i,j}(t) + \epsilon(D_{yy}\phi)_{i,j}(t),
\quad\quad\quad\quad\,
 i=1,\dots, n_x; \, j =1,\dots, n_y,\, t\in [0,T],\\
\phi_{i,j}(0) = g(x_i,y_j), \quad\quad\quad\quad\quad\quad\quad\quad\quad\quad\quad\quad  i=1,\dots, n_x; \, j =1,\dots, n_y.
\end{dcases}
\end{equation}
where $D_x^+$, $D_x^-$, $D_y^+$, $D_y^-$, $D_{xx}$, and $D_{yy}$ are finite difference operators defined by
\begin{align*}
    (D_x^+\phi)_{i,j} &= \frac{\phi_{i+1,j}(t) - \phi_{i,j}(t)}{\Delta x}, \quad
    (D_x^-\phi)_{i,j} = \frac{\phi_{i,j} (t)- \phi_{i-1,j}(t)}{\Delta x}, \\(D_y^+\phi)_{i,j} &= \frac{\phi_{i,j+1}(t) - \phi_{i,j}(t)}{\Delta y}, \quad
    (D_y^-\phi)_{i,j} = \frac{\phi_{i,j} (t)- \phi_{i,j-1}(t)}{\Delta y}, \\(D_{xx}\phi)_{i,j} &= \frac{\phi_{i+1,j}(t) - 2 \phi_{i,j}(t) + \phi_{i-1,j}(t)}{\Delta x^2}, \quad
    (D_{yy}\phi)_{i,j} = \frac{\phi_{i,j+1}(t) - 2 \phi_{i,j}(t) + \phi_{i,j-1}(t)}{\Delta y^2}.
\end{align*}
Analogous to the one-dimensional case discussed in Section~\ref{sec:semidisc_1d}, we require that the numerical Hamiltonian $\hat H$ exhibits both consistency, meaning that $\hat H(x,y,t,p_1,p_1,p_2,p_2) = H(x,y,t,p_1,p_2)$, and monotonicity, meaning that $\hat H(x,y,t,p_1^+,p_1^-,p_2^+,p_2^-)$ is non-increasing with respect to $p_1^+$ and $p_2^+$, and non-decreasing with respect to $p_1^-$ and $p_2^-$.

This semi-discrete equation can be solved using the following saddle point formulation:
\begin{equation}\label{eqt:saddle_pt_general_semidisc_2d}
{\small
\begin{split}
    \min_{\substack{\phi_{i,j} \forall i,j\\ \phi_{i,j}(0)=g(x_i, y_j)}}\max_{\substack{\rho_{i,j}, v_{i,j}^+, v_{i,j}^-,\\ w_{i,j}^+, w_{i,j}^-:\rho_{i,j}\geq 0}}  \int_0^T \sum_{i=1}^{n_x}\sum_{j=1}^{n_y}\rho_{i,j}(t)\Big(\dot \phi_{i,j}(t) + v^+_{i,j}(t)(D_x^+\phi)_{i,j}(t) + v^-_{i,j}(t)(D_x^-\phi)_{i,j}(t) \\
    + w^+_{i,j}(t)(D_y^+\phi)_{i,j}(t) 
    + w^-_{i,j}(t)(D_y^-\phi)_{i,j}(t)
    - \epsilon (D_{xx} \phi)_{i,j}(t) - \epsilon (D_{yy} \phi)_{i,j}(t) \\
    -\hat H^*(x_i, y_j, t,v_{i,j}^+(t), v_{i,j}^-(t), w_{i,j}^+(t), w_{i,j}^-(t)) \Big)dt 
    -c\sum_{i=1}^{n_x}\sum_{j=1}^{n_y}\phi_{i,j}(T),
\end{split}
}
\end{equation}
where $c>0$ is a hyper-parameter, and $\hat H^*(x_i,y_j,t, \cdot,\cdot, \cdot,\cdot)$ is the Fenchel-Legendre transform of $\hat H(x_i,y_j, t,\cdot, \cdot, \cdot,\cdot)$.
For more properties of this saddle point problem and the corresponding algorithm, see Appendix~\ref{appendix:semidisc_2d} and Algorithm~\ref{alg:pdhg_semi_2d}.

\subsection{Fully-discrete formulation} \label{sec:fulldisc_2d}
In this section, we adopt implicit time discretization to derive a fully-discrete formulation, enabling the selection of larger time steps to bypass the CFL condition. We represent the number of grid points in the interval $[0,T]$ as $n_t$, and the spacing between consecutive grid points as $\Delta t = \frac{T}{n_t-1}$.
The fully-discrete HJ PDE with a numerical Hamiltonian $\hat H$ is given by 
\begin{equation}\label{eqt:disc_HJ_equation_2d}
\begin{dcases}
(D_t^- \phi)_{i,j,k} + \hat H\left(x_i, y_j,t_k, (D_x^+\phi)_{i,j,k}, (D_x^-\phi)_{i,j,k}, (D_y^+\phi)_{i,j,k}, (D_y^-\phi)_{i,j,k}\right) \\
\quad\quad\quad\quad= \epsilon (D_{xx} \phi)_{i,j,k} + \epsilon (D_{yy} \phi)_{i,j,k}, \quad\quad\quad\quad i=1,\dots, n_x; j=1,\dots, n_y; k=2,\dots, n_t,\\
\phi_{i,j,1} = g(x_i,y_j), \quad\quad\quad\quad \quad\quad\quad\quad\quad\quad\quad\quad\quad\,\,\,
i=1,\dots, n_x; j=1,\dots, n_y,
\end{dcases}
\end{equation}
and the saddle point problem becomes
\begin{equation}\label{eqt:saddle_pt_general_fulldisc_2d}
{\scriptsize
\begin{split}
\min_{\substack{\phi_{i,j,k}\forall i,j,k\\ \phi_{i,j,1}=g(x_i,y_j)}}\max_{\substack{\rho_{i,j,k}, v_{i,j,k}^+, v_{i,j,k}^-,\\ w_{i,j,k}^+, w_{i,j,k}^-:\rho_{i,j,k}\geq 0}}  \sum_{i=1}^{n_x}\sum_{j=1}^{n_y}\sum_{k=1}^{n_t-1} \rho_{i,j,k}\Bigg((D_t^- \phi)_{i,j,k+1}+ v^+_{i,j,k}(D_x^+\phi)_{i,j,k+1} +v^-_{i,j,k}(D_x^-\phi)_{i,j,k+1} + w^+_{i,j,k}(D_y^+\phi)_{i,j,k+1}   \\
+w^-_{i,j,k}(D_y^-\phi)_{i,j,k+1} - \epsilon (D_{xx} \phi)_{i,j,k+1} - \epsilon (D_{yy} \phi)_{i,j,k+1} - \hat H^*\left(x_i, y_j,t_{k+1},v^+_{i,j,k}, v^-_{i,j,k}, w^+_{i,j,k}, w^-_{i,j,k}\right) \Bigg) - \frac{c}{\Delta t}\sum_{i=1}^{n_x}\sum_{j=1}^{n_y}\phi_{i,j,n_t}.
\end{split}
}
\end{equation}
This saddle point problem is solved using PDHG method. For more details about this problem and the algorithm, we refer readers to Section~\ref{appendix:fulldisc_2d} and Algorithm~\ref{alg:pdhg_full_2d}.

Just like the one-dimensional scenario, the objective function's linearity concerning $\phi$ or $\rho$ leads to explicit update formulas for $\phi$ and $\rho$. On the other hand, updating $(v^+,v^-,w^+,w^-)$ necessitates solving the proximal point of $(v^+,v^-,w^+,w^-)\mapsto \sum_{i,j,k} \hat H^*(x_i, y_j,t_{k+1},v^+_{i,j,k}, v^-_{i,j,k}, w^+_{i,j,k}, w^-_{i,j,k})$, which can be conducted in parallel for each point $(x_i,y_j,t_{k+1})$. Furthermore, if $\hat H$ is separable, i.e., it can be expressed as $\hat H_1^1(x_i,y_j,t_{k+1}, v^+_{i,j,k}) + \hat H_2^1(x_i,y_j,t_{k+1}, v^-_{i,j,k}) + \hat H_1^2(x_i,y_j,t_{k+1}, w^+_{i,j,k}) + \hat H_2^2(x_i,y_j,t_{k+1}, w^-_{i,j,k})$, then updating $(v^+_{i,j,k},v^-_{i,j,k},w^+_{i,j,k},w^-_{i,j,k})$ can be further accomplished in parallel.

Analogous to the scenario in the one-dimensional case, if $H$ is convex, separable, and $1$-homogeneous with respect to $p$, there exists a specific formula for jointly updating $\rho, v^+, v^-, w^+, w^-$. This, in turn, eliminates the need for the inner loop. For more details, refer to Appendix~\ref{appendix:L1_2d}.

\section{Numerical results}\label{sec:numerics}
In this section, we display a range of numerical results that evaluate the performance of our proposed method. We initially use two simple experiments in Sections~\ref{sec:eg1} and~\ref{sec:eg2} to present error tables, which confirm that our method yields first-order accuracy in computing the viscosity solution for these examples. Subsequently, we utilize more intricate cases (examples 3 and 4 in Sections~\ref{sec:eg3} and~\ref{sec:eg4}) to demonstrate the benefits of using larger time steps. These experiments highlight the ability of our method to handle Hamiltonians that depend on $(x,t)$ and exhibit non-smooth behaviors.

Among these four examples, the second one involves Hamiltonians that are shifted $1$-homogeneous with respect to $p$, i.e., $H(x,t,\cdot) - H(x,t,0)$ is $1$-homogeneous for any $(x,t)$. For this case, we implement Algorithms~\ref{alg:pdhg_full_m_1d_L1} and~\ref{alg:pdhg_full_m_2d_L1} explained in Appendices~\ref{appendix:L1_1d} and~\ref{appendix:L1_2d} to eliminate the necessity of the inner loop in the proposed method. On the other hand, in the other three examples, we use Algorithms~\ref{alg:pdhg_full_1d} and~\ref{alg:pdhg_full_2d}.
In these experiments, for one-dimensional cases, we apply Engquist-Osher scheme~\cite{osher1988fronts,engquist1980stable,engquist1981one} and set the numerical Hamiltonian $\hat H(x,t,p^+,p^-) = H_-(x,t,p^+) + H_+(x,t,p^-)$, where $H_+$ is non-decreasing with respect to $p^-$, $H_-$ is non-increasing with respect to $p^+$, and they satisfy $H_+ + H_- = H$. For two-dimensional cases, if the Hamiltonian $H$ can be written as $H(x,t,p_1,p_2) = H_1(x,t,p_1) + H_2(x,t,p_2)$ for some functions $H_1$ and $H_2$, we handle each dimension separately and set the numerical Hamiltonian $\hat H(x,t,p_1^+,p_1^-,p_2^+,p_2^-) = H_{1,-}(x,t,p_1^+) + H_{1,+}(x,t,p_1^-) + H_{2,-}(x,t,p_2^+) + H_{2,+}(x,t,p_2^-)$, where $H_{i,+}$ ($i=1,2$) is non-decreasing with respect to $p_i^-$, $H_{i,-}$ ($i=1,2$) is non-increasing with respect to $p_i^+$, and they satisfy $H_{i,+}+H_{i,-} = H_i$. 
Note that our approach also works for non-separable Hamiltonians. For these situations, we define the numerical Hamiltonians based on references that will be specified for each case.

\subsection{First experiment: quadratic Hamiltonian}\label{sec:eg1}
In the first experiment, we solve the following HJ PDE:
    \begin{equation}\label{eqt:eg1_burgers}
    \begin{dcases}
        \frac{\partial \phi(x,t)}{\partial t} + \frac{1}{2}\|\nabla_x \phi(x,t)\|^2 = 0, & x\in [0,2]^n, t\in [0,1],\\
        \phi(x,0) = \frac{1}{2}\|x-1\|^2, & x\in [0,2]^n.
    \end{dcases}
    \end{equation}
We apply Algorithms~\ref{alg:pdhg_full_1d} and~\ref{alg:pdhg_full_2d} to solve this problem. We use the Engquist-Osher scheme in the saddle point formulation. To be specific, we set the numerical Hamiltonian $\hat H(x,t,p^+,p^-) = H_-(p^+) + H_+(p^-)$ for the one-dimensional case and $\hat H(x,t,p_1^+,p_1^-, p_2^+,p_2^-) = H_-(p_1^+) + H_+(p_1^-) + H_-(p_2^+) + H_+(p_2^-)$ for the two-dimensional case, where $H_-(p^+) = \frac{1}{2}\min\{p^+,0\}^2$ and $H_+(p^-) = \frac{1}{2}\max\{p^-,0\}^2$.
We provide the error tables for both one-dimensional and two-dimensional scenarios in Table~\ref{tab:eg10_1d} and Table~\ref{tab:eg10_2d} respectively.

\begin{table} [h!]
    \footnotesize
    \centering
    \begin{tabular}{c|c|c|c|c}
    \hline
         $n_x\times n_t$ & $20\times 11$ & $40\times 21$ & $80\times 41$ & $160\times 81$ \\
\hline
Averaged absolute residual of HJ PDE & 9.99E-07 & 9.99E-07 & 9.89E-07 & 9.82E-07 \\
\hline
$\ell^1$ relative error & 5.81E-02 & 3.24E-02 & 1.68E-02 & 8.27E-03 \\
       \hline
    \end{tabular}
    \caption{Error table illustrating the performance of our proposed method for solving the one-dimensional HJ PDE~\eqref{eqt:eg1_burgers}}
    \label{tab:eg10_1d}
\end{table}

\begin{table} [h!]
    \footnotesize
    \centering
    \begin{tabular}{c|c|c|c|c}
    \hline
         $n_x\times n_y\times n_t$ & $20\times 20\times 11$ & $40\times 40\times 21$ & $80\times 80\times 41$ & $160\times 160\times 81$ \\
\hline
Averaged absolute residual of HJ PDE & 1.00E-06 & 1.00E-06 & 9.99E-07 & 1.00E-06 \\
\hline
$\ell^1$ relative error & 5.52E-02 & 3.00E-02 & 1.46E-02 & 6.07E-03 \\
       \hline
    \end{tabular}
    \caption{Error table illustrating the performance of our proposed method for solving the two-dimensional HJ PDE~\eqref{eqt:eg1_burgers}}
    \label{tab:eg10_2d}
\end{table}

\begin{figure}[htbp]
    \centering
    \begin{subfigure}{0.45\textwidth}
        \centering \includegraphics[width=\textwidth]{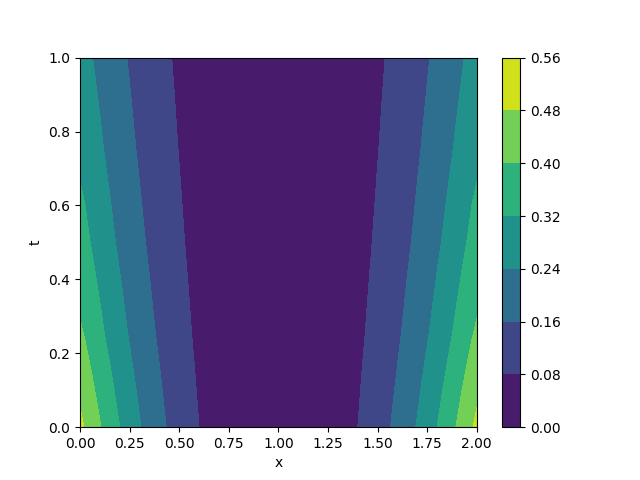}
        \caption{1D, $n_x = 80$, $n_t = 5$}
    \end{subfigure}
    \hfill
    \begin{subfigure}{0.45\textwidth}
        \centering \includegraphics[width=\textwidth]{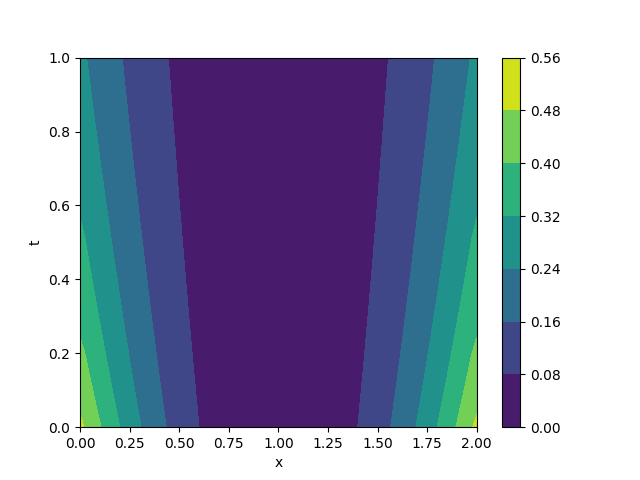}
        \caption{1D, $n_x = 80$, $n_t = 41$}
    \end{subfigure}

    \begin{subfigure}{0.45\textwidth}
        \centering \includegraphics[width=\textwidth]{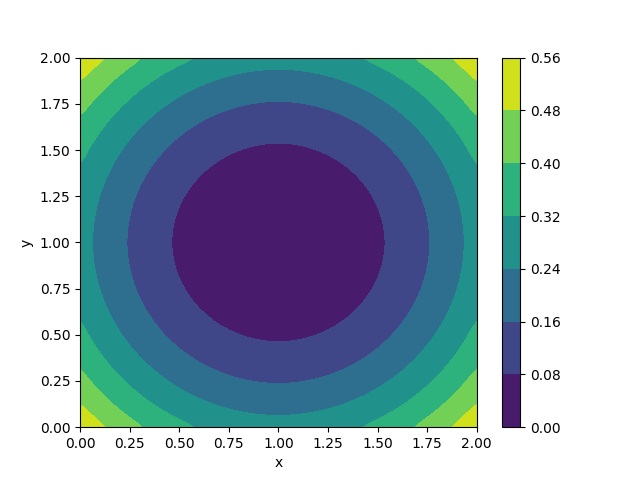}
        \caption{2D, $t=T$, $n_x = n_y= 80$, $n_t = 5$}
    \end{subfigure}
    \hfill
    \begin{subfigure}{0.45\textwidth}
        \centering \includegraphics[width=\textwidth]{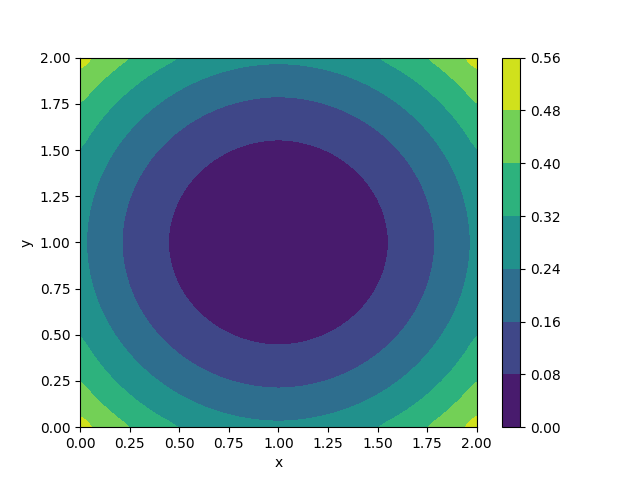}
        \caption{2D, $t=T$, $n_x=n_y = 80$, $n_t = 41$}
    \end{subfigure}
    
    \caption{
    Contours illustrating the solution of the one-dimensional HJ PDE~\eqref{eqt:eg1_burgers} in (a) and (b), along with the level sets for the two-dimensional HJ PDE in (c) and (d).
    In (a) and (c), a relatively larger time step of $\Delta t = 0.25$ is applied, whereas in (b) and (d), a smaller time step of $\Delta t = 0.025$ is employed.} 
    \label{fig:eg10_epsl0}
\end{figure}

Each table contains two rows of error measurements. In the first row, we showcase the average absolute value of the PDE residuals. This error is calculated using the following formula for one-dimensional cases:
\begin{equation*}
\frac{1}{n_x(n_t-1)}\sum_{i=1}^{n_x}\sum_{k=2}^{n_t}|(D_t^- \phi)_{i,k} + \hat H(x_i, t_k,(D_x^+ \phi)_{i,k}, (D_x^- \phi)_{i,k}) - \epsilon (D_{xx} \phi)_{i,k}|
\end{equation*}
and for two-dimensional cases:
\begin{equation*}
{\small
\begin{split}
\frac{1}{n_xn_y(n_t-1)}\sum_{i=1}^{n_x}\sum_{j=1}^{n_y}\sum_{k=2}^{n_t}|(D_t^- \phi)_{i,j,k} + \hat H(x_i,y_j, t_{k}, (D_x^+ \phi)_{i,j,k}, (D_x^- \phi)_{i,j,k}, (D_y^+ \phi)_{i,j,k}, (D_y^- \phi)_{i,j,k}) - \epsilon (D_{xx} \phi + D_{yy} \phi)_{i,j,k}|.
\end{split}
}
\end{equation*}
The residual errors observed across all cases are consistently below $10^{-6}$. Remarkably, we employ the residual error as the termination condition for the proposed method, setting the threshold $\delta$ to $10^{-6}$. These errors confirm the convergence of the our algorithm in terms of reaching an error below the specified threshold for all tested grid sizes.

In the second row, we conduct a comparison between the numerical solution obtained through our proposed method and the reference solution, which we term the ``ground truth". This reference solution, denoted as $\phi^{gt}$, is generated using either the Lax-Oleinik formula~\cite{bardi1984hopf,evans1998partial,Hopf1965} or the explicit Engquist-Osher scheme~\cite{osher1988fronts,engquist1980stable,engquist1981one} on a finely discretized grid. Consequently, the $\ell^1$ relative errors are presented. The calculation of the error is as follows:
\begin{itemize}
    \item
For one-dimensional cases:
\begin{equation*}
    \frac{\frac{1}{n_x n_t} \sum_{i=1}^{n_x}\sum_{k=1}^{n_t} |\phi_{i,k} - \phi_{i,k}^{gt}|}{\max\left\{\frac{1}{n_x n_t} \sum_{i=1}^{n_x}\sum_{k=1}^{n_t} |\phi_{i,k}^{gt}|, 1.0\right\}}
\end{equation*}

\item For two-dimensional cases:
\begin{equation*}
    \frac{\frac{1}{n_x n_y n_t} \sum_{i=1}^{n_x}\sum_{j=1}^{n_y}\sum_{k=1}^{n_t} |\phi_{i,j,k} - \phi_{i,j,k}^{gt}|}{\max\left\{\frac{1}{n_x n_y n_t} \sum_{i=1}^{n_x}\sum_{j=1}^{n_y}\sum_{k=1}^{n_t} |\phi_{i,j,k}^{gt}|, 1.0\right\}}
\end{equation*}
\end{itemize}
This evaluation provides insights into the accuracy of our method when compared to well-established reference solutions. 
It's noticeable that the error approximately halves when the grid size is doubled. This behavior demonstrates a first-order error reduction rate corresponding to the increase in grid size.

Furthermore, we present the one-dimensional solution in Fig~\ref{fig:eg10_epsl0} (a)-(b) and the two-dimensional solution in Fig~\ref{fig:eg10_epsl0} (c)-(d). In both this example and the subsequent examples, we visualize the one-dimensional solution through the representation of its level sets, where the $x$-axis corresponds to the spatial domain and the $y$-axis denotes the time domain. For the two-dimensional solution, we depict the level sets at a specific time instant (in this example, we select $t=1$), employing the $x$ and $y$ axes to denote the two dimensions of the spatial domain.

In figures (a) and (c), we utilize a larger time step of $\Delta t = \frac{T}{n_t-1} = 0.25$, while in figures (b) and (d), we opt for a smaller time step of $\Delta t = \frac{T}{n_t-1} = 0.025$. Throughout all cases, we maintain a spatial discretization count of $n_x = 80$ for the one-dimensional case and $n_x=n_y=80$ for the two-dimensional scenario. These visual representations reveal a satisfactory performance achieved by employing the larger time step discretization, indicating that the proposed method is not constrained by the CFL condition, thanks to the utilization of implicit time discretization.

\subsection{Second experiment: $L^1$ Hamiltonian}\label{sec:eg2}
In the second example, we solve the following HJ PDE in one-dimension and two-dimensions
\begin{equation}\label{eqt:eg2_L1}
    \begin{dcases}
        \frac{\partial \phi(x,t)}{\partial t} + \|\nabla_x \phi(x,t)\|_1 = 0, & x\in [0,2]^n, t\in[0,1],\\
        \phi(x,0) = \sum_{i=1}^n \sin \pi x_i, & x\in [0,2]^n.
    \end{dcases}
\end{equation}
The Hamiltonian in this scenario exhibits non-smooth properties. 
We solve this problem using Algorithms~\ref{alg:pdhg_full_m_1d_L1} and~\ref{alg:pdhg_full_m_2d_L1}.
Employing the same methodology as illustrated in the previous example, we calculate the error tables. Specifically, the one-dimensional error table is presented in Table~\ref{tab:eg0_1d}, while the two-dimensional error table is displayed in Table~\ref{tab:eg0_2d}.

\begin{table} [h!]
    \footnotesize
    \centering
    \begin{tabular}{c|c|c|c|c}
    \hline
         $n_x\times n_t$ & $20\times 11$ & $40\times 21$ & $80\times 41$ & $160\times 81$ \\
\hline
Averaged absolute residual of HJ PDE & 5.72E-07 & 7.64E-07 & 5.96E-07 & 6.53E-07 \\
\hline
$\ell^1$ relative error & 1.03E-01 & 5.90E-02 & 3.20E-02 & 1.67E-02 \\
       \hline
    \end{tabular}
    \caption{Error table illustrating the performance of our proposed method for solving the one-dimensional HJ PDE~\eqref{eqt:eg2_L1}}
    \label{tab:eg0_1d}
\end{table}

\begin{table} [h!] 
    \footnotesize
    \centering
    \begin{tabular}{c|c|c|c|c}
    \hline
         $n_x\times n_y\times n_t$ & $20\times 20\times 11$ & $40\times 40\times 21$ & $80\times 80\times 41$ & $160\times 160\times 81$ \\
\hline
Averaged absolute residual of HJ PDE & 8.62E-07 & 9.68E-07 & 9.54E-07 & 9.82E-07 \\
\hline
$\ell^1$ relative error & 1.03E-01 & 5.74E-02 & 2.93E-02 & 1.36E-02 \\
       \hline
    \end{tabular}
    \caption{Error table illustrating the performance of our proposed method for solving the two-dimensional HJ PDE~\eqref{eqt:eg2_L1}}
    \label{tab:eg0_2d}
\end{table}

\begin{figure}[htbp]
    \centering
    \begin{subfigure}{0.45\textwidth}
        \centering \includegraphics[width=\textwidth]{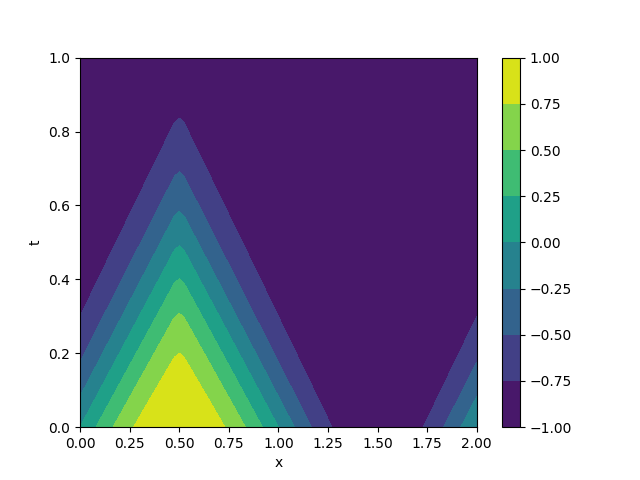}
        \caption{1D, $n_t = 41$, $n_x = 80$}
    \end{subfigure}
    \begin{subfigure}{0.45\textwidth}
        \centering \includegraphics[width=\textwidth]{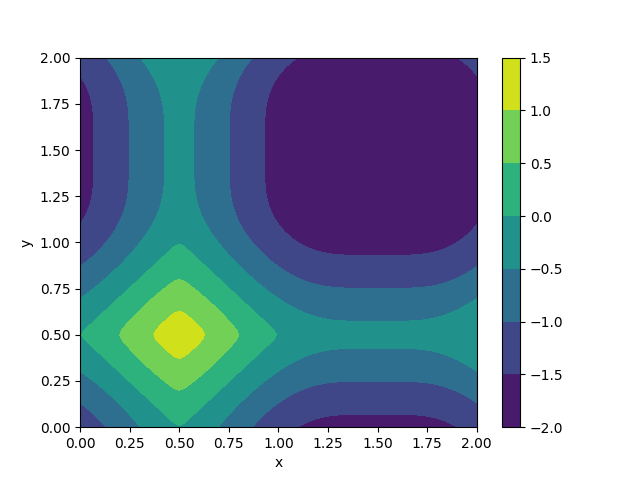}
        \caption{2D, $t=0.25$, $n_t = 41$, $n_x = n_y = 80$}
    \end{subfigure}

    \begin{subfigure}{0.45\textwidth}
        \centering \includegraphics[width=\textwidth]{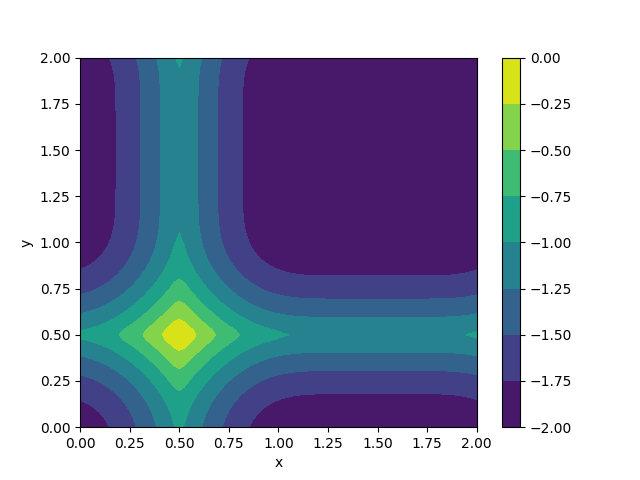}
        \caption{2D, $t=0.5$, $n_t = 41$, $n_x = n_y = 80$}
    \end{subfigure}
    \begin{subfigure}{0.45\textwidth}
        \centering \includegraphics[width=\textwidth]{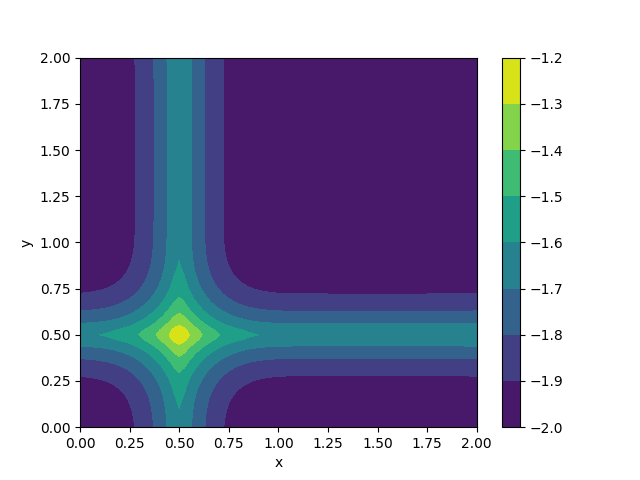}
        \caption{2D, $t=0.75$, $n_t = 41$, $n_x = n_y = 80$}
    \end{subfigure}
    
    \caption{Visualization of the solution to the one-dimensional HJ PDE~\eqref{eqt:eg2_L1} in (a), and the evolution of level sets for the two-dimensional HJ PDE solution at various time instances ($t=0.25$, $t=0.5$, and $t=0.75$) in (b)-(d).}
    \label{fig:eg0_epsl0}
\end{figure}

The errors featured in the first rows of both error tables represent the averaged residual errors of the HJ PDE. Notably, all these errors remain below the $10^{-6}$ threshold. This indicates the convergence of our method in terms of achieving an HJ PDE residual beneath the predetermined threshold of $10^{-6}$.
Turning to the second rows of the error tables, they depict the $\ell^1$ relative errors when compared to the reference solution. We observe that these errors reduce by approximately half as we increase the grid size by a factor of $2$. This phenomenon is in alignment with the utilization of a first-order Engquist-Osher scheme.

Furthermore, the level sets of the solution for the one-dimensional case are depicted in Fig~\ref{fig:eg0_epsl0} (a), while the level sets of the solution at distinct time points ($t=0.25$, $0.5$, $0.75$) are presented in Fig~\ref{fig:eg0_epsl0} (b)-(d). The results of the experiments and the error analysis indicate a higher error magnitude compared to the previous example. Consequently, we adopt a finer grid to enhance the visual representation in the figures.

\subsection{Third experiment: spatially-dependent Hamiltonian} \label{sec:eg3}
In this experiment, we solve the following HJ PDE whose Hamiltonian depends on the spatial variable $x$:
\begin{equation}\label{eqt:eg3_L1_xdep}
    \begin{dcases}
        \frac{\partial \phi(x,t)}{\partial t} + \|\nabla_x \phi(x,t)\| + f(x) = 0, & x\in [0,2]^n, t\in [0,1],\\
        \phi(x,0) = \sum_{i=1}^n \sin \pi x_i, & x\in [0,2]^n,
    \end{dcases}
    \end{equation}
where $f$ is defined by $f(x) = 3\exp(-4 \|x-1\|^2) + 1$.
We utilize Algorithms~\ref{alg:pdhg_full_1d} and~\ref{alg:pdhg_full_2d} to solve this problem. Notably, the Hamiltonian represented in~\eqref{eqt:eg3_L1_xdep} is non-separable, meaning it cannot be expressed as $\sum_{i=1}^n H_i(x,t,p_i)$ for certain functions $H_1, \dots, H_n$. Under this circumstance, we adopt the numerical Hamiltonian introduced in~\cite{osher1988fronts}. For one-dimensional instances, the numerical Hamiltonian is defined as $(\max\{p^-, 0\}^2 + \min\{p^+, 0\}^2)^{1/2} + f(x)$, whereas for two-dimensional scenarios, it is defined as $(\max\{p_1^-, 0\}^2 + \min\{p_1^+, 0\}^2 + \max\{p_2^-, 0\}^2 + \min\{p_2^+, 0\}^2)^{1/2} + f(x)$.

\begin{figure}[htbp]
    \centering
    \begin{subfigure}{0.45\textwidth}
        \centering \includegraphics[width=\textwidth]{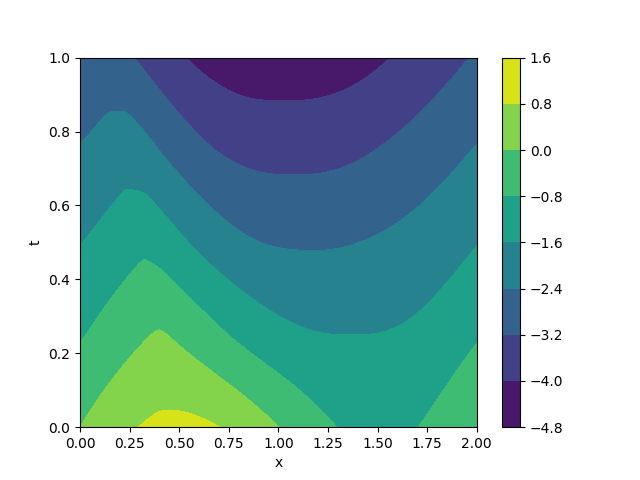}
        \caption{1D, $n_x = 80$, $n_t = 5$}
    \end{subfigure}
    \hfill
    \begin{subfigure}{0.45\textwidth}
        \centering \includegraphics[width=\textwidth]{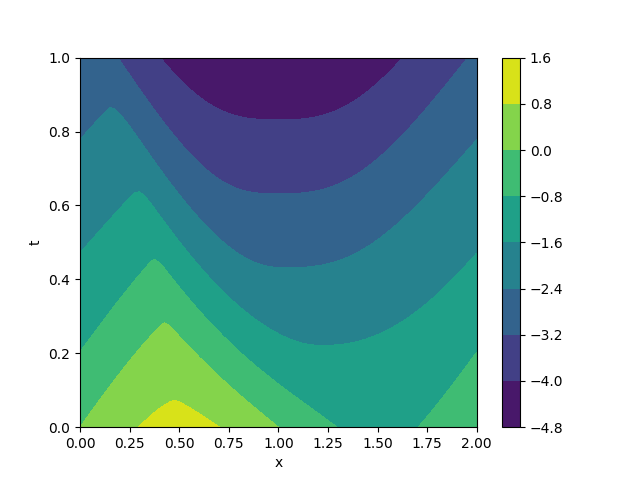}
        \caption{1D, $n_x = 80$, $n_t = 41$}
    \end{subfigure}

    \begin{subfigure}{0.45\textwidth}
        \centering \includegraphics[width=\textwidth]{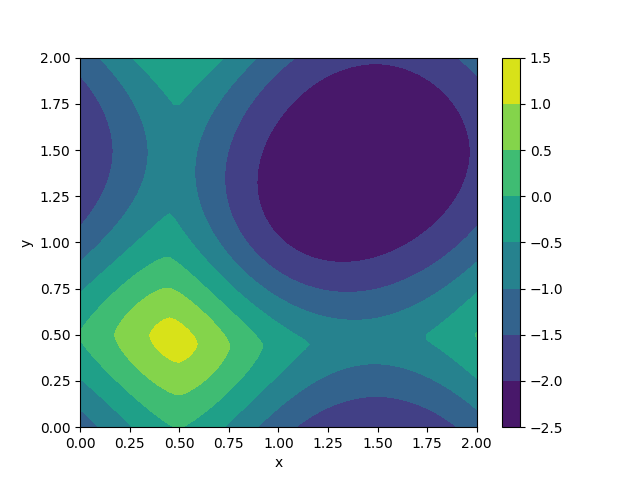}
        \caption{2D, $t=0.25$, $n_x = n_y=80$, $n_t = 41$}
    \end{subfigure}
    \hfill
    \begin{subfigure}{0.45\textwidth}
        \centering \includegraphics[width=\textwidth]{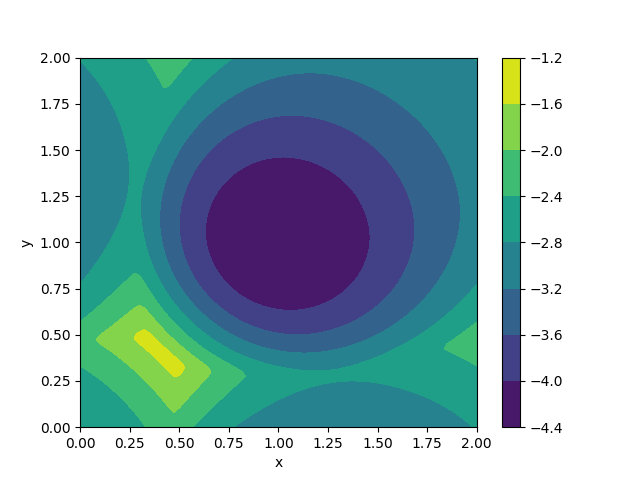}
        \caption{2D, $t=0.75$, $n_x = n_y=80$, $n_t = 41$}
    \end{subfigure}
    \caption{Contours illustrating the solution of the one-dimensional HJ PDE~\eqref{eqt:eg3_L1_xdep} in (a) and (b), along with the level sets for the two-dimensional HJ PDE in (c) and (d).
    In (a), a relatively larger time step of $\Delta t = 0.25$ is applied, whereas in (b)-(d), a smaller time step of $\Delta t = 0.025$ is employed.}
    \label{fig:eg2}
\end{figure}

The outcomes derived through our proposed approach are depicted in Figure~\ref{fig:eg2}.
In Figure~\ref{fig:eg2} (a), we present the level sets of the one-dimensional solution achieved using a relatively larger time step of $\Delta t = 0.25$, while in Figure~\ref{fig:eg2} (b), we display the level sets of the one-dimensional solution attained using a finer time step of $\Delta t = 0.025$. The outcomes of the two-dimensional solution, obtained using a time step of $\Delta t = 0.025$, are depicted in Figures~\ref{fig:eg2} (c) and (d). Figure~\ref{fig:eg2} (c) illustrates the solution at $t=0.25$, while Figure~\ref{fig:eg2} (d) displays the solution at $t=0.75$.
This example demonstrates the capability of our proposed method to effectively manage some complex Hamiltonians with spatial dependency.

\subsection{Fourth experiment: spatiotemporally-dependent Hamiltonian } \label{sec:eg4}

\begin{figure}[htbp]
    \centering
    
    \begin{subfigure}{0.45\textwidth}
        \centering \includegraphics[width=\textwidth]{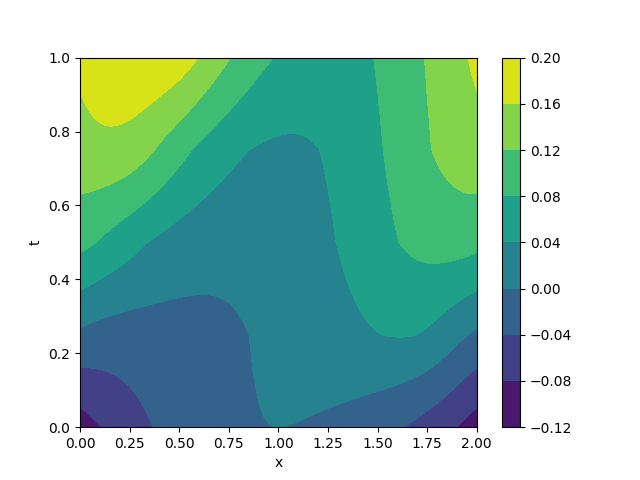}
        \caption{1D, $n_x=80$, $n_t =5$}
    \end{subfigure}
    \hfill
    \begin{subfigure}{0.45\textwidth}
        \centering \includegraphics[width=\textwidth]{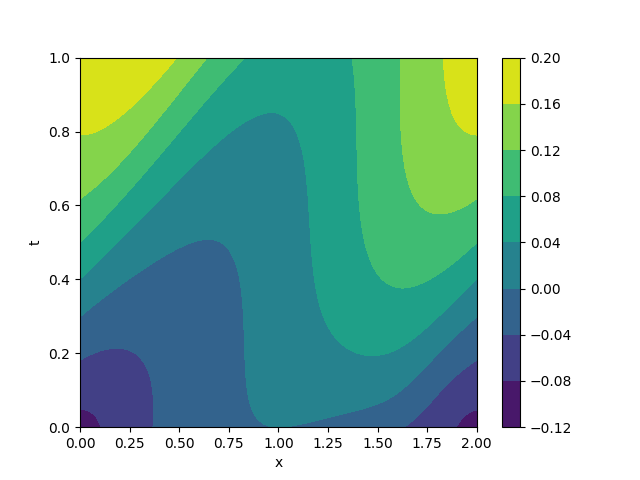}
        \caption{1D, $n_x=80$, $n_t = 41$}
    \end{subfigure}
    
    \begin{subfigure}{0.45\textwidth}
        \centering \includegraphics[width=\textwidth]{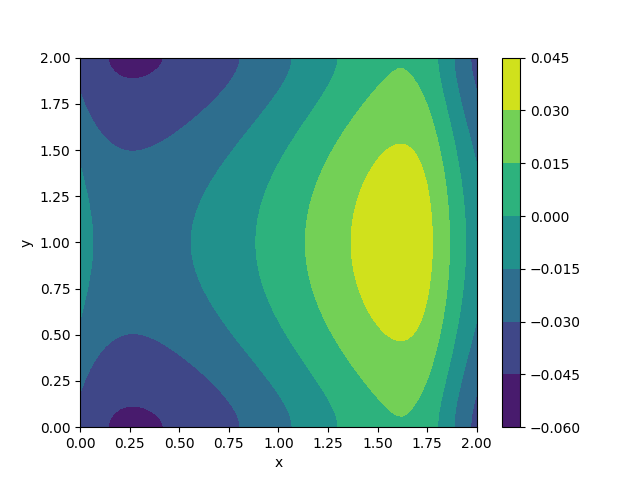}
        \caption{2D, $t=0.25$, $n_x=n_y=80$, $n_t=5$}
    \end{subfigure}
    \hfill 
    \begin{subfigure}{0.45\textwidth}
        \centering \includegraphics[width=\textwidth]{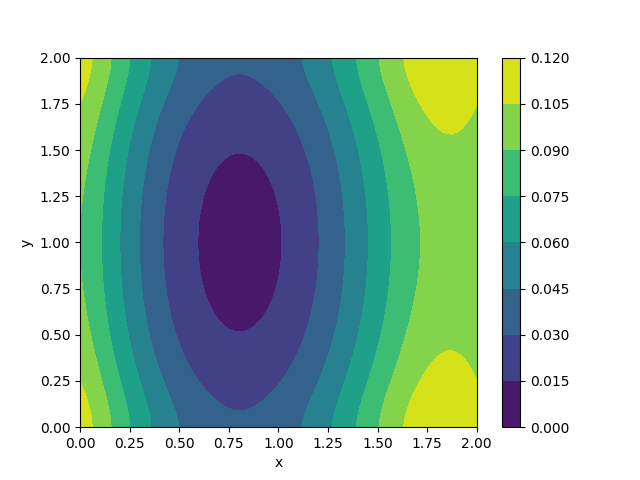}
        \caption{2D, $t=0.5$, $n_x=n_y=80$, $n_t=5$}
    \end{subfigure}

    \begin{subfigure}{0.45\textwidth}
        \centering \includegraphics[width=\textwidth]{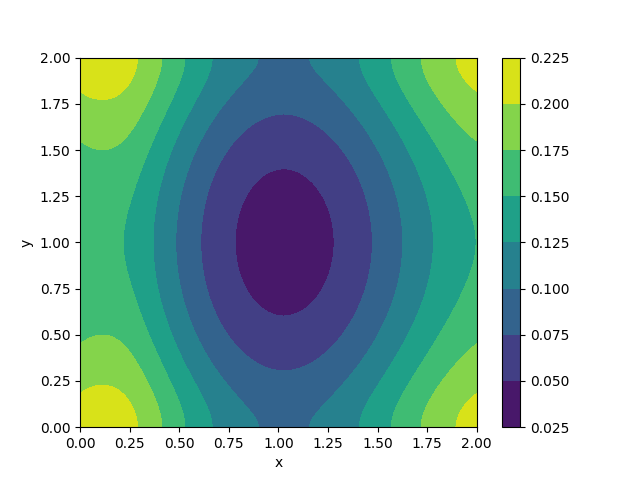}
        \caption{2D, $t=0.75$, $n_x=n_y=80$, $n_t=5$}
    \end{subfigure}
    \hfill 
    \begin{subfigure}{0.45\textwidth}
        \centering \includegraphics[width=\textwidth]{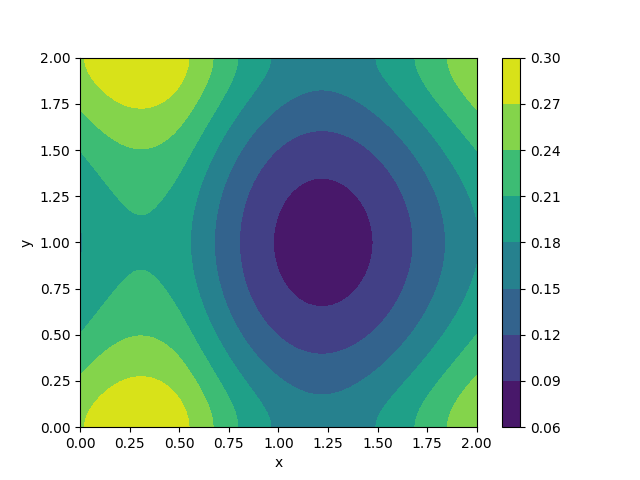}
        \caption{2D, $t=1$, $n_x=n_y=80$, $n_t=5$}
    \end{subfigure}
    
    \caption{Numeric solutions for the one-dimensional and two-dimensional HJ PDE~\eqref{eqt:eg4_xtdep} are depicted. In (a), we display the level sets of the one-dimensional solution obtained using a larger time step of $\Delta t = 0.25$, whereas in (b), we showcase the one-dimensional solution achieved with a smaller time step of $\Delta t = 0.025$. Furthermore, level sets of the two-dimensional solution are presented at distinct time points ($t=0.25$, $0.5$, $0.75$, $1.0$) in (c)-(f).}
    \label{fig:eg12}
\end{figure}

In this experiment, we solve the following viscous HJ PDE whose Hamiltonian depends on both the spatial variable $x$ and the time variable $t$:
    \begin{equation}\label{eqt:eg4_xtdep}
    \begin{dcases}
        \frac{\partial \phi(x,t)}{\partial t} + \frac{1}{2}\|\nabla_x \phi(x,t)\|^2 + f(x,t) = 0.1 \Delta_x \phi(x,t), & x\in [0,2]^n, t\in [0,1],\\
        \phi(x,0) = -\frac{1}{10}\|x-1\|^2, & x\in [0,2]^n,
    \end{dcases}
    \end{equation}
    where $f(x,t) = -\frac{1}{2}\min\{(x_1- t - 0.5)^2, (x_1- t + 1.5)^2, (x_1- t - 2.5)^2\} - \frac{1}{4}\sum_{j=2}^n(x_j-1)^2$.
We apply Algorithms~\ref{alg:pdhg_full_1d} and~\ref{alg:pdhg_full_2d} to solve this problem. We use the Engquist-Osher scheme in the saddle point formulation. To be specific, we set the numerical Hamiltonian $\hat H(x,t,p^+,p^-) = H_-(p^+) + H_+(p^-) + f(x,t)$ for the one-dimensional case and $\hat H(x,t,p_1^+,p_1^-, p_2^+,p_2^-) = H_-(p_1^+) + H_+(p_1^-) + H_-(p_2^+) + H_+(p_2^-) + f(x,t)$ for the two-dimensional case, where $H_-(p^+) = \frac{1}{2}\min\{p^+,0\}^2$ and $H_+(p^-) = \frac{1}{2}\max\{p^-,0\}^2$.

The solution is illustrated in Fig~\ref{fig:eg12}. In (a)-(b), we present the level sets depicting the solution to the one-dimensional HJ PDE. In (a), a larger time step of $\Delta t = 0.25$ is utilized, while in (b), a smaller time step of $\Delta t = 0.025$ is employed. Additionally, we provide contour plots showcasing the solution to the two-dimensional PDE at distinct time instances ($t=0.25$, $0.5$, $0.75$, $1.0$) in (c)-(f). Due to computational complexity, we adopt a larger time step of $\Delta t = 0.25$ in the two-dimensional case. This example demonstrates the capability of our proposed method to handle intricate Hamiltonians dependent on $(x,t)$.

\section{Summary}\label{sec:summary}

This paper solves HJ PDEs using a saddle point formulation, which is solved using PDHG method. We provide numerical validation that this method can compute the viscosity solutions with errors related to the grid size. We can handle certain Hamiltonians which depend on $(x,t)$. Moreover, we use implicit time discretization, which circumvents the restrictive CFL time step restriction for explicit methods. In other words, we can choose big time step to speed up the computation. 
The merit of this approach lies in the simplicity of the saddle point formulation. This simple formulation is achieved by capitalizing on the Fenchel-Legendre transform and the duality inherent in HJ PDEs. This simplicity facilitates updates within our method to have either explicit formulations or be conducive to parallel computation. In a special case where the Hamiltonian $H(x,t,p)$ is separable and 1-homogeneous with respect to $p$, the saddle point formulation takes the standard form in~\cite{chambolle2011first}, ensuring the convergence of the proposed algorithms through PDHG theory (see Remark~\ref{rem:B1_convergencePDHG}).

Although it is a first-order method, it has the potential to serve as an initialization for more accurate methods, particularly in applications that demand smaller errors. It may be also an interesting future direction to combine this method with high order schemes.
Moreover, this method converts an equation to a saddle point problem which can fit pretty well under the framework of machine learning. For solving HJ PDEs or problems related to HJ PDEs, the formulas provided in this paper can provide some ideas on design of loss functions.

\section*{Acknowledgement}
We would like to express our thanks to Dr. Liu Yang for providing the Jax code of Algorithm~\ref{alg:pdhg_full_m_1d_L1}, which greatly assisted in the numerical implementation of our research. 

\bibliographystyle{siam} 
\bibliography{references}

\appendix
\section{More details about the algorithms}
\subsection{More discussion about the continuous saddle point formula}\label{appendix:cont}
The first order optimization condition for the saddle point problem in~\eqref{eqt:saddle_cont_appendix} is
\begin{equation}
\begin{dcases}
\frac{\partial \phi(x,t)}{\partial t} + H(x, t,\nabla_x \phi(x,t)) = \epsilon\Delta_x \phi(x,t), & x\in \Omega, t\in[0,T],\\
\frac{\partial \rho(x,t)}{\partial t} + \nabla_x \cdot (\rho(x,t)\nabla_p H(x, t,\nabla_x \phi(x,t))) + \epsilon\Delta_x \rho(x,t)=0, & x\in \Omega, t\in[0,T],\\
\phi(x,0) = g(x),\quad \rho(x,T) = c, & x\in \Omega,
\end{dcases}
\end{equation}
which is similar to the mean-field control problem, except that we do not restrict $\rho$ to be non-negative. 
There is a gap between~\eqref{eqt:saddle_cont_appendix} and~\eqref{eqt:saddle_cont}, namely, whether we require $\rho$ to be non-negative or not. Although the theoretical understanding of when this condition holds is lacking, numerically we have found that our algorithm works and computes the viscosity solution.

Note that in the last line of~\eqref{eqt:saddle_cont_appendix}, $v$ is a vector, while $v$ takes the form of a function in~\eqref{eqt:saddle_cont}. We abuse the notation $v$ here. In~\eqref{eqt:saddle_cont_appendix}, there are infinitely many finite-dimensional optimization problems, which depends on $(x,t)$. Upon their consolidation into a single optimization problem, the variable transforms into a function that depends on $(x,t)$, denoted as $v(x,t)$ in~\eqref{eqt:saddle_cont}, whose value corresponds to the original variable $v$ in~\eqref{eqt:saddle_cont_appendix}.

\subsection{One-dimensional semi-discrete method}\label{sec:appendix_semidisc_1d}
To solve the semi-discrete equation~\eqref{eqt:semidisc_HJ_equation}, we propose to solve the saddle point problem~\eqref{eqt:saddle_pt_general_semidisc}, whose first order optimality condition is
\begin{equation}
\begin{dcases}
\dot \phi_i(t) + \hat H(x_i, t, (D_x^+ \phi)_i(t), (D_x^- \phi)_i(t)) \leq \epsilon (D_{xx} \phi)_i(t), & i=1,\dots, n_x; t\in[0,T],\\
\dot\rho_i(t) + (D_x^- (\rho v^+))_i(t) + (D_x^+ (\rho v^-))_i(t) + \epsilon (D_{xx} \rho)_i(t)=0, & i=1,\dots, n_x; t\in[0,T],\\
(v_i^+(t), v_i^-(t)) = \nabla_{(p^+, p^-)} \hat H(x_i, t, (D_x^+ \phi)_i(t), (D_x^- \phi)_i(t)),& i=1,\dots, n_x; t\in[0,T],\\
\phi_i(0) = g(x_i),\quad \rho_i(T) = c, & i=1,\dots, n_x,\\
\rho_i(t)\left(\dot \phi_i(t) + \hat H(x_i, t, (D_x^+ \phi)_i(t), (D_x^- \phi)_i(t)) - \epsilon (D_{xx} \phi)_i(t)\right)=0, & i=1,\dots, n_x; t\in[0,T].
\end{dcases}
\end{equation}

Similar to the continuous setting, if the corresponding mean-field control problem possesses a solution with a positive $\rho$, the inequality in the first row becomes an equality. This implies that the proposed saddle point problem addresses the semi-discrete HJ PDE~\eqref{eqt:semidisc_HJ_equation}. This inclusively covers situations where a positive diffusion coefficient $\epsilon$ is present.
It's important to note that by ensuring an appropriate discretization for the two PDEs in~\eqref{eqt:2pde_cont}, both the discretize-then-optimize and optimize-then-discretize approaches yield the same method.
These comments hold true for the following sections and we will not repeat them.

The algorithm for the proposed one-dimensional semi-discrete method is outlined in Algorithm~\ref{alg:pdhg_semi_1d}, where we denote the objective function in~\eqref{eqt:saddle_pt_general_semidisc} as $\mathcal{L}_{semi}$.

\begin{algorithm}[htbp]
\SetAlgoLined
\SetKwInOut{Input}{Inputs}
\SetKwInOut{Output}{Outputs}
\Input{Stepsize $\tau, \sigma>0$, error tolerance $\delta>0$, inner maximal iteration number $N_{inner}$ and outer maximal iteration number $N_{outer}$.}
\Output{Solution to the saddle point problem~\eqref{eqt:saddle_pt_general_semidisc}.}
For all $i=1,\dots, n_x$, initialize the functions by $\phi_i^0(t)=g(x_i), \forall t\in[0,T]$, $\rho_i^0 \equiv c$, $v_i^{0,+}=v_i^{0,-}\equiv 0$.

 \For{$\ell = 0,1,\dots,N_{outer}-1$}{
 Update the functions $\phi_i\colon [0,T]\to \R$ for all $i=1,\dots, n_x$ by
 \begin{equation}
 {\small
 \begin{split}
 (\phi_i^{\ell+1})_i &= \argmin_{\phi_i \forall i : \phi_i(0)=g(x_i)} \mathcal{L}_{semi}((\phi_i)_i, (\rho^\ell_i)_i, (v^{\ell,+}_i)_i, (v^{\ell,-}_i)_i) + \frac{1}{2\tau} \sum_{i=1}^{n_x}\left(\|\dot \phi_i - \dot \phi_i^\ell\|^2 +\| (D_x^+\phi)_i - (D_x^+\phi^\ell)_i\|^2\right) \\
 &= (\phi_i^\ell)_i + \tau (-\partial_t^2 - D_{xx})^{-1}\left(\dot\rho^\ell + D_x^- (\rho^\ell v^{\ell,+}) + D_x^+ (\rho^\ell v^{\ell,-}) + \epsilon D_{xx} \rho^\ell\right),
 \end{split}
 }
 \end{equation}
 where $(-\partial_t^2 - D_{xx})^{-1}(f_i)_i$ (here each $f_i$ is a function of $t$) denotes the solution $(u_i)_i$ to the ODE system $- \ddot u_i - (D_{xx} u)_i=f_i$ in $[0,T]$ for all $i=1,\dots, n_x$ with periodic spatial condition, Dirichlet initial condition $u_i(0)=0$, and Neumann terminal condition $\dot u_i(T)=0$.
 
 \If{$\sum_{i=1}^{n_x}\|\dot \phi_i^{\ell+1}(t) + \hat H(x_i, t, (D_x^+ \phi^{\ell+1})_i(t), (D_x^- \phi^{\ell+1})_i(t)) - \epsilon (D_{xx} \phi^{\ell+1})_i(t)\|_1 \leq \delta$}{
   Return $(\phi_i^{\ell+1})_i$.
   }
 
 Set $\bar\phi_i^{\ell+1} = 2\phi_i^{\ell+1} - \phi_i^{\ell}$ for all $i=1,\dots, n_x$.
 
 Set $v_i^{\ell+1,0,+} = v_i^{\ell,+}$, $v_i^{\ell+1,0,-} = v_i^{\ell,-}$, $\rho_i^{\ell+1,0} = \rho_i^{\ell}$ for all $i=1,\dots, n_x$.

  \For{$m=0,1,\dots, N_{inner}-1$}{
 Update $v_i^+$ and $v_i^-$ for all $i=1,\dots, n_x$ by
 \begin{equation}
 {\scriptsize
 \begin{split}
     &(v_i^{\ell+1,m+1,+}, v_i^{\ell+1,m+1,-})_i = \argmax_{v_i^+, v_i^- \forall i} \Bigg\{\mathcal{L}_{semi}((\bar\phi^{\ell+1}_i)_i, (\rho^{\ell+1,m}_i)_i, (v^+_i)_i, (v^-_i)_i)\\
     &\quad\quad\quad\quad\quad\quad\quad\quad
     - \frac{1}{2\sigma} \sum_{i=1}^{n_x}\left(\|\rho_i^{\ell+1,m} (v_i^+ - v_i^{\ell+1,m,+})\|^2 + \|\rho_i^{\ell+1,m} (v_i^- - v_i^{\ell+1,m,-})\|^2\right)\Bigg\}  \\
    &= \argmin_{v_i^+, v_i^- \forall i} \Bigg\{\sum_{i=1}^{n_x}\int_0^T \hat H^*(x_i,t, v_i^+(t), v_i^-(t)) + \frac{\rho^{\ell+1,m}_i(t)}{2\sigma} \left(v_i^+(t) - v_i^{\ell+1,m,+}(t) - \sigma \frac{(D_x^+ \bar\phi^{\ell+1})_i(t)}{\rho^{\ell+1,m}_i(t)} \right)^2\\
    &\quad\quad\quad\quad\quad\quad\quad\quad +\frac{\rho^{\ell+1,m}_i(t)}{2\sigma} \left(v_i^-(t) - v_i^{\ell+1,m,-}(t) - \sigma \frac{(D_x^- \bar\phi^{\ell+1})_i(t)}{\rho^{\ell+1,m}_i(t)} \right)^2 dt\Bigg\}.
    \end{split}
    }
 \end{equation}
 Update $\rho_i$ for all $i=1,\dots, n_x$ by 
 \begin{equation}
 {\small
 \begin{split}
     (\rho_i^{\ell+1,m+1})_i &= \argmax_{\rho_i\forall i: \rho_i \geq 0} \mathcal{L}_{semi}((\bar\phi^{\ell+1}_i)_i, (\rho_i)_i, (v^{\ell+1,m+1,+}_i)_i, (v^{\ell+1,m+1,-}_i)_i) - \frac{1}{2\sigma} \sum_{i=1}^{n_x}\|\rho_i - \rho_i^{\ell+1,m}\|^2\\
    &= (\max\{\mu_i^{\ell+1,m+1}, 0\})_i,
    \end{split}
    }
 \end{equation}
 where $\mu_i^{\ell+1,m+1}$ is a function defined by $\mu_i^{\ell+1,m+1} = \rho_i^{\ell+1,m} + \sigma (\dot {\bar\phi}_i^{\ell+1}(t) + v_i^{\ell+1,m+1,+}(t)(D_x^+\bar\phi^{\ell+1})_i(t) + v_i^{\ell+1,m+1,-}(t)(D_x^-\bar\phi^{\ell+1})_i(t) - \hat H^*(x_i, t, v_i^{\ell+1,m+1,+}(t), v_i^{\ell+1,m+1,-}(t)) - \epsilon (D_{xx} \bar\phi^{\ell+1})_i(t))$.
 }
 Set $v^{\ell+1,+}_i = v^{\ell+1,N_{inner},+}_i$, $v^{\ell+1,-}_i = v^{\ell+1,N_{inner},-}_i$, $\rho^{\ell+1}_i = \rho^{\ell+1,N_{inner}}_i$ for all $i=1,\dots, n_x$.
 }
 Return $(\phi_i^{N_{outer}})_i$.
 \caption{The proposed algorithm for solving~\eqref{eqt:saddle_pt_general_semidisc}\label{alg:pdhg_semi_1d}}
\end{algorithm}

\subsection{One-dimensional fully-discrete method}\label{appendix:fulldisc_1d}
To solve the fully-discrete equation~\eqref{eqt:disc_HJ_equation}, we propose solving the saddle point problem~\eqref{eqt:saddle_pt_general_fulldisc}, the first-order optimality condition of which is given by:
\begin{equation} \label{eqt:2pdes_full_discrete_1d}
\begin{dcases}
(D_t^- \phi)_{i,k} + \hat H(x_i, t_k,(D_x^+ \phi)_{i,k}, (D_x^- \phi)_{i,k}) \leq \epsilon (D_{xx} \phi)_{i,k}, & i=1,\dots, n_x; k = 2,\dots, n_t,\\
(D_t^+\rho)_{i,k} + (D_x^- (\rho v^+))_{i,k} + (D_x^+ (\rho v^-))_{i,k} + \epsilon (D_{xx} \rho)_{i,k}=0, & i=1,\dots, n_x; k = 1,\dots, n_t-1,\\
(v_{i,k}^+, v_{i,k}^-) = \nabla_{(p^+, p^-)} \hat H(x_i, t_{k+1},(D_x^+ \phi)_{i,k+1}, (D_x^- \phi)_{i,k+1}),& i=1,\dots, n_x; k = 1,\dots, n_t-1,\\
\phi_{i,1} = g(x_i),\quad  \rho_{i,n_t} = c, & i=1,\dots, n_x,\\
\rho_{i,k-1}\left((D_t^- \phi)_{i,k} + \hat H(x_i, t_k,(D_x^+ \phi)_{i,k}, (D_x^- \phi)_{i,k}) - \epsilon (D_{xx} \phi)_{i,k}\right) = 0, & i=1,\dots, n_x; k = 2,\dots, n_t.
\end{dcases}
\end{equation}
The proposed algorithm for solving~\eqref{eqt:saddle_pt_general_fulldisc} is presented in Algorithm~\ref{alg:pdhg_full_1d}. Within this algorithm, we employ $\mathcal{L}_{full}$ to represent the objective function in~\eqref{eqt:saddle_pt_general_fulldisc}. During the updating of $\phi$, the term $(D_t^+\rho)_{i,n_t-1}$ represents $\frac{c - \rho_{i,n_t-1}}{\Delta t}$, aligning with the terminal condition for $\rho$ in~\eqref{eqt:2pdes_full_discrete_1d}.

\begin{algorithm}[htbp]
\SetAlgoLined
\SetKwInOut{Input}{Inputs}
\SetKwInOut{Output}{Outputs}
\Input{Stepsize $\tau, \sigma>0$, error tolerance $\delta>0$, inner maximal iteration number $N_{inner}$ and outer maximal iteration number $N_{outer}$.}
\Output{Solution to the to the saddle point problem~\eqref{eqt:saddle_pt_general_fulldisc}.}
For each $i=1,\dots, n_x$, initialize the matrices by $\phi_{i,k}^0=g(x_i)$ for $k=1,\dots, n_t$, $\rho_{i,k}^0 = c$, $v_{i,k}^{0,+} = v_{i,k}^{0,-}= 0$ for $k=1,\dots, n_t-1$.

 \For{$\ell = 0,1,\dots,N_{outer}-1$}{
 Update the matrix $\phi_{i,k}$ for $i=1,\dots, n_x$; $k=1,\dots, n_t$ by
 \begin{equation}
 {\small
 \begin{split}
 (\phi_{i,k}^{\ell+1})_{i,k} &= \argmin_{\substack{\phi_{i,k} \forall i,k \\ \phi_{i,1}=g(x_i)}} \mathcal{L}_{full}((\phi_{i,k})_{i,k}, (\rho^\ell_{i,k})_{i,k}, (v^{\ell,+}_{i,k})_{i,k}, (v^{\ell,-}_{i,k})_{i,k})  \\
 & \quad\quad\quad\quad
 + \frac{1}{2\tau} \sum_{i=1}^{n_x}\sum_{i=2}^{n_t}\left(((D_t^- \phi)_{i,k} - (D_t^- \phi^\ell)_{i,k})^2 +( (D_x^+\phi)_{i,k} - (D_x^+\phi^\ell)_{i,k})^2 \right) \\
 &= (\phi_{i,k}^\ell)_{i,k} + \tau (-D_{tt} - D_{xx})^{-1}\left(D_t^+\rho^\ell + D_x^- (\rho^\ell v^{\ell,+}) + D_x^+ (\rho^\ell v^{\ell,-}) + \epsilon D_{xx} \rho^\ell\right),
 \end{split}
 }
 \end{equation}
 where $(-D_{tt} - D_{xx})^{-1}(f_{i,k})_{i,k}$ (for a matrix $f$ with elements $f_{i,k}$ and the linear operator $D_{tt}f = (\frac{f_{i,k-1} -2f_{i,k} + f_{i,k+1}}{\Delta t^2})_{i,k}$) denotes the solution $u$ to the linear system $- (D_{tt} u)_{i,k+1} - (D_{xx} u)_{i,k+1} = f_{i,k}$ for all $i=1,\dots, n_x$; $k=1,\dots, n_t-1$ with periodic spatial condition, Dirichlet initial condition $u_{i,1}=0$, and Neumann terminal condition $u_{i,n_t+1}=u_{i,n_t}$. 
 
 \If{$\sum_{i=1}^{n_x}\sum_{k=2}^{n_t}|(D_t^- \phi^{\ell+1})_{i,k} + \hat H(x_i, t_k,(D_x^+ \phi^{\ell+1})_{i,k}, (D_x^- \phi^{\ell+1})_{i,k}) - \epsilon (D_{xx} \phi^{\ell+1})_{i,k}| \leq \delta$}{
   Return $(\phi_{i,k}^{\ell+1})_{i,k}$.
   }
 
 Set $\bar\phi_{i,k}^{\ell+1} = 2\phi_{i,k}^{\ell+1} - \phi_{i,k}^{\ell}$ for all $i=1,\dots, n_x$; $k=2,\dots, n_t$.
 
 Set $v_{i,k}^{\ell+1,0,+} = v_{i,k}^{\ell,+}$, $v_{i,k}^{\ell+1,0,-} = v_{i,k}^{\ell,-}$, $\rho_{i,k}^{\ell+1,0} = \rho_{i,k}^{\ell}$ for all $i=1,\dots, n_x$; $k=1,\dots, n_t-1$.

  \For{$m=0,1,\dots, N_{inner}-1$}{
 Update $v_{i,k}^+$ and $v_{i,k}^-$ for all $i=1,\dots, n_x$; $k=1,\dots, n_t-1$ by
 \begin{equation}
 {\scriptsize
 \begin{split}
     &((v_{i,k}^{\ell+1,m+1,+}, v_{i,k}^{\ell+1,m+1,-}))_{i,k} = \argmax_{v_{i,k}^+, v_{i,k}^- \forall i,k} \Bigg\{\mathcal{L}_{full}((\bar\phi^{\ell+1}_{i,k})_{i,k}, (\rho^{\ell+1,m}_{i,k})_{i,k}, (v^+_{i,k})_{i,k}, (v^-_i)_{i,k})\\
     &\quad\quad\quad\quad\quad\quad\quad\quad
     - \frac{1}{2\sigma} \sum_{i=1}^{n_x}\sum_{k=1}^{n_t-1}\left((\rho_{i,k}^{\ell+1,m} (v_{i,k}^+ - v_{i,k}^{\ell+1,m,+}))^2 + (\rho_{i,k}^{\ell+1,m} (v_{i,k}^- - v_{i,k}^{\ell+1,m,-}))^2\right)\Bigg\}  \\
    &= \argmin_{v_{i,k}^+, v_{i,k}^- \forall i,k} \Bigg\{\sum_{i=1}^{n_x}\sum_{k=1}^{n_t-1} \Bigg(\hat H^*(x_i,t_{k+1},v_{i,k}^+, v_{i,k}^-) + \frac{\rho^{\ell+1,m}_{i,k}}{2\sigma} \left(v_{i,k}^+ - v_{i,k}^{\ell+1,m,+} - \sigma \frac{(D_x^+ \bar\phi^{\ell+1})_{i,k+1}}{\rho^{\ell+1,m}_{i,k}} \right)^2\\
    &\quad\quad\quad\quad\quad\quad\quad\quad +\frac{\rho^{\ell+1,m}_{i,k}}{2\sigma} \left(v_{i,k}^- - v_{i,k}^{\ell+1,m,-} - \sigma \frac{(D_x^- \bar\phi^{\ell+1})_{i,k+1}}{\rho^{\ell+1,m}_{i,k}} \right)^2 \Bigg)\Bigg\}.
    \end{split}
    }
 \end{equation}
 Update $\rho_{i,k}$ for all $i=1,\dots, n_x$; $k=1,\dots, n_t-1$ by 
 \begin{equation}
 {\small
 \begin{split}
     (\rho_{i,k}^{\ell+1,m+1})_{i,k} &= \argmax_{\rho_{i,k}\forall i,k: \rho_{i,k} \geq 0} \mathcal{L}_{full}((\bar\phi^{\ell+1}_{i,k})_{i,k}, (\rho_{i,k})_{i,k}, (v^{\ell+1,m+1,+}_{i,k})_{i,k}, (v^{\ell+1,m+1,-}_{i,k})_{i,k}) \\
     &\quad\quad\quad\quad\quad\quad\quad\quad
     - \frac{1}{2\sigma} \sum_{i=1}^{n_x}\sum_{k=1}^{n_t-1}(\rho_{i,k} - \rho_{i,k}^{\ell+1,m})^2\\
    &= (\max\{\mu_{i,k}^{\ell+1,m+1}, 0\})_{i,k},
    \end{split}
    }
 \end{equation}
 where $\mu_{i,k}^{\ell+1,m+1}$ is a number defined by $\mu_{i,k}^{\ell+1,m+1} = \rho_{i,k}^{\ell+1,m} + \sigma ((D_t^- \bar\phi^{\ell+1})_{i,k+1} + v_{i,k}^{\ell+1,m+1,+}(D_x^+\bar\phi^{\ell+1})_{i,k+1} + v_{i,k}^{\ell+1,m+1,-}(D_x^-\bar\phi^{\ell+1})_{i,k+1} - \hat H^*(x_i, t_{k+1},v_{i,k}^{\ell+1,m+1,+}, v_{i,k}^{\ell+1,m+1,-}) - \epsilon (D_{xx} \bar\phi^{\ell+1})_{i,k+1})$.
 }
 Set $v^{\ell+1,+}_{i,k} = v^{\ell+1,N_{inner},+}_{i,k}$, $v^{\ell+1,-}_{i,k} = v^{\ell+1,N_{inner},-}_{i,k}$, $\rho^{\ell+1}_{i,k} = \rho^{\ell+1,N_{inner}}_{i,k}$ for all $i=1,\dots, n_x$; $k=1,\dots, n_t-1$.
 }
 Return $(\phi_{i,k}^{N_{outer}})_{i,k}$.
 \caption{The proposed algorithm for solving~\eqref{eqt:saddle_pt_general_fulldisc}\label{alg:pdhg_full_1d}}
\end{algorithm}

\subsection{Two-dimensional semi-discrete method}\label{appendix:semidisc_2d}
To solve the semi-discrete equation~\eqref{eqt:semidisc_HJ_equation_2d}, we propose solving the saddle point problem~\eqref{eqt:saddle_pt_general_semidisc_2d}, whose first order optimality condition is
\begin{equation}
{\small
\begin{dcases}
\dot \phi_{i,j}(t) + \hat H(x_i,y_j, t, (D_x^+ \phi)_{i,j}(t), (D_x^- \phi)_{i,j}(t), (D_y^+ \phi)_{i,j}(t), (D_y^- \phi)_{i,j}(t)) \\
\quad\quad\quad\quad \leq \epsilon (D_{xx} \phi)_{i,j}(t)+\epsilon (D_{yy} \phi)_{i,j}(t), & i=1,\dots, n_x; j=1,\dots, n_y; t\in[0,T],\\
\dot\rho_{i,j}(t) + (D_x^- (\rho v^+))_{i,j}(t) + (D_x^+ (\rho v^-))_{i,j}(t) + (D_y^- (\rho w^+))_{i,j}(t) \\
\quad\quad\quad\quad + (D_y^+ (\rho w^-))_{i,j}(t) + \epsilon (D_{xx} \rho)_{i,j}(t) + \epsilon (D_{yy} \rho)_{i,j}(t) =0, & i=1,\dots, n_x; j=1,\dots, n_y; t\in[0,T],\\
(v_{i,j}^+(t), v_{i,j}^-(t), w_{i,j}^+(t), w_{i,j}^-(t)) = \nabla_{(p_1^+, p_1^-, p_2^+, p_2^-)} \hat H(x_i,y_j, t,  \\
\quad\quad\quad\quad (D_x^+ \phi)_{i,j}(t),(D_x^- \phi)_{i,j}(t), (D_y^+ \phi)_{i,j}(t), (D_y^- \phi)_{i,j}(t)),& i=1,\dots, n_x; j=1,\dots, n_y; t\in[0,T],\\
\phi_{i,j}(0) = g(x_i,y_j),\quad \rho_{i,j}(T) = c, & i=1,\dots, n_x; j=1,\dots, n_y,\\
\rho_{i,j}(t)\Big(\dot \phi_{i,j}(t) + \hat H(x_i,y_j, t, (D_x^+ \phi)_{i,j}(t), (D_x^- \phi)_{i,j}(t), (D_y^+ \phi)_{i,j}(t), \\
\quad\quad\quad\quad  (D_y^- \phi)_{i,j}(t))  -\epsilon (D_{xx} \phi)_{i,j}(t)-\epsilon (D_{yy} \phi)_{i,j}(t) \Big)=0, & i=1,\dots, n_x; j=1,\dots, n_y; t\in[0,T].
\end{dcases}
}
\end{equation}

The algorithm for the proposed two-dimensional semi-discrete method is outlined in Algorithm~\ref{alg:pdhg_semi_2d}. For the sake of simplicity, we use $\alpha$ to represent $(\alpha_{i,j})_{i,j}$ for any given quantity $\alpha$. This algorithm closely resembles the one-dimensional case, and thus certain details have been omitted. We adopt the quadratic function as the penalty term for all functions except $\phi$. For $\phi$, our choice of penalty term for preconditioning is $\sum_{i=1}^{n_x}\sum_{j=1}^{n_y}\left(\|\dot \phi_{i,j} - \dot \phi_{i,j}^\ell\|^2 +\| (D_x^+\phi)_{i,j} - (D_x^+\phi^\ell)_{i,j}\|^2 +\| (D_y^+\phi)_{i,j} - (D_y^+\phi^\ell)_{i,j}\|^2\right)$.

\begin{algorithm}[htbp]
\SetAlgoLined
\SetKwInOut{Input}{Inputs}
\SetKwInOut{Output}{Outputs}
\Input{Stepsize $\tau, \sigma>0$, error tolerance $\delta>0$, inner maximal iteration number $N_{inner}$ and outer maximal iteration number $N_{outer}$.}
\Output{Solution to the saddle point problem~\eqref{eqt:saddle_pt_general_semidisc_2d}.}
For each $i=1,\dots, n_x$; $j=1,\dots, n_y$, initialize the functions by $\phi_{i,j}^0(t)=g(x_i,y_j)$ for all $t\in[0,T]$, $\rho_{i,j}^0 \equiv c$, $v_{i,j}^{0,+}=v_{i,j}^{0,-}=w_{i,j}^{0,+}=w_{i,j}^{0,-}\equiv 0$.

 \For{$\ell = 0,1,\dots,N_{outer}-1$}{
 Update the functions $\phi_{i,j}\colon [0,T]\to \R$ for all $i=1,\dots, n_x$; $j=1,\dots, n_y$ by
 \begin{equation}
 {\scriptsize
 \begin{split}
 \phi^{\ell+1} 
 = \phi^\ell + \tau (-\partial_t^2 - D_{xx}- D_{yy})^{-1}\left(\dot\rho^\ell + D_x^- (\rho^\ell v^{\ell,+}) + D_x^+ (\rho^\ell v^{\ell,-}) + D_y^- (\rho^\ell w^{\ell,+}) + D_y^+ (\rho^\ell w^{\ell,-}) + \epsilon (D_{xx} + D_{yy}) \rho^\ell\right),
 \end{split}
 }
 \end{equation}
 where $(-\partial_t^2 - D_{xx}- D_{yy})^{-1}f$ (here each $f_{i,j}$ is a function of $t$) denotes the solution $u$ to the ODE system $- \ddot u_{i,j} - (D_{xx} u + D_{yy} u)_{i,j}=f_{i,j}$ in $[0,T]$ for all $i=1,\dots, n_x$; $j=1,\dots, n_y$ with periodic spatial condition, Dirichlet initial condition $u_{i,j}(0)=0$, and Neumann terminal condition $\dot u_{i,j}(T)=0$.
 
 \If{$\sum_{i=1}^{n_x}\sum_{j=1}^{n_y}\|\dot \phi_{i,j}^{\ell+1}(t) + \hat H(x_i,y_j, t, (D_x^+ \phi^{\ell+1})_{i,j}(t), (D_x^- \phi^{\ell+1})_{i,j}(t), (D_y^+ \phi^{\ell+1})_{i,j}(t), (D_y^- \phi^{\ell+1})_{i,j}(t)) - \epsilon (D_{xx} \phi^{\ell+1} + D_{yy} \phi^{\ell+1})_{i,j}(t) \|_1 \leq \delta$}{
   Return $\phi^{\ell+1}$.
   }
 
 Set $\bar\phi_{i,j}^{\ell+1} = 2\phi_{i,j}^{\ell+1} - \phi_{i,j}^{\ell}$ for all $i=1,\dots, n_x$; $j=1,\dots, n_y$.
 
 Set $v_{i,j}^{\ell+1,0,+} = v_{i,j}^{\ell,+}$, $v_{i,j}^{\ell+1,0,-} = v_{i,j}^{\ell,-}$, $w_{i,j}^{\ell+1,0,+} = w_{i,j}^{\ell,+}$, $w_{i,j}^{\ell+1,0,-} = w_{i,j}^{\ell,-}$, $\rho_{i,j}^{\ell+1,0} = \rho_{i,j}^{\ell}$ for all $i=1,\dots, n_x$; $j=1,\dots, n_y$.

  \For{$m=0,1,\dots, N_{inner}-1$}{
 Update $v_{i,j}^+$, $v_{i,j}^-$, $w_{i,j}^+$ and $w_{i,j}^-$ for all $i=1,\dots, n_x$; $j=1,\dots, n_y$ by
 \begin{equation}
 {\scriptsize
 \begin{split}
     &(v^{\ell+1,m+1,+}, v^{\ell+1,m+1,-}, w^{\ell+1,m+1,+}, w^{\ell+1,m+1,-}) \\
    =& \argmin_{v_{i,j}^+, v_{i,j}^-, w_{i,j}^+, w_{i,j}^- \forall i,j} \Bigg\{\sum_{i=1}^{n_x}\sum_{j=1}^{n_y}\int_0^T \hat H^*(x_i,t, v_{i,j}^+(t), v_{i,j}^-(t), w_{i,j}^+(t), w_{i,j}^-(t)) \\
    &\quad\quad + \frac{\rho^{\ell+1,m}_{i,j}(t)}{2\sigma} \Bigg(\Big(v_{i,j}^+(t) - v_{i,j}^{\ell+1,m,+}(t) - \sigma \frac{(D_x^+ \bar\phi^{\ell+1})_{i,j}(t)}{\rho^{\ell+1,m}_{i,j}(t)} \Big)^2
    + \Big(v_{i,j}^-(t) - v_{i,j}^{\ell+1,m,-}(t) - \sigma \frac{(D_x^- \bar\phi^{\ell+1})_{i,j}(t)}{\rho^{\ell+1,m}_{i,j}(t)} \Big)^2 \\
    &\quad\quad +  \Big(w_{i,j}^+(t) - w_{i,j}^{\ell+1,m,+}(t) - \sigma \frac{(D_y^+ \bar\phi^{\ell+1})_{i,j}(t)}{\rho^{\ell+1,m}_{i,j}(t)} \Big)^2
    + \Big(w_{i,j}^-(t) - w_{i,j}^{\ell+1,m,-}(t) - \sigma \frac{(D_y^- \bar\phi^{\ell+1})_{i,j}(t)}{\rho^{\ell+1,m}_{i,j}(t)} \Big)^2\Bigg)
    dt\Bigg\}.
    \end{split}
    }
 \end{equation}
 Update $\rho_{i,j}$ for all $i=1,\dots, n_x$; $j=1,\dots, n_y$ by 
 \begin{equation}
 {\scriptsize
 \begin{split}
     \rho^{\ell+1,m+1} 
    &= (\max\{\mu_{i,j}^{\ell+1,m+1}, 0\})_{i,j},
    \end{split}
    }
 \end{equation}
 where $\mu_{i,j}^{\ell+1,m+1}$ is a function defined by 
 \begin{equation*}
 {\scriptsize
 \begin{split}
     \mu_{i,j}^{\ell+1,m+1} &= \rho_{i,j}^{\ell+1,m} + \sigma \Big(\dot {\bar\phi}_{i,j}^{\ell+1}(t) + v_{i,j}^{\ell+1,m+1,+}(t)(D_x^+\bar\phi^{\ell+1})_{i,j}(t) + v_{i,j}^{\ell+1,m+1,-}(t)(D_x^-\bar\phi^{\ell+1})_{i,j}(t) \\
    & + w_{i,j}^{\ell+1,m+1,+}(t)(D_y^+\bar\phi^{\ell+1})_{i,j}(t)
    + w_{i,j}^{\ell+1,m+1,-}(t)(D_y^-\bar\phi^{\ell+1})_{i,j}(t) - \hat H^*(x_i,y_j, t, v_{i,j}^{\ell+1,m+1,+}(t), \\
    &v_{i,j}^{\ell+1,m+1,-}(t), w_{i,j}^{\ell+1,m+1,+}(t), w_{i,j}^{\ell+1,m+1,-}(t)) - \epsilon (D_{xx} \bar\phi^{\ell+1} + D_{yy} \bar\phi^{\ell+1})_{i,j}(t)\Big).
 \end{split}
 }
 \end{equation*}
 }
 Set $v^{\ell+1,+} = v^{\ell+1,N_{inner},+}$, $v^{\ell+1,-} = v^{\ell+1,N_{inner},-}$, 
 $w^{\ell+1,+} = w^{\ell+1,N_{inner},+}$, $w^{\ell+1,-} = w^{\ell+1,N_{inner},-}$, $\rho^{\ell+1} = \rho^{\ell+1,N_{inner}}$.
 }
 Return $\phi^{N_{outer}}$.
 \caption{The proposed algorithm for solving~\eqref{eqt:saddle_pt_general_semidisc_2d}\label{alg:pdhg_semi_2d}}
\end{algorithm}

\subsection{Two-dimensional fully-discrete method}\label{appendix:fulldisc_2d}
To solve the fully-discrete equation~\eqref{eqt:disc_HJ_equation_2d}, we propose solving the saddle point problem~\eqref{eqt:saddle_pt_general_fulldisc_2d}, whose first order optimality condition is
\begin{equation}\label{eqt:2pdes_full_discrete_2d}
{\small
\begin{dcases}
(D_t^- \phi)_{i,j,k} + \hat H(x_i,y_j, t_k,(D_x^+ \phi)_{i,j,k}, (D_x^- \phi)_{i,j,k}, (D_y^+ \phi)_{i,j,k}, (D_y^- \phi)_{i,j,k}) \leq \epsilon (D_{xx} \phi)_{i,j,k}+\epsilon (D_{yy} \phi)_{i,j,k}, \\
\quad\quad\quad\quad\quad\quad\quad\quad
\quad\quad\quad\quad\quad\quad\quad\quad
i=1,\dots, n_x;j=1,\dots, n_y; k = 2,\dots, n_t,\\
(D_t^+\rho)_{i,j,k} + (D_x^- (\rho v^+))_{i,j,k} + (D_x^+ (\rho v^-))_{i,j,k} + (D_y^- (\rho w^+))_{i,j,k} + (D_y^+ (\rho w^-))_{i,j,k} + \epsilon (D_{xx} \rho)_{i,j,k} + \epsilon (D_{yy} \rho)_{i,j,k}=0, \\
\quad\quad\quad\quad\quad\quad\quad\quad
\quad\quad\quad\quad\quad\quad\quad\quad
i=1,\dots, n_x;j=1,\dots, n_y; k = 1,\dots, n_t-1,\\
(v_{i,j,k}^+, v_{i,j,k}^-, w_{i,j,k}^+, w_{i,j,k}^-) = \nabla_{(p_1^+, p_1^-, p_2^+, p_2^-)} \hat H(x_i,y_j, t_{k+1},(D_x^+ \phi)_{i,j,k+1}, (D_x^- \phi)_{i,j,k+1}, (D_y^+ \phi)_{i,j,k+1}, (D_y^- \phi)_{i,j,k+1}),\\
\quad\quad\quad\quad\quad\quad\quad\quad
\quad\quad\quad\quad\quad\quad\quad\quad
i=1,\dots, n_x;j=1,\dots, n_y; k = 1,\dots, n_t-1,\\
\phi_{i,j,1} = g(x_i, y_j),\quad  \rho_{i,j,n_t} = c, \quad\quad\,\,\,\, i=1,\dots, n_x;j=1,\dots, n_y,\\
\rho_{i,k-1}\Big((D_t^- \phi)_{i,j,k} + \hat H(x_i,y_j, t_k,(D_x^+ \phi)_{i,j,k}, (D_x^- \phi)_{i,j,k}, (D_y^+ \phi)_{i,j,k}, (D_y^- \phi)_{i,j,k}) - \epsilon (D_{xx} \phi)_{i,j,k}-\epsilon (D_{yy} \phi)_{i,j,k}\Big) = 0, \\
\quad\quad\quad\quad\quad\quad\quad\quad
\quad\quad\quad\quad\quad\quad\quad\quad
i=1,\dots, n_x;j=1,\dots, n_y; k = 2,\dots, n_t.
\end{dcases}
}
\end{equation}
The proposed algorithm for solving~\eqref{eqt:saddle_pt_general_fulldisc_2d} is shown in Algorithm~\ref{alg:pdhg_full_2d}. In the algorithm, we use $\mathcal{L}_{full}$ to denote the objective function in~\eqref{eqt:saddle_pt_general_fulldisc}. When updating $\phi$, the term $(D_t^+\rho)_{i,j,n_t-1}$ denotes $\frac{c - \rho_{i,j,n_t-1}}{\Delta t}$, which is consistent with the terminal condition for $\rho$ in~\eqref{eqt:2pdes_full_discrete_2d}. For simplicity, we use $\alpha$ to denote $(\alpha_{i,j,k})_{i,j,k}$ for any quantity $\alpha$. This algorithm is similar to the one-dimensional case, so we omitted some details. We choose the quadratic function for the penalty term of all functions besides $\phi$. For $\phi$, our choice of penalty term for preconditioning is $\sum_{i=1}^{n_x}\sum_{j=1}^{n_y}\sum_{i=2}^{n_t}\Big(((D_t^- \phi)_{i,j,k} - (D_t^- \phi^\ell)_{i,j,k})^2 +( (D_x^+\phi)_{i,j,k} - (D_x^+\phi^\ell)_{i,j,k})^2 + ( (D_y^+\phi)_{i,j,k} - (D_y^+\phi^\ell)_{i,j,k})^2 \Big)$.

\begin{algorithm}[htbp]
\SetAlgoLined
\SetKwInOut{Input}{Inputs}
\SetKwInOut{Output}{Outputs}
\Input{Stepsize $\tau, \sigma>0$, error tolerance $\delta>0$, inner maximal iteration number $N_{inner}$ and outer maximal iteration number $N_{outer}$.}
\Output{Solution to the saddle point problem~\eqref{eqt:saddle_pt_general_fulldisc_2d}.}
For each $i=1,\dots, n_x$; $j=1,\dots, n_y$, initialize the matrices by $\phi_{i,j,k}^0=g(x_i,y_j)$ for $k=1,\dots, n_t$, $\rho_{i,j,k}^0 = c$, $v_{i,j,k}^{0,+} = v_{i,j,k}^{0,-} = w_{i,j,k}^{0,+} = w_{i,j,k}^{0,-} = 0$ for $k=1,\dots, n_t-1$.

 \For{$\ell = 0,1,\dots,N_{outer}-1$}{
 Update the tensor element $\phi_{i,j,k}$ for $i=1,\dots, n_x$; $j=1,\dots, n_y$; $k=1,\dots, n_t$ by
 \begin{equation}
 {\small
 \begin{split}
 \phi^{\ell+1} 
 &= \phi^\ell + \tau (-D_{tt} - D_{xx} - D_{yy})^{-1}\Big(D_t^+\rho^\ell + D_x^- (\rho^\ell v^{\ell,+}) + D_x^+ (\rho^\ell v^{\ell,-})  \\
 & \quad\quad\quad\quad\quad\quad\quad\quad + D_y^- (\rho^\ell w^{\ell,+}) + D_y^+ (\rho^\ell w^{\ell,-}) +\epsilon D_{xx} \rho^\ell + \epsilon D_{yy} \rho^\ell\Big),
 \end{split}
 }
 \end{equation}
 where $(-D_{tt} - D_{xx} - D_{yy})^{-1}f$ (for a tensor $f$ with elements $f_{i,j,k}$ and the linear operator $D_{tt}f = (\frac{f_{i,j,k-1} -2f_{i,j,k} + f_{i,j,k+1}}{\Delta t^2})_{i,j,k}$) denotes the solution $u$ to the linear system $- (D_{tt} u)_{i,j,k+1} - (D_{xx} u + D_{yy} u)_{i,j,k+1} = f_{i,j,k}$ for all $i=1,\dots, n_x$; $j=1,\dots, n_y$; $k=1,\dots, n_t-1$ with periodic spatial condition, Dirichlet initial condition $u_{i,j,1}=0$, and Neumann terminal condition $u_{i,j,n_t+1}=u_{i,j,n_t}$. 
 
 \If{$\sum_{i=1}^{n_x}\sum_{j=1}^{n_y}\sum_{k=2}^{n_t}|(D_t^- \phi^{\ell+1})_{i,j,k} + \hat H(x_i, y_j, t_k,(D_x^+ \phi^{\ell+1})_{i,j,k}, (D_x^- \phi^{\ell+1})_{i,j,k}, (D_y^+ \phi^{\ell+1})_{i,j,k}, (D_y^- \phi^{\ell+1})_{i,j,k}) - \epsilon (D_{xx} \phi^{\ell+1} + D_{yy} \phi^{\ell+1})_{i,j,k}| \leq \delta$}{
   Return $\phi^{\ell+1}$.
   }
 
 Set $\bar\phi^{\ell+1} = 2\phi^{\ell+1} - \phi^{\ell}$.
 
 Set $v^{\ell+1,0,+} = v^{\ell,+}$, $v^{\ell+1,0,-} = v^{\ell,-}$, $w^{\ell+1,0,+} = w^{\ell,+}$, $w^{\ell+1,0,-} = w^{\ell,-}$, $\rho^{\ell+1,0} = \rho^{\ell}$.

  \For{$m=0,1,\dots, N_{inner}-1$}{
 Update $v_{i,j,k}^+$, $v_{i,j,k}^-$, $w_{i,j,k}^+$ and $w_{i,j,k}^-$ for all $i=1,\dots, n_x$; $j=1,\dots, n_y$; $k=1,\dots, n_t-1$ by
 \begin{equation}
 {\scriptsize
 \begin{split}
    & (v^{\ell+1,m+1,+}, v^{\ell+1,m+1,-}, w^{\ell+1,m+1,+}, w^{\ell+1,m+1,-}) \\
    &= \argmin_{v_{i,j,k}^+, v_{i,j,k}^-, w_{i,j,k}^+, w_{i,j,k}^- \forall i,j,k} \Bigg\{\sum_{i=1}^{n_x}\sum_{j=1}^{n_y}\sum_{k=1}^{n_t-1} \Bigg(\hat H^*(x_i, y_j, t_{k+1},v_{i,j,k}^+, v_{i,j,k}^-, w_{i,j,k}^+, w_{i,j,k}^-) \\
    &\quad\quad\quad\quad + \frac{\rho^{\ell+1,m}_{i,j,k}}{2\sigma} \Bigg(\Big(v_{i,j,k}^+ - v_{i,j,k}^{\ell+1,m,+} - \sigma \frac{(D_x^+ \bar\phi^{\ell+1})_{i,j,k+1}}{\rho^{\ell+1,m}_{i,j,k}} \Big)^2 + \Big(v_{i,j,k}^- - v_{i,j,k}^{\ell+1,m,-} - \sigma \frac{(D_x^- \bar\phi^{\ell+1})_{i,j,k+1}}{\rho^{\ell+1,m}_{i,j,k}} \Big)^2 \\
    &\quad\quad\quad\quad + \Big(w_{i,j,k}^+ - w_{i,j,k}^{\ell+1,m,+} - \sigma \frac{(D_y^+ \bar\phi^{\ell+1})_{i,j,k+1}}{\rho^{\ell+1,m}_{i,j,k}} \Big)^2 + \Big(w_{i,j,k}^- - w_{i,j,k}^{\ell+1,m,-} - \sigma \frac{(D_y^- \bar\phi^{\ell+1})_{i,j,k+1}}{\rho^{\ell+1,m}_{i,j,k}} \Big)^2 \Bigg)\Bigg\}.
    \end{split}
    }
 \end{equation}
 
 Update $\rho_{i,j,k}$ for all $i=1,\dots, n_x$; $j=1,\dots, n_y$; $k=1,\dots, n_t-1$ by 
 \begin{equation}
 {\small
 \begin{split}
     \rho_{i,j,k}^{\ell+1,m+1} 
    &= \max\{\mu_{i,j,k}^{\ell+1,m+1}, 0\},
    \end{split}
    }
 \end{equation}
 where $\mu_{i,j,k}^{\ell+1,m+1}$ is a number defined by
 \begin{equation*}
 {\scriptsize
    \begin{split}
        \mu_{i,j,k}^{\ell+1,m+1} &= \rho_{i,j,k}^{\ell+1,m} + \sigma ((D_t^- \bar\phi)_{i,j,k+1}^{\ell+1} + v_{i,j,k}^{\ell+1,m+1,+}(D_x^+\bar\phi^{\ell+1})_{i,j,k+1} + v_{i,j,k}^{\ell+1,m+1,-}(D_x^-\bar\phi)_{i,j,k+1}^{\ell+1} \\
 &\quad\quad\quad\quad 
 + w_{i,j,k}^{\ell+1,m+1,+}(D_y^+\bar\phi^{\ell+1})_{i,j,k+1} + w_{i,j,k}^{\ell+1,m+1,-}(D_y^-\bar\phi)_{i,j,k+1}^{\ell+1} - \hat H^*(x_i, y_j, t_{k+1},\\
 &\quad\quad\quad\quad  
 v_{i,j,k}^{\ell+1,m+1,+}, v_{i,j,k}^{\ell+1,m+1,-}, 
 w_{i,j,k}^{\ell+1,m+1,+}, w_{i,j,k}^{\ell+1,m+1,-}) 
 - \epsilon (D_{xx} \bar\phi + D_{yy} \bar\phi)_{i,j,k+1}^{\ell+1}).
    \end{split}
    }
 \end{equation*}
 }
 Set $v^{\ell+1,+} = v^{\ell+1,N_{inner},+}$, $v^{\ell+1,-} = v^{\ell+1,N_{inner},-}$,
 $w^{\ell+1,+} = w^{\ell+1,N_{inner},+}$, $w^{\ell+1,-} = w^{\ell+1,N_{inner},-}$,
 $\rho^{\ell+1} = \rho^{\ell+1,N_{inner}}$.
 }
 Return $\phi^{N_{outer}}$.
 \caption{The proposed algorithm for solving~\eqref{eqt:saddle_pt_general_fulldisc_2d}\label{alg:pdhg_full_2d}}
\end{algorithm}

\section{An equivalent formulation}
Through the substitution of the variable $m = \rho v$, we obtain an alternative expression for equation~\eqref{eqt:saddle_cont} as follows:
\begin{equation} \label{eqt:saddle_cont_m}
{\scriptsize
\begin{split}
    \min_{\substack{\phi\\ \phi(x,0)=g(x)}}\max_{\substack{\rho,m\\ \rho \geq 0}}  \int_0^T \int_\Omega \rho(x,t)\left(\frac{\partial \phi(x,t)}{\partial t} - \epsilon \Delta_x \phi(x,t)\right)+ m(x,t)\nabla_x\phi(x,t) - \rho(x,t)H^*\left(x,t, \frac{m(x,t)}{\rho(x,t)}\right)  dxdt - c\int \phi(x,T)dx.
\end{split}
}
\end{equation}

In contrast to the saddle point formulation~\eqref{eqt:saddle_cont}, this revised equation conforms to the conventional PDHG-applicable saddle point formulation (as detailed in~\cite{chambolle2011first}), which thereby ensures convergence for the corresponding algorithm. Analogously, the semi-discrete and fully discrete counterparts discussed in Sections~\ref{appendix:L1_1d} and~\ref{appendix:L1_2d} also follow this structure, thereby endowing their corresponding algorithms with theoretical convergence guarantees. For more details, we refer readers to the following remark and also Remark~\ref{rem:B3_theory}.

\begin{remark}\label{rem:B1_convergencePDHG}
The convergence of the PDHG algorithm for equation~\eqref{eqt:saddle_cont_m}, as well as its associated semi-discrete and fully-discrete counterparts in the specialized context of separable and shifted $1$-homogeneous Hamiltonians described in Appendices~\ref{appendix:L1_1d} and~\ref{appendix:L1_2d}, can be established based on established theory~\cite{chambolle2011first}. These specific saddle point problems adhere to the conventional format outlined in~\cite{chambolle2011first}, characterized by bilinear interactions between primal and dual variables, while other terms exhibit convexity concerning $\phi$ or concavity concerning $(\rho, m)$. To be more precise, in scenarios where $H$ exhibits separability and shifted $1$-homogeneity with respect to $p$, the primal and dual optimization problems in equation~\eqref{eqt:saddle_m_full_discrete_L1} can be reformulated as linear programming problems involving $(\rho, m)$ or $\phi$. Additionally, explicit updating formulas for both primal and dual variables can be derived in this specific case (for more information, refer to Appendices~\ref{appendix:L1_1d} and~\ref{appendix:L1_2d}).

These explicit formulas, which remove the inner loop in Algorithm~\ref{alg:pdhg_cont}, are attainable due to the simplification of the term $\rho(x,t)H^*\left(x,t, \frac{m(x,t)}{\rho(x,t)}\right)$ (and its corresponding nonlinear terms in both semi-discrete and fully-discrete scenarios), resulting in linear objective functions and linear constraints. Nevertheless, in more general cases where the intricate nonlinear term cannot be simplified, there is no feasible approach to derive an explicit updating formula for $(\rho,m)$. In such instances, Algorithm~\ref{alg:pdhg_cont} and its semi-discrete and fully-discrete variations are more suitable.
\end{remark}

We will now provide a summary of the relevant semi-discrete and fully-discrete formulations, along with their corresponding algorithms, for both one-dimensional and two-dimensional cases in the upcoming sections.

\subsection{One-dimensional semi-discrete and fully-discrete methods}\label{appendix:semidisc_disc_1d_m}
For the one-dimensional case, we first discretize the spatial domain $[a,b]$ and get the following semi-discrete saddle point problem
\begin{equation} \label{eqt:saddle_semidiscrete_m_1d}
{\scriptsize
\begin{split}
\min_{\substack{\phi_i \forall i\\ \phi_i(0)=g(x_i)}} \max_{\substack{\rho_i, m_i \forall i\\ \rho_i\geq 0}}  \int_0^T \sum_{i=1}^{n_x}\left(\rho_i(t)\dot \phi_i(t) + m_i(t) (D_x^+\phi)_i(t) - \epsilon\rho_i(t) (D_{xx}\phi)_i(t)\right) - \hat L((x_i)_i, t,(\rho_i(t))_i, (m_i(t))_i) dt - c\sum_{i=1}^{n_x}\phi_i(T),
\end{split}
}
\end{equation}
where $\hat L\colon \R^{n_x}\times [0,+\infty)\times \R^{n_x}\times \R^{n_x}\to\R$ is defined by
\begin{equation*}
\hat L((x_i)_i, t, (\rho_i)_i, (m_i)_i) = \max_{d_i \forall i} \sum_{i=1}^{n_x} (m_id_i - \rho_i \hat H(x_i, t, d_i, d_{i-1})).
\end{equation*}
The derivation of this formula is given as follows
\begin{equation}
{\scriptsize
\begin{split}
&\min_{\phi \text{ satisfying~\eqref{eqt:semidisc_HJ_equation}}} -c \sum_{i=1}^{n_x} \phi_i(T)\\
=& \min_{\substack{\phi_i \forall i\\ \phi_i(0)=g(x_i)}} \max_{\rho_i\forall i} \int_0^T \sum_{i=1}^{n_x} \rho_i(t)\left(\dot\phi_i(t)+ \hat H(x_i, t, (D_x^+\phi)_{i}(t), (D_x^-\phi)_{i}(t)) - \epsilon (D_{xx} \phi)_i(t)\right) dt -c \sum_{i=1}^{n_x} \phi_i(T)\\
=& \min_{\substack{\phi_i, d_i \forall i\\ \phi_i(0)=g(x_i)\\
d_i(t) = (D_x^+\phi)_{i}(t)}} \max_{\rho_i \forall i} \int_0^T \sum_{i=1}^{n_x} \rho_i(t)\left(\dot\phi_i(t)+ \hat H(x_i, t, d_i(t), d_{i-1}(t)) - \epsilon (D_{xx} \phi)_i(t)\right) dt -c \sum_{i=1}^{n_x} \phi_i(T)\\
=& \min_{\substack{\phi_i, d_i \forall i\\ \phi_i(0)=g(x_i)}} \max_{\rho_i , m_i \forall i} \int_0^T \sum_{i=1}^{n_x} \rho_i(t)\left(\dot\phi_i(t)+ \hat H(x_i, t, d_i(t), d_{i-1}(t)) - \epsilon (D_{xx} \phi)_i(t)\right) + m_i(t)((D_x^+\phi)_{i}(t)- d_i(t)) dt -c \sum_{i=1}^{n_x} \phi_i(T)\\
=& \min_{\substack{\phi_i\forall i\\ \phi_i(0)=g(x_i)}} \max_{\rho_i , m_i \forall i} \int_0^T \sum_{i=1}^{n_x} \rho_i(t)\left(\dot\phi_i(t) - \epsilon (D_{xx} \phi)_i(t)\right) + m_i(t)(D_x^+\phi)_{i}(t) dt -c \sum_{i=1}^{n_x} \phi_i(T)\\
&\quad\quad\quad\quad \quad\quad\quad\quad - \int_0^T \max_{d_i(t)} \sum_{i=1}^{n_x}(m_i(t)d_i(t) - \rho_i(t)\hat H(x_i, t, d_i(t), d_{i-1}(t))) dt\\
=& \min_{\substack{\phi_i\forall i\\ \phi_i(0)=g(x_i)}} \max_{\rho_i , m_i \forall i} \int_0^T \sum_{i=1}^{n_x} \rho_i(t)\left(\dot\phi_i(t) - \epsilon (D_{xx} \phi)_i(t)\right) + m_i(t)(D_x^+\phi)_{i}(t) dt -c \sum_{i=1}^{n_x} \phi_i(T) - \int_0^T \hat L((x_i(t))_i, t, (\rho_i(t))_i, (m_i(t))_i) dt.
\end{split}
}
\end{equation}
As we explained in Remark~\ref{rem:relation_mfc} and Appendix~\ref{sec:appendix_semidisc_1d}, when the corresponding mean-field control problem has a solution including a positive density $\rho$, we can add the constraint $\rho \geq 0$ to the saddle point problem, from which we obtain~\eqref{eqt:saddle_semidiscrete_m_1d}.

\begin{remark}\label{rem:B2_separable_semidisc_m_1d}
When $\hat H(x, t,p^+, p^-)$ takes on a ``separable" form, denoted by $\hat H(x, t,p^+, p^-) = \hat H_1(x,t, p^+) + \hat H_2(x, t,p^-)$, the following equality holds:
\begin{equation}
\begin{split}
&\max_{d\in\R^{n_x}} \sum_{i=1}^{n_x} (m_id_i - \rho_i \hat H(x_i,t, d_i, d_{i-1}))
= \max_{d\in\R^{n_x}} \sum_{i=1}^{n_x} ( m_id_i - \rho_i (\hat H_1(x_i,t, d_i)+\hat H_2(x_i,t, d_{i-1}))\\
=& \max_{d\in\R^{n_x}} \sum_{i=1}^{n_x} m_id_i - (\rho_i \hat H_1(x_i,t, d_i)+ \rho_{i+1}\hat H_2(x_{i+1},t, d_{i}))
= \sum_{i=1}^{n_x} (\rho_i \hat H_{1,i,t}+ \rho_{i+1}\hat H_{2,i+1,t})^*(m_i),
\end{split}
\end{equation}
where $\hat H_{1,i,t}$ and $\hat H_{2,i,t}$ denotes the functions $p\mapsto \hat H_{1}(x_i,t, p)$ and $p\mapsto \hat H_{2}(x_i,t, p)$, respectively.
The corresponding saddle point problem becomes
\begin{equation}\label{eqt:remB2_semidisc_m_1d_EO}
{\scriptsize
\begin{split}
\min_{\substack{\phi_i \forall i\\ \phi_i(0)=g(x_i)}} \max_{\substack{\rho_i, m_i \forall i\\ \rho_i\geq 0}}  \int_0^T \sum_{i=1}^{n_x}\left( \rho_i(t)\dot \phi_i(t)+ m_i(t) (D_x^+\phi)_i(t) - \epsilon\rho_i(t) (D_{xx}\phi)_i(t) - (\rho_i(t)\hat H_{1,i,t}+ \rho_{i+1}(t)\hat H_{2,i+1,t})^*(m_i(t))
\right)dt - c\sum_{i=1}^{n_x}\phi_i(T).
\end{split}
}
\end{equation}
For instance, in the context of applying the Engquist-Osher scheme~\cite{osher1988fronts,engquist1980stable,engquist1981one}, the function $\hat H_{1,i,t}$ corresponds to the non-increasing part of $H(x,t,\cdot)$, while $\hat H_{2,i,t}$ corresponds to the non-decreasing part of $H(x,t,\cdot)$. Under this framework, the link between~\eqref{eqt:saddle_pt_general_semidisc} and~\eqref{eqt:remB2_semidisc_m_1d_EO} involves the transformation of variables:
$\rho_iv_i^+ = \min\{m_i, 0\}$ and $\rho_iv_i^- = \max\{m_{i-1},0\}$.
Furthermore, if $H$ exhibits 1-homogeneity and convexity with respect to $p$, the function $(\rho_i(t)\hat H_{1,i,t}+ \rho_{i+1}(t)\hat H_{2,i+1,t})^*$ becomes associated with the indicator function, which notably simplifies the computation process. For further details, please refer to Appendix~\ref{appendix:L1_1d}.
\end{remark}

Next, we apply implicit time discretization and employ $(D_t^- \phi)_{i,j,k}$ to estimate $\dot \phi_{i,j}(t)$. Consequently, the saddle point formulation~\eqref{eqt:saddle_semidiscrete_m_2d} transforms into the following expression (It's important to observe that the index $k$ for $\rho$, $m^1$, and $m^2$ ranges from $1$ to $n_t-1$, while for $\phi$, it spans from $1$ to $n_t$.)
\begin{equation}\label{eqt:saddle_m_full_discrete_1d}
{\scriptsize
\begin{split}
\min_{\substack{\phi_{i,k} \forall i,k\\ \phi_{i,1}=g(x_i)}} \max_{\substack{\rho_{i,k}, m_{i,k} \forall i,k\\ \rho_{i,k}\geq 0}}  
\sum_{i=1}^{n_x} \sum_{k=1}^{n_t-1}\left(\rho_{i,k}(D_t^-\phi)_{i,k+1} + m_{i,k} (D_x^+\phi)_{i,k+1} - \epsilon\rho_{i,k} (D_{xx}\phi)_{i,k+1}\right) - \hat L((x_i)_i,(t_k)_k, (\rho_{i,k})_{i,k}, (m_{i,k})_{i,k}) - \frac{c}{\Delta t}\sum_{i=1}^{n_x}\phi_{i,n_t},
\end{split}
}
\end{equation}
where we abuse the notation $\hat L$ and define $\hat L\colon \R^{n_x}\times \R^{n_t}\times \R^{n_x\times (n_t-1)}\times \R^{n_x \times (n_t-1)}\to\R$ by
\begin{equation*}
\hat L((x_i)_i, (t_k)_k, (\rho_{i,k})_{i,k}, (m_{i,k})_{i,k}) = \max_{d_{i,k}\in\R\forall i,k} \sum_{k=1}^{n_t-1} \sum_{i=1}^{n_x} (m_{i,k}d_{i,k} - \rho_{i,k} \hat H(x_i, t_{k+1},d_{i,k}, d_{i-1,k})).
\end{equation*}
The discussion in Remark~\ref{rem:B2_separable_semidisc_m_1d} also applies to this fully-discrete formula with straightforward modifications. 

While it is possible to decouple the formula for $\hat L$ with respect to $k$, this decoupling cannot be achieved with respect to $i$.
This discrepancy constitutes a significant limitation of this formulation when compared to~\eqref{eqt:saddle_pt_general_fulldisc}, wherein the nonlinear term within the objective function can be effortlessly decoupled concerning both $i$ and $k$. This property facilitates the parallel computation of the updates for $\rho, v^+, v^-$ in Algorithm~\ref{alg:pdhg_full_1d}. 
However, in certain special scenarios (such as when $H$ is separable and shifted 1-homogeneous with respect to $p$, as discussed in the subsequent section), the function $\hat L$ takes on a particular structure, enabling the derivation of an explicit joint updating formula for $(\rho,m)$.

\subsection{One-dimensional special case: 1-homogeneous Hamiltonian}\label{appendix:L1_1d}
In this section, we consider the one-dimensional HJ PDE whose Hamiltonian $H$ is in the form of
\begin{equation} \label{eqt:H_L1_1d}
H(x,t,p) = \cinh(x,t)|p| + f(x,t),
\end{equation}
with the assumption that $\cinh(x,t) > 0$ for any $x\in\R$, $t\in [0,T]$.
We apply the Engquist-Osher scheme~\cite{osher1988fronts,engquist1980stable,engquist1981one}, in which we have $\hat H(x, t,p_+, p_-) = \hat H_-(x, t,p_+) + \hat H_+(x,t, p_-) + f(x,t)$, where $\hat H_-(x, t,p) = \cinh(x,t)\max\{-p,0\}$ and $\hat H_+(x, t,p) = \cinh(x,t)\max\{p,0\}$. 
(Notably, these functions $\hat H_+$ and $\hat H_+$ draw parallels to the ReLU activation function within the domain of neural networks.) In this scenario, we arrive at the following formula for the function $\hat L$ in~\eqref{eqt:saddle_m_full_discrete_1d}:
\begin{equation}
{\small
\begin{split}
\hat L((x_i)_i,(t_k)_k, (\rho_{i,k})_{i,k}, (m_{i,k})_{i,k}) 
&= \max_{d_{i,k}\in\R\forall i,k}  \sum_{i=1}^{n_x} \sum_{k=1}^{n_t-1} (m_{i,k}d_{i,k} - \rho_{i,k} \hat H(x_i, t_{k+1}, d_{i,k}, d_{i-1,k}))\\
&= \sum_{i=1}^{n_x} \sum_{k=1}^{n_t-1} (-\rho_{i,k}f(x_i,t_{k+1}) + \ind_{[-\rho_{i,k}\cinh(x_i,t_{k+1}), \rho_{i+1,k}\cinh(x_{i+1},t_{k+1})]}(m_{i,k})),
\end{split}
}
\end{equation}
where $\ind_{C}$ represents the indicator function of set $C$, which is $0$ when the variable is within $C$, and $+\infty$ otherwise.
Then, the saddle point problem becomes
\begin{equation}\label{eqt:saddle_m_full_discrete_L1}
{\small
\begin{split}
\min_{\substack{\phi_{i,k} \forall i,k\\ \phi_{i,1}=g(x_i)}} \max_{\substack{\rho_{i,k}, m_{i,k} \forall i,k\\ \rho_{i,k}\geq 0}} \Bigg\{\sum_{i=1}^{n_x}\sum_{k=1}^{n_t-1} \left( \rho_{i,k} (D_t^-\phi)_{i,k+1}+ m_{i,k} (D_x^+\phi)_{i,k+1} - \epsilon\rho_{i,k} (D_{xx}\phi)_{i,k+1} - \rho_{i,k}f(x_i,t_{k+1})
\right) - \frac{c}{\Delta t}\sum_{i=1}^{n_x}\phi_{i,n_t}
\colon\\
  -\rho_{i,k}\cinh(x_i,t_{k+1})\leq m_{i,k}\leq  \rho_{i+1,k}\cinh(x_{i+1},t_{k+1}), \forall i=1,\dots, n_x; k=1,\dots, n_t-1\Bigg\}.
\end{split}
}
\end{equation}
As elucidated in Remark~\ref{rem:B1_convergencePDHG}, the objective function within the saddle point problem~\eqref{eqt:saddle_m_full_discrete_L1} adheres to the conventional form within the theory in~\cite{chambolle2011first}. Furthermore, the sole nonlinear term in this context is the indicator function, resulting in explicit formulas for updating $\phi$, $\rho$, and $m$. 
It's important to note that this approach can be extended to encompass any Hamiltonian that is convex and shifted $1$-homogeneous with respect to~$p$.

We employ the PDHG method to tackle~\eqref{eqt:saddle_m_full_discrete_L1}. A comprehensive outline of the algorithm is encapsulated in Algorithm~\ref{alg:pdhg_full_m_1d_L1}. In the algorithm, we conveniently employ the notation $\mathcal{L}_{full}$ to represent the objective function in~\eqref{eqt:saddle_m_full_discrete_L1}. During the update of $\phi$, the component $(D_t^+\rho)_{i,n_t-1}$ represents $\frac{c - \rho_{i,n_t-1}}{\Delta t}$, aligning with the terminal condition for $\rho$.
As we possess an explicit formula for the joint update of $(\rho,m)$, there is no need for an inner loop. This method is supported by a theoretical guarantee, as explained in Remark~\ref{rem:B1_convergencePDHG}, and this is further detailed in the following remark.

\begin{remark}\label{rem:B3_theory}
According to the theory in~\cite{chambolle2011first}, the algorithm converges when the chosen step sizes $\tau$ and $\sigma$ adhere to the condition $\tau \sigma \|K\|^2 < 1$, where $K$ represents the bilinear operator. In our case, we have
\begin{equation*}
\begin{split}
    |K((\rho,m), \phi)| &= \left|\sum_{i=1}^{n_x}\sum_{k=1}^{n_t-1} \left( \rho_{i,k} (D_t^-\phi)_{i,k+1}+ m_{i,k} (D_x^+\phi)_{i,k+1} - \epsilon\rho_{i,k} (D_{xx}\phi)_{i,k+1}
\right)\right|\\
&\leq \|\rho\|\|D_t^-\phi\| + \|m\| \|D_x^+\phi\| + \frac{2\epsilon}{\Delta x} \|\rho\|\|D_x^+\phi\|\\
&\leq 2\left(1+ \frac{2\epsilon}{\Delta x}\right) \|(\rho, m)\| \|(D_t^-\phi, D_x^+\phi)\|.
\end{split}
\end{equation*}
Consequently, we derive $\|K\|\leq 2\left(1+ \frac{2\epsilon}{\Delta x}\right)$, allowing for the selection of step sizes $\tau$ and $\sigma$ that satisfy $\tau \sigma< (2 + \frac{4\epsilon}{\Delta x})^{-2}$.
\end{remark}

\begin{remark}
As elaborated in Remark~\ref{rem:B2_separable_semidisc_m_1d}, the connection between~\eqref{eqt:saddle_pt_general_fulldisc} and~\eqref{eqt:saddle_m_full_discrete_1d} involves the variable transformations $\rho_{i,k}v_{i,k}^+ = \min\{m_{i,k}, 0\}$ and $v_{i,k}^- = \max\{m_{i-1,k},0\}$. Algorithm~\ref{alg:pdhg_full_m_1d_L1} is closely related to Algorithm~\ref{alg:pdhg_full_1d}, though there exist some distinctions.

The penalties pertaining to $\rho$ in Algorithms~\ref{alg:pdhg_full_m_1d_L1} and~\ref{alg:pdhg_full_1d} are identical. The primary focus is directed towards the penalty related to $m$, which is $\sum_{i,k}(m_{i,k} - m_{i,k}^\ell)^2$. If both $m_{i,k}$ and $m^\ell_{i,k}$ are positive, the $(i,k)$-th term in the  summation becomes $(\rho_{i,k}v_{i,k}^- - \rho^\ell_{i,k}v_{i,k}^{\ell,-})^2$, closely resembling the penalty term $(\rho^\ell_{i,k}v_{i,k}^- - \rho^\ell_{i,k}v_{i,k}^{\ell,-})^2$ in Algorithm~\ref{alg:pdhg_full_1d}. Similarly, when both $m_{i,k}$ and $m^\ell_{i,k}$ are negative, $(m_{i,k} - m_{i,k}^\ell)^2$ becomes $(\rho_{i,k}v_{i,k}^+ - \rho^\ell_{i,k}v_{i,k}^{\ell,+})^2$, analogous to the penalty term for $v^+$ in Algorithm~\ref{alg:pdhg_full_1d}.

Therefore, if the sign of $m$ remains constant between the previous and current steps, Algorithms~\ref{alg:pdhg_full_m_1d_L1} and~\ref{alg:pdhg_full_1d} exhibit similarities.
\end{remark}

\begin{algorithm}[htbp]
\SetAlgoLined
\SetKwInOut{Input}{Inputs}
\SetKwInOut{Output}{Outputs}
\Input{Stepsize $\tau, \sigma>0$, error tolerance $\delta>0$, and the maximal iteration number $N$.}
\Output{Solution to the saddle point problem~\eqref{eqt:saddle_m_full_discrete_L1}.}
For each $i=1,\dots, n_x$, initialize the matrices by $\phi_{i,k}^0=g(x_i)$ for $k=1,\dots, n_t$, $\rho_{i,k}^0 = c$, $m_{i,k}^{0} = 0$ for $k=1,\dots, n_t-1$.

 \For{$\ell = 0,1,\dots,N-1$}{
 Update the matrix $\phi_{i,k}$ for $i=1,\dots, n_x$; $k=1,\dots, n_t$ by
 \begin{equation}
 {\small
 \begin{split}
 (\phi_{i,k}^{\ell+1})_{i,k} &= \argmin_{\substack{\phi_{i,k} \forall i,k \\ \phi_{i,1}=g(x_i)}} \mathcal{L}_{full}((\phi_{i,k})_{i,k}, (\rho^\ell_{i,k})_{i,k}, (v^{\ell,+}_{i,k})_{i,k}, (v^{\ell,-}_{i,k})_{i,k})  \\
 & \quad\quad\quad\quad
 + \frac{1}{2\tau} \sum_{i=1}^{n_x}\sum_{k=2}^{n_t}\left(((D_t^- \phi)_{i,k} - (D_t^- \phi^\ell)_{i,k})^2 +( (D_x^+\phi)_{i,k} - (D_x^+\phi^\ell)_{i,k})^2 \right) \\
 &= (\phi_{i,k}^\ell)_{i,k} + \tau (-D_{tt} - D_{xx})^{-1}\left(D_t^+\rho^\ell + D_x^- m^\ell + \epsilon D_{xx} \rho^\ell\right),
 \end{split}
 }
 \end{equation}
 where $(-D_{tt} - D_{xx})^{-1}(f_{i,k})_{i,k}$ (for a matrix $f$ with elements $f_{i,k}$ and the linear operator $D_{tt}f = (\frac{f_{i,k-1} -2f_{i,k} + f_{i,k+1}}{\Delta t^2})_{i,k}$) denotes the solution $u$ to the linear system $- (D_{tt} u)_{i,k+1} - (D_{xx} u)_{i,k+1} = f_{i,k}$ for all $i=1,\dots, n_x$; $k=1,\dots, n_t-1$ with periodic spatial condition, Dirichlet initial condition $u_{i,1}=0$, and Neumann terminal condition $u_{i,n_t+1}=u_{i,n_t}$.
 
 \If{$\sum_{i=1}^{n_x}\sum_{k=2}^{n_t}|(D_t^- \phi^{\ell+1})_{i,k} + \hat H(x_i, t_{k}, (D_x^+ \phi^{\ell+1})_{i,k}, (D_x^- \phi^{\ell+1})_{i,k}) - \epsilon (D_{xx} \phi^{\ell+1})_{i,k}| \leq \delta$}{
   Return $(\phi_{i,k}^{\ell+1})_{i,k}$.
   }
 
 Set $\bar\phi_{i,k}^{\ell+1} = 2\phi_{i,k}^{\ell+1} - \phi_{i,k}^{\ell}$ for all $i=1,\dots, n_x$; $k=2,\dots, n_t$.
 
Update $\rho_{i,k}$ for all $i=1,\dots, n_x$; $k=1,\dots, n_t-1$ by
\begin{equation*}
{\scriptsize
\begin{split}
    \rho^{\ell+1}_{i,k} = \argmin_{y\in \mathcal{A}_{i,k}^\ell}
\Bigg\{(y - \alpha_{i,k}^\ell)^2 + \chara_{\{z_{i,k}^\ell\leq 0, \, -y\cinh(x_i, t_{k+1}) \geq z_{i,k}^\ell\}}\left(y\cinh(x_i, t_{k+1}) + z_{i,k}^\ell\right)^2\\
+ \chara_{\{z_{i-1,k}^\ell > 0, \,y\cinh(x_{i}, t_{k+1}) \leq z_{i-1,k}^\ell\}}\left(y\cinh(x_{i}, t_{k+1}) - z_{i-1,k}^\ell\right)^2\Bigg\},
\end{split}
}
\end{equation*}
where $z^\ell_{i,k} = m^\ell_{i,k} + \sigma (D_x^+\bar\phi^{\ell+1})_{i, k+1}$, 
$\alpha_{i,k}^\ell = \rho^{\ell}_{i,k} + \sigma ((D_t^-\bar\phi^{\ell+1})_{i,k+1} + f(x_i, t_{k+1}) - \epsilon (D_{xx}\bar\phi^{\ell+1})_{i, k+1})$, and $\delta_C$ is a function which takes the value $1$ if the constraint in $C$ is satisfied, and $0$ otherwise.
This optimization problem can be solved by comparing the function values of several candidates in $\mathcal{A}_{i,k}^\ell$ defined in~\eqref{eqt:udpating_rho_1d_L1_2}. For details, see Section~\ref{sec:update_rho_m_L1_1d}.

Update $m_{i,k}$ for all $i=1,\dots, n_x$; $k=1,\dots, n_t-1$ by
\begin{equation*}
m^{\ell+1}_{i,k} = \max\left\{\min\left\{z^\ell_{i,k}, \rho^{\ell+1}_{i+1,k}\cinh(x_{i+1}, t_{k+1})\right\}, -\rho^{\ell+1}_{i,k}\cinh(x_i, t_{k+1})\right\},
\end{equation*}
where $z_{i,k}^\ell$ is defined in the last step.
}

 Return $(\phi_{i,k}^{N})_{i,k}$.
 \caption{The proposed algorithm for solving~\eqref{eqt:saddle_m_full_discrete_L1}\label{alg:pdhg_full_m_1d_L1}}
\end{algorithm}

The update procedure for $\phi$ remains as previously described. We will now detail the process of deriving updates for $\rho$ and $m$.

\subsubsection{The detailed derivation of updating $\rho$ and $m$ in Algorithm~\ref{alg:pdhg_full_m_1d_L1}} \label{sec:update_rho_m_L1_1d}

When $H$ adheres to the conditions specified in~\eqref{eqt:H_L1_1d}, the updates for $\rho$ and $m$ require solving the subsequent problem:
\begin{equation}\label{eqt:secB21_updating_rho_m_L1_1d}
{\scriptsize
\begin{split}
\min_{\substack{\rho_{i,k}, m_{i,k} \forall i,k\\ \rho_{i,k}\geq 0}} \Bigg\{-\sum_{i=1}^{n_x}\sum_{k=1}^{n_t-1} \left( \rho_{i,k} (D_t^-\bar\phi^{\ell+1})_{i,k+1}+ m_{i,k} (D_x^+\bar\phi^{\ell+1})_{i,k+1} - \epsilon\rho_{i,k} (D_{xx}\bar\phi^{\ell+1})_{i,k+1} + \rho_{i,k}f(x_i,t_{k+1})
\right) + \frac{1}{2\sigma}\sum_{k=1}^{n_t-1}\sum_{i=1}^{n_x}(\rho_{i,k} - \rho_{i,k}^\ell)^2\\
+ \frac{1}{2\sigma}\sum_{k=1}^{n_t-1}\sum_{i=1}^{n_x} (m_{i,k} - m_{i,k}^\ell)^2
\colon
  -\rho_{i,k}\cinh(x_i,t_{k+1})\leq m_{i,k}\leq  \rho_{i+1,k}\cinh(x_{i+1},t_{k+1}), \forall i=1,\dots, n_x; k=1,\dots, n_t-1\Bigg\}.
\end{split}
}
\end{equation}
Define $z$ and $\alpha$ as in Algorithm~\ref{alg:pdhg_full_m_1d_L1}, i.e.,
$z^\ell_{i,k} = m^\ell_{i,k} + \sigma (D_x^+\bar\phi^{\ell+1})_{i, k+1}$, and 
$\alpha_{i,k}^\ell = \rho^{\ell}_{i,k} + \sigma ((D_t^-\bar\phi^{\ell+1})_{i,k+1} + f(x_i, t_{k+1}) - \epsilon (D_{xx}\bar\phi^{\ell+1})_{i, k+1})$. Then, we can rewrite~\eqref{eqt:secB21_updating_rho_m_L1_1d} as
\begin{equation}\label{eqt:secB21_updating_rho_m_L1_1d_simplified}
{\scriptsize
\begin{split}
\min_{\substack{\rho_{i,k}, m_{i,k} \forall i,k\\ \rho_{i,k}\geq 0}} \Bigg\{\frac{1}{2\sigma}\sum_{k=1}^{n_t-1}\sum_{i=1}^{n_x}\left((\rho_{i,k} - \alpha_{i,k}^\ell)^2
+ (m_{i,k} - z_{i,k}^\ell)^2\right)
\colon
  -\rho_{i,k}\cinh(x_i,t_{k+1})\leq m_{i,k}\leq  \rho_{i+1,k}\cinh(x_{i+1},t_{k+1}), \forall i=1,\dots, n_x; k=1,\dots, n_t-1\Bigg\}.
\end{split}
}
\end{equation}

For a fixed $\rho$, the minimization for $m$ is separable and in the following form
\begin{equation} 
\begin{split}
m^{\ell+1}_{i,k}
& = \argmin_{y\in\R} \left\{(y - z_{i,k}^\ell)^2 \colon y\in [-\rho_{i,k}\cinh(x_i, t_{k+1}),\rho_{i+1,k}\cinh(x_{i+1}, t_{k+1})]\right\},
\end{split}
\end{equation}
whose solution is given by the truncation of $z_{i,k}^\ell$ in the interval $[-\rho_{i,k}\cinh(x_i, t_{k+1}), \rho_{i+1, k}\cinh(x_{i+1}, t_{k+1})]$:
\begin{equation}\label{eqt:update_m}
m^{\ell+1}_{i,k} = \max\{\min\{z_{i,k}^\ell, \rho_{i+1, k}\cinh(x_{i+1}, t_{k+1})\}, -\rho_{i,k}\cinh(x_i, t_{k+1})\}.
\end{equation}

Note that the right hand side of~\eqref{eqt:update_m} can be seen as a function of $\rho_{i,k}$ and $\rho_{i+1,k}$. After some computation, we get 
\begin{itemize}
    \item If $z_{i,k}^\ell >0$, we have $m^{\ell+1}_{i,k} = \min\{z_{i,k}^\ell, \rho_{i+1,k}\cinh(x_{i+1}, t_{k+1})\}$, and hence $m^{\ell+1}_{i,k}$ only depends on $\rho_{i+1,k}$ in the following form (denoting this function by $G_{i,k}^+$)
    \begin{equation}
m^{\ell+1}_{i,k} = G_{i,k}^+(\rho_{i+1,k}) := \begin{dcases}
\rho_{i+1,k}\cinh(x_{i+1}, t_{k+1}), & \text{if }\rho_{i+1}\cinh(x_{i+1}, t_{k+1}) \leq z_{i,k}^\ell,\\
z_{i,k}^\ell, & \text{if }\rho_{i+1}\cinh(x_{i+1}, t_{k+1}) > z_{i,k}^\ell.
\end{dcases}
\end{equation}

\item If $z_{i,k}^\ell \leq 0$, we have $m^{\ell+1}_{i,k} = \max\{z_{i,k}^\ell, -\rho_{i,k}\cinh(x_i, t_{k+1})\}$, and hence $m^{\ell+1}_{i,k}$ only depends on $\rho_{i,k}$ in the following form (denoting this function by $G_{i,k}^-$)
    \begin{equation}
m^{\ell+1}_{i,k} = G_{i,k}^-(\rho_{i,k}) := \begin{dcases}
-\rho_{i,k}\cinh(x_i, t_{k+1}), & \text{if }-\rho_{i,k}\cinh(x_i, t_{k+1}) \geq z_{i,k}^\ell,\\
z_{i,k}^\ell, & \text{if } -\rho_{i,k}\cinh(x_i, t_{k+1}) < z_{i,k}^\ell.
\end{dcases}
\end{equation}
\end{itemize}
In other words, depending on the sign of $z_{i,k}^\ell$, $m^{\ell+1}_{i,k}$ is a piecewise linear (ReLU-shape) function which depends either on $\rho_{i,k}$ or on $\rho_{i+1,k}$.

We plug these formulas into~\eqref{eqt:secB21_updating_rho_m_L1_1d_simplified} to get a minimization problem of $\rho$ only:
\begin{equation}
{\small
\begin{split}
\rho^{\ell+1}_{i,k}
&= \argmin_{y\geq 0} 
(y - \alpha_{i,k}^\ell)^2
+ \chara_{\{z_{i,k}^\ell\leq 0\}}(G_{i,k}^-(y) - z_{i,k}^\ell)^2 + \chara_{\{z_{i-1,k}^\ell > 0\}} (G_{i-1,k}^+(y) - z_{i-1,k}^\ell)^2
\\
& = \argmin_{y\geq 0}
(y - \alpha_{i,k}^\ell)^2 + \chara_{\{z_{i,k}^\ell\leq 0, \, -y\cinh(x_i, t_{k+1}) \geq z_{i,k}^\ell\}}\left(y\cinh(x_i, t_{k+1}) + z_{i,k}^\ell\right)^2\\
&\quad\quad\quad\quad\quad\quad\quad\quad 
+ \chara_{\{z_{i-1,k}^\ell > 0, \,y\cinh(x_{i-1}, t_{k+1}) \leq z_{i-1,k}^\ell\}}\left(y\cinh(x_{i-1}, t_{k+1}) - z_{i-1,k}^\ell\right)^2.
\end{split}
}
\end{equation}
Consequently, the objective function takes on the form of a piecewise quadratic function. With each $\delta$-function capable of taking only values of $1$ or $0$, there exist four potential combinations. We consider the minimization of each potential combination as a candidate. These candidates encompass $\alpha_{i,k}^\ell$, $\frac{\alpha_{i,k}^\ell -\cinh(x_i,t_{k+1}) z_{i,k}^\ell}{\cinh(x_i,t_{k+1})^2+1}$, $\frac{\alpha_{i,k}^\ell +\cinh(x_{i-1},t_{k+1}) z_{i-1,k}^\ell}{\cinh(x_{i-1}, t_{k+1})^2+1}$, and $\frac{\alpha_{i,k}^\ell -\cinh(x_i,t_{k+1}) z_{i,k}^\ell +\cinh(x_{i-1},t_{k+1}) z_{i-1,k}^\ell}{\cinh(x_i,t_{k+1})^2 + \cinh(x_{i-1}, t_{k+1})^2+1}$. 
We also incorporate $-\frac{z_{i,k}^\ell}{\cinh(x_i, t_{k+1})}$ and $\frac{z_{i-1,k}^\ell}{\cinh(x_{i-1}, t_{k+1})}$ into the candidates, which correspond to the points at which the $\delta$-functions transition between values.
Furthermore, considering the constraint $\rho \geq 0$, we truncate all these candidates to the interval $[0,+\infty)$ and additionally include $0$ as another candidate.
In essence, we update $\rho$ employing the subsequent formula:
\begin{equation}\label{eqt:udpating_rho_1d_L1}
{\scriptsize
\begin{split}
\rho^{\ell+1}_{i,k}
= \argmin_{y\in \mathcal{A}_{i,k}^\ell} 
(y - \alpha_{i,k}^\ell)^2 + \chara_{\{z_{i,k}^\ell\leq 0, \, -y\cinh(x_i, t_{k+1}) \geq z_{i,k}^\ell\}}\left(y\cinh(x_i, t_{k+1}) + z_{i,k}^\ell\right)^2\\
+ \chara_{\{z_{i-1,k}^\ell > 0, \,y\cinh(x_{i-1}, t_{k+1}) \leq z_{i-1,k}^\ell\}}\left(y\cinh(x_{i-1}, t_{k+1}) - z_{i-1,k}^\ell\right)^2,
\end{split}
}
\end{equation}
where $\mathcal{A}_{i,k}^\ell$ designates the set encompassing all the candidates, defined by
\begin{equation}\label{eqt:udpating_rho_1d_L1_2}
{\scriptsize
\begin{split}
\mathcal{A}_{i,k}^\ell = \Bigg\{0, \left(\alpha_{i,k}^\ell\right)_+, 
\left(\frac{\alpha_{i,k}^\ell -\cinh(x_i,t_{k+1}) z_{i,k}^\ell}{\cinh(x_i,t_{k+1})^2+1}\right)_+, 
\left(\frac{\alpha_{i,k}^\ell +\cinh(x_{i-1},t_{k+1}) z_{i-1,k}^\ell}{\cinh(x_{i-1}, t_{k+1})^2+1}\right)_+, \left(\frac{\alpha_{i,k}^\ell -\cinh(x_i,t_{k+1}) z_{i,k}^\ell +\cinh(x_{i-1},t_{k+1}) z_{i-1,k}^\ell}{\cinh(x_i,t_{k+1})^2 + \cinh(x_{i-1}, t_{k+1})^2+1}\right)_+,\\
\left(-\frac{z_{i,k}^\ell}{\cinh(x_i, t_{k+1})}\right)_+, \left(\frac{z_{i-1,k}^\ell}{\cinh(x_{i-1}, t_{k+1})}\right)_+\Bigg\}.
\end{split}
}
\end{equation}
Here, $y_+$ denotes the positive part of the number $y$, denoted as $y_+ := \max\{y,0\}$.

\subsection{Two-dimensional semi-discrete and fully-discrete method}
\label{appendix:semidisc_disc_2d_m}
Within this section, we present the semi-discrete and fully-discrete counterparts of the two-dimensional version of~\eqref{eqt:saddle_cont_m}. The derivations closely mirror those in the one-dimensional cases, and as such, we have omitted certain particulars. Initially, we proceed with spatial domain discretization, leading to the following formulation:
\begin{equation} \label{eqt:saddle_semidiscrete_m_2d}
{\scriptsize
\begin{split}
\min_{\substack{\phi_{i,j} \forall i,j\\ \phi_{i,j}(0)=g(x_i, y_j)}} \max_{\substack{\rho_{i,j}, m_{i,j}^1, m_{i,j}^2 \forall i,j\\ \rho_{i,j}\geq 0}}  \int_0^T \sum_{i=1}^{n_x}\sum_{j=1}^{n_y}\left(\rho_{i,j}(t)\dot \phi_{i,j}(t) + m^1_{i,j}(t) (D_x^+\phi)_{i,j}(t) + m^2_{i,j}(t) (D_y^+\phi)_{i,j}(t) - \epsilon\rho_{i,j}(t) (D_{xx}\phi + D_{yy}\phi)_{i,j}(t)\right) \\
- \hat L\left((x_{i,j})_{i,j}, t, (\rho_{i,j}(t))_{i,j}, (m_{i,j}^1(t))_{i,j}, (m_{i,j}^2(t))_{i,j}\right) dt - c\sum_{i=1}^{n_x} \sum_{j=1}^{n_y}\phi_{i,j}(T).
\end{split}
}
\end{equation}
where $\hat L\colon \R^{n_x}\times \R^{n_y}\times [0,+\infty) \times \R^{n_x\times n_y}\times \R^{n_x\times n_y}\times \R^{n_x\times n_y}\to\R$ is defined by
\begin{equation*}
{\scriptsize
\begin{split}
\hat L\left((x_{i})_{i}, (y_{j})_{j}, t, (\rho_{i,j})_{i,j}, (m_{i,j}^1)_{i,j}, (m_{i,j}^2)_{i,j}\right) = \max_{d_{i,j}^1, d^2_{i,j} \forall i,j} \sum_{i=1}^{n_x} \sum_{j=1}^{n_y}\left(m_{i,j}^1d_{i,j}^1 + m_{i,j}^2d_{i,j}^2 - \rho_{i,j} \hat H(x_{i}, y_{j},t, d_{i,j}^1, d_{i-1,j}^1, d_{i,j}^2, d_{i,j-1}^2)\right).
\end{split}
}
\end{equation*}
The derivation of this formula is similar as the one-dimensional case, and hence omitted here.

\begin{remark}\label{rem:separable_semidisc_2d}
In the case when $\hat H(x, y, t, p^{1,+}, p^{1,-}, p^{2,+}, p^{2,-})$ is ``separable'', i.e., it can be written as $\hat H_1^1(x,y, t, p^{1,+}) + \hat H_2^1(x,y, t, p^{1,-}) + \hat H_1^2(x,y, t, p^{2,+}) + \hat H_2^2(x,y, t, p^{2,-})$, we have
\begin{equation}
{\scriptsize
\begin{split}
&\max_{d_{i,j}^1, d_{i,j}^2\forall i,j} \sum_{i=1}^{n_x}\sum_{j=1}^{n_y} (m_{i,j}^1d_{i,j}^1 + m_{i,j}^2d_{i,j}^2 - \rho_{i,j} \hat H(x_i,y_j,t, d_{i,j}^1, d_{i-1,j}^1, d_{i,j}^2, d_{i,j-1}^2))\\
= &\max_{d_{i,j}^1, d_{i,j}^2\forall i,j} \sum_{i=1}^{n_x}\sum_{j=1}^{n_y} (m_{i,j}^1d_{i,j}^1 + m_{i,j}^2d_{i,j}^2 - \rho_{i,j} \hat H_1^1(x_i,y_j, t, d_{i,j}^1) - \rho_{i+1,j} \hat H_2^1(x_{i+1},y_j, t, d_{i,j}^1) - \rho_{i,j}\hat H_1^2(x_i,y_j, t, d_{i,j}^2) - \rho_{i,j+1}\hat H_2^2(x_i,y_{j+1}, t, d_{i,j}^2))\\
=& \sum_{i=1}^{n_x}\sum_{j=1}^{n_y} \left((\rho_{i,j}\hat H_{1,i,j,t}^1 + \rho_{i+1,j}\hat H_{2,i+1,j,t}^1)^*(m_{i,j}^1) + (\rho_{i,j}\hat H_{1,i,j,t}^2 + \rho_{i,j+1}\hat H_{2,i,j+1,t}^2)^*(m_{i,j}^2)\right),
\end{split}
}
\end{equation}
where $\hat H_{a,i,j,t}^b$ denotes the functions $p\mapsto \hat H_{a}^b(x_i,y_j, t, p)$.
The corresponding saddle point problem is
\begin{equation}
{\scriptsize
\begin{split}
\min_{\substack{\phi_{i,j} \forall i,j\\ \phi_{i,j}(0)=g(x_i,y_j)}} \max_{\substack{\rho_{i,j}, m_{i,j}^1, m_{i,j}^2 \forall i,j\\ \rho_{i,j}\geq 0}}  \int_0^T \sum_{i=1}^{n_x}\sum_{j=1}^{n_y}\Big(\rho_{i,j}(t)\dot \phi_{i,j}(t) + m^1_{i,j}(t) (D_x^+\phi)_{i,j}(t) + m^2_{i,j}(t) (D_y^+\phi)_{i,j}(t) - \epsilon\rho_{i,j}(t) (D_{xx}\phi + D_{yy}\phi)_{i,j}(t) \\
- \left((\rho_{i,j}(t)\hat H_{1,i,j,t}^1 + \rho_{i+1,j}(t)\hat H_{2,i+1,j,t}^1)^*(m_{i,j}^1(t)) + (\rho_{i,j}\hat H_{1,i,j,t}^2 + \rho_{i,j+1}(t)\hat H_{2,i,j+1,t}^2)^*(m_{i,j}^2(t))\right)
\Big)dt - c\sum_{i=1}^{n_x}\sum_{j=1}^{n_y}\phi_{i,j}(T).
\end{split}
}
\end{equation}
\end{remark}

Now, we perform the implicit time discretization and use $(D_t^- \phi)_{i,j,k}$ to approximate $\dot \phi_{i,j}(t)$. Then, the saddle point formula~\eqref{eqt:saddle_semidiscrete_m_2d} becomes (note that the index $k$ for $\rho$, $m^1$, $m^2$ are from $1$ to $n_t-1$, while the range of $k$ for $\phi$ is from $1$ to $n_t$)
\begin{equation}\label{eqt:saddle_m_full_discrete_2d}
{\scriptsize
\begin{split}
\min_{\substack{\phi_{i,j,k} \forall i,j,k\\ \phi_{i,j,1}=g(x_i,y_j)}} \max_{\substack{\rho_{i,j,k}, m_{i,j,k}^1, m_{i,j,k}^2 \forall i,j,k\\ \rho_{i,j,k}\geq 0}}  
\sum_{i=1}^{n_x}\sum_{j=1}^{n_y} \sum_{k=1}^{n_t-1}\Big(\rho_{i,j,k}(D_t^-\phi)_{i,j,k+1} + m_{i,j,k}^1 (D_x^+\phi)_{i,j,k+1} + m_{i,j,k}^2 (D_y^+\phi)_{i,j,k+1}\quad\quad\quad\quad  \\
- \epsilon\rho_{i,j,k} (D_{xx}\phi + D_{yy}\phi)_{i,j,k+1}\Big) - \hat L\left((x_i)_i, (y_j)_j, (t_k)_k, (\rho_{i,j,k})_{i,j,k}, (m^1_{i,j,k})_{i,j,k}, (m_{i,j,k}^2)_{i,j,k}\right) - \frac{c}{\Delta t}\sum_{i=1}^{n_x} \sum_{j=1}^{n_y} \phi_{i,j,n_t},
\end{split}
}
\end{equation}
where we abuse the notation $\hat L$ and let $\hat L\colon \R^{n_x}\times \R^{n_y}\times \R^{n_t} \times \left(\R^{n_x\times (n_t-1)}\right)^3\to\R$ be defined by
\begin{equation*}
{\scriptsize
\begin{split}
\hat L\left((x_i)_i, (y_j)_j, (t_k)_k, (\rho_{i,j,k})_{i,j,k}, (m^1_{i,j,k})_{i,j,k}, (m_{i,j,k}^2)_{i,j,k}\right) = \max_{d_{i,j,k}^1, d_{i,j,k}^2\forall i,j,k} \sum_{k=1}^{n_t-1} \sum_{i=1}^{n_x} \sum_{j=1}^{n_y}  \Big(m^1_{i,j,k}d^1_{i,j,k} + m^2_{i,j,k}d^2_{i,j,k} \\
- \rho_{i,j,k} \hat H(x_i, y_j, t_{k+1}, d_{i,j,k}^1, d_{i-1,j,k}^1, d_{i,j,k}^2, d_{i,j-1,k}^2)\Big).
\end{split}
}
\end{equation*}

Similar to the one-dimensional case, this formula is particularly well-suited for Hamiltonians that possess specific structures, such as being separable and shifted $1$-homogeneous with respect to $p$, as we will elaborate on in the subsequent section.

\subsection{Two-dimensional special case: separable and shifted 1-homogeneous Hamiltonian}\label{appendix:L1_2d}
In this section, we consider the two-dimensional HJ PDE whose Hamiltonian is in the form of
\begin{equation}\label{eqt:H_L1_2d}
H(x,y,t,p) = \cinh(x,y,t)\|p\|_1 + f(x,y,t) = \cinh(x,y,t)(|p_1| + |p_2|) + f(x,y,t),
\end{equation}
with the assumption that $\cinh(x,y,t) > 0$ for any $(x,y)\in\Omega$, $t\in [0,T]$.
Although we consider this specific case, the derivation can be generalized to a Hamiltonian which is separable (see Remark~\ref{rem:separable_semidisc_2d}), convex, and 1-homogeneous with respect to $p$.

We adopt the Engquist-Osher scheme~\cite{osher1988fronts,engquist1980stable,engquist1981one}, where we set $\hat H(x,y,t, p^{1,+}, p^{1,-}, p^{2,+}, p^{2,-})$ to be $\hat H_{-}(x,y,t, p^{1,+}) + \hat H_{+}(x,y,t, p^{1,-}) + \hat H_{-}(x,y,t, p^{2,+}) + \hat H_{+}(x,y,t, p^{2,-}) + f(x,y,t)$, with $\hat H_-(x,y,t, p) = \cinh(x,y,t)\max\{-p,0\}$ and $\hat H_+(x,y, t,p) = \cinh(x,y,t)\max\{p,0\}$. 
In this context, the computation takes the following form (where $\rho_{i,j}$, $m^1_{i,j}$, $m^2_{i,j}$, $d^1_{i,j}$, $d^2_{i,j}$ are used instead of $\rho_{i,j,k}$, $m^1_{i,j,k}$, $m^2_{i,j,k}$, $d^1_{i,j,k}$, $d^2_{i,j,k}$ for simplified notation):
\begin{equation}
{\scriptsize
\begin{split}
&\max_{d^1,d^2\in \R^{n_x\times n_y}} \sum_{i,j} m^1_{i,j}d^1_{i,j} + m^2_{i,j}d^2_{i,j} - \rho_{i,j}\hat H(x_i, y_j, t_{k+1}, d^1_{i,j}, d^1_{i-1,j}, d^2_{i,j}, d^2_{i,j-1})\\
=& \max_{d^1,d^2\in \R^{n_x\times n_y}} \sum_{i,j} \Big(m^1_{i,j}d^1_{i,j} + m^2_{i,j}d^2_{i,j} - \rho_{i,j}\hat H_-(x_i, y_j, t_{k+1}, d^1_{i,j}) - \rho_{i+1,j}\hat H_+(x_{i+1}, y_j, t_{k+1}, d^1_{i,j}) - \rho_{i,j}\hat H_-(x_i, y_j, t_{k+1}, d^2_{i,j}) \\
&\quad\quad\quad\quad\quad\quad\quad\quad\quad\quad\quad\quad\quad\quad\quad\quad
- \rho_{i,j+1}\hat H_+(x_{i}, y_{j+1}, t_{k+1}, d^2_{i,j}) - f(x_i, y_j, t_{k+1})\Big)\\
= & -\sum_{i,j}f(x_i, y_j, t_{k+1}) + \ind_{[-\rho_{i,j}\cinh(x_i,y_j,t_{k+1}), \rho_{i+1,j}\cinh(x_{i+1},y_j,t_{k+1})]}(m_{i,j}^1)
+ \ind_{[-\rho_{i,j}\cinh(x_i,y_j,t_{k+1}), \rho_{i,j+1}\cinh(x_{i},y_{j+1},t_{k+1})]}(m_{i,j}^2).
\end{split}
}
\end{equation}
Then, the saddle point problem~\eqref{eqt:saddle_m_full_discrete_2d} becomes
\begin{equation}\label{eqt:saddle_m_full_discrete_L1_2d}
{\scriptsize
\begin{split}
\min_{\substack{\phi\\ \phi_{i,j,1}=g(x_i,y_j)}} \max_{\substack{\rho, m^1, m^2 \\ \rho_{i,j,k}\geq 0}}  
\Bigg\{\sum_{i=1}^{n_x}\sum_{j=1}^{n_y} \sum_{k=1}^{n_t-1}\Big(\rho_{i,j,k}(D_t^-\phi)_{i,j,k+1} + m_{i,j,k}^1 (D_x^+\phi)_{i,j,k+1} + m_{i,j,k}^2 (D_y^+\phi)_{i,j,k+1}
- \epsilon\rho_{i,j,k} (D_{xx}\phi + D_{yy}\phi)_{i,j,k+1}\\ - \rho_{i,j,k} f(x_i,y_j,t_k)\Big) - \frac{c}{\Delta t}\sum_{i=1}^{n_x} \sum_{j=1}^{n_y} \phi_{i,j,n_t} \colon 
-\rho_{i,j,k}\cinh(x_i, y_j, t_{k+1})\leq m^1_{i,j,k}\leq  \rho_{i+1,j,k}\cinh(x_{i+1},y_j,t_{k+1}), \\
-\rho_{i,j,k}\cinh(x_i, y_j, t_{k+1})\leq m^2_{i,j,k}\leq  \rho_{i,j+1,k}\cinh(x_{i},y_{j+1}, t_{k+1}), \forall  i=1,\dots, n_x; j=1,\dots, n_y; k=1,\dots, n_t-1\Bigg\}.
\end{split}
}
\end{equation}
The update procedure for $\phi$ remains as previously described. We will now detail the process of deriving updates for $\rho$ and $m$.

\begin{algorithm}[htbp]
\SetAlgoLined
\SetKwInOut{Input}{Inputs}
\SetKwInOut{Output}{Outputs}
\Input{Stepsize $\tau, \sigma>0$, error tolerance $\delta>0$, and the maximal iteration number $N$.}
\Output{Solution to the saddle point problem~\eqref{eqt:saddle_m_full_discrete_L1_2d}.}
For each $i=1,\dots, n_x$; $j=1,\dots, n_y$, initialize the matrices by $\phi_{i,j,k}^0=g(x_i,y_j)$ for $k=1,\dots, n_t$, $\rho_{i,j,k}^0 = c$, $m_{i,j,k}^{1,0} = m_{i,j,k}^{2,0} = 0$ for $k=1,\dots, n_t-1$.

 \For{$\ell = 0,1,\dots,N-1$}{
 Update the matrix $\phi_{i,k}$ for $i=1,\dots, n_x$; $j=1,\dots, n_y$; $k=1,\dots, n_t$ by
 \begin{equation}
 {\small
 \begin{split}
 \phi^{\ell+1}
 &= \phi^\ell + \tau (-D_{tt} - D_{xx} - D_{yy})^{-1}\left(D_t^+\rho^\ell + D_x^- m^{1,\ell} + D_y^- m^{2,\ell} + \epsilon D_{xx} \rho^\ell + \epsilon D_{yy} \rho^\ell\right),
 \end{split}
 }
 \end{equation}
 where $(-D_{tt} - D_{xx}-D_{yy})^{-1}f$ (for a tensor $f$ with elements $f_{i,j,k}$ and the linear operator $D_{tt}f = (\frac{f_{i,j,k-1} -2f_{i,j,k} + f_{i,j,k+1}}{\Delta t^2})_{i,j,k}$) denotes the solution $u$ to the linear system $- (D_{tt} u)_{i,j,k+1} - (D_{xx} u + D_{yy}u)_{i,j,k+1} = f_{i,j,k}$ for all $i=1,\dots, n_x$; $j=1,\dots, n_y$; $k=1,\dots, n_t-1$ with periodic spatial condition, Dirichlet initial condition $u_{i,j,1}=0$, and Neumann terminal condition $u_{i,j,n_t+1}=u_{i,j,n_t}$.
 
\If{$\sum_{i=1}^{n_x}\sum_{j=1}^{n_y}\sum_{k=2}^{n_t}|(D_t^- \phi^{\ell+1})_{i,j,k} + \hat H(x_i,y_j, t_{k}, (D_x^+ \phi^{\ell+1})_{i,j,k}, (D_x^- \phi^{\ell+1})_{i,j,k}, (D_y^+ \phi^{\ell+1})_{i,j,k}, (D_y^- \phi^{\ell+1})_{i,j,k}) - \epsilon (D_{xx} \phi^{\ell+1} + D_{yy} \phi^{\ell+1})_{i,j,k}| \leq \delta$}{
   Return $\phi^{\ell+1}$.
   }
 
 Set $\bar\phi^{\ell+1} = 2\phi^{\ell+1} - \phi^{\ell}$.
 
Update $\rho$ by
\begin{equation*}
{\scriptsize
\begin{split}
    \rho^{\ell+1}_{i,k} &= \argmin_{y\in \mathcal{A}_{i,j,k}^\ell}
    (y-\alpha_{i,j,k}^\ell)^2 + \chara_{\{z^{1,\ell}_{i,j,k}\leq 0, -y\cinh(x_i, y_j,t_{k+1}) \geq z^{1,\ell}_{i,j,k}\}}(y\cinh(x_i, y_j,t_{k+1}) + z^{1,\ell}_{i,j,k})^2\\
&\quad\quad\quad\quad\quad\quad\quad\quad + \chara_{\{z^{1,\ell}_{i-1,j,k} > 0, y\cinh(x_{i-1}, y_j, t_{k+1}) \leq z^{1,\ell}_{i-1,j,k}\}}(y\cinh(x_{i-1}, y_j, t_{k+1}) - z^{1,\ell}_{i-1,j,k})^2\\
&\quad\quad\quad\quad\quad\quad\quad\quad + \chara_{\{z^{2,\ell}_{i,j,k}\leq 0, -y\cinh(x_i, y_j, t_{k+1}) \geq z^{2,\ell}_{i,j,k}\}}(y\cinh(x_i, y_j, t_{k+1}) + z^{2,\ell}_{i,j,k})^2\\
&\quad\quad\quad\quad\quad\quad\quad\quad + \chara_{\{z^{2,\ell}_{i,j-1,k} > 0, y\cinh(x_i, y_{j-1}, t_{k+1}) \leq z^{2,\ell}_{i,j-1,k}\}}(y\cinh(x_i, y_{j-1}, t_{k+1}) - z^{2,\ell}_{i,j-1,k})^2,
\end{split}
}
\end{equation*}
where $z^{1,\ell}_{i,j,k} = m^{1,\ell}_{i,j,k} + \sigma (D_x^+\bar\phi^{\ell+1})_{i, j,k+1}$, $z^{2,\ell}_{i,j,k} = m^{2,\ell}_{i,j,k} + \sigma (D_y^+\bar\phi^{\ell+1})_{i, j,k+1}$, 
$\alpha_{i,j,k}^\ell = \rho^{\ell}_{i,j,k} + \sigma ((D_t^-\bar\phi^{\ell+1})_{i,j,k+1} + f(x_i, y_j,t_{k+1}) - \epsilon (D_{xx}\bar\phi^{\ell+1} + D_{yy}\bar\phi^{\ell+1})_{i,j, k+1})$, and $\delta_C$ denotes a function which takes the value $1$ if the constraint in $C$ is satisfied, and $0$ otherwise.
This optimization problem can be solved by comparing the function values of several candidates in $\mathcal{A}_{i,j,k}^\ell$ defined in~\eqref{eqt:candidates_rho_L1_2d}. For details, see Section~\ref{sec:update_rho_m_L1_2d}.

Update $m_{i,j,k}^1$ and $m_{i,j,k}^2$ for all $i,j,k$ ($k<n_t$) by
\begin{equation*}
\begin{split}
m^{1,\ell+1}_{i,j,k} &= \max\left\{\min\left\{z^{1,\ell}_{i,j,k}, \rho^{\ell+1}_{i+1,j,k}\cinh(x_{i+1},y_j,t_{k+1})\right\}, -\rho^{\ell+1}_{i,j,k}\cinh(x_i, y_j,t_{k+1})\right\},\\
m^{2,\ell+1}_{i,j,k} &= \max\left\{\min\left\{z^{2,\ell}_{i,j,k}, \rho^{\ell+1}_{i,j+1,k}\cinh(x_{i},y_{j+1}, t_{k+1})\right\}, -\rho^{\ell+1}_{i,j,k}\cinh(x_i, y_j,t_{k+1})\right\}.
\end{split}
\end{equation*}
where $z_{i,j,k}^{1,\ell}$ and $z_{i,j,k}^{2,\ell}$ are defined in the last step.
}

 Return $\phi^{N}$.
 \caption{The proposed algorithm for solving~\eqref{eqt:saddle_m_full_discrete_L1_2d}\label{alg:pdhg_full_m_2d_L1}}
\end{algorithm}

\subsubsection{The detailed derivation of updating $\rho$ and $m$ in Algorithm~\ref{alg:pdhg_full_m_2d_L1}} \label{sec:update_rho_m_L1_2d}
Define $z$ and $\alpha$ as in Algorithm~\ref{alg:pdhg_full_m_2d_L1}. With a similar computation as in Section~\ref{sec:update_rho_m_L1_2d}, we get explicit formulas for $m^{\ell+1}_{1,i,j}$ and $m^{\ell+1}_{2,i,j}$ using $\rho$:
\begin{equation}\label{eqt:update_m_2d}
\begin{split}
m^{1,\ell+1}_{i,j,k} &= \max\left\{\min\left\{z^{1,\ell}_{i,j,k}, \rho_{i+1,j,k}\cinh(x_{i+1},y_j,t_{k+1})\right\}, -\rho_{i,j,k}\cinh(x_i, y_j,t_{k+1})\right\},\\
m^{2,\ell+1}_{i,j,k} &= \max\left\{\min\left\{z^{2,\ell}_{i,j,k}, \rho_{i,j+1,k}\cinh(x_{i},y_{j+1}, t_{k+1})\right\}, -\rho_{i,j,k}\cinh(x_i, y_j,t_{k+1})\right\}.
\end{split}
\end{equation}
In other words, we get
\begin{itemize}
    \item If $z^{1,\ell}_{i,j,k} >0$, $m^{1,\ell+1}_{i,j,k} = \min\{z^{1,\ell}_{i,j,k}, \rho_{i+1,j,k}\cinh(x_{i+1},y_j, t_{k+1})\}$, and hence $m^{1,\ell+1}_{i,j,k}$ only depends on $\rho_{i+1,j,k}$ in the following form (denoting this function by $G^{1,+}_{i,j,k}$)
    \begin{equation}
m^{1,\ell+1}_{i,j,k} = G^{1,+}_{i,j,k}(\rho_{i+1,j,k}) := \begin{dcases}
\rho_{i+1,j,k}\cinh(x_{i+1},y_j, t_{k+1}), & \text{if }\rho_{i+1,j,k}\cinh(x_{i+1},y_j, t_{k+1}) \leq z^{1,\ell}_{i,j,k},\\
z^{1,\ell}_{i,j,k}, & \text{if }\rho_{i+1,j,k}\cinh(x_{i+1},y_j, t_{k+1}) > z^{1,\ell}_{i,j,k}.
\end{dcases}
\end{equation}

\item If $z^{1,\ell}_{i,j,k} \leq 0$, we have $m^{1,\ell+1}_{i,j,k} = \max\{z^{1,\ell}_{i,j,k}, -\rho_{i,j,k}\cinh(x_i, y_j, t_{k+1})\}$, and hence $m^{1,\ell+1}_{i,j,k}$ only depends on $\rho_{i,j,k}$ in the following form (denoting this function by $ G^{1,-}_{i,j,k}$)
    \begin{equation}
m^{1,\ell+1}_{i,j,k} = G^{1,-}_{i,j,k}(\rho_{i,j,k}) := \begin{dcases}
-\rho_{i,j,k}\cinh(x_i, y_j, t_{k+1}), & \text{if }-\rho_{i,j,k}\cinh(x_i, y_j, t_{k+1}) \geq z^{1,\ell}_{i,j,k},\\
z^{1,\ell}_{i,j,k}, & \text{if } -\rho_{i,j,k}\cinh(x_i, y_j, t_{k+1}) < z^{1,\ell}_{i,j,k}.
\end{dcases}
\end{equation}

\item If $z^{2,\ell}_{i,j,k} >0$, we have $m^{2,\ell+1}_{i,j,k} = \min\{z^{2,\ell}_{i,j,k}, \rho_{i,j+1,k}\cinh(x_{i},y_{j+1},t_{k+1})\}$, and hence $m^{2,\ell+1}_{i,j,k}$ only depends on $\rho_{i,j+1,k}$ in the following form (denoting this function by $G^{2,+}_{i,j,k}$)
    \begin{equation}
m^{2,\ell+1}_{i,j,k} = G^{2,+}_{i,j,k}(\rho_{i,j+1,k}) := \begin{dcases}
\rho_{i,j+1,k}\cinh(x_{i},y_{j+1},t_{k+1}), & \text{if }\rho_{i,j+1,k}\cinh(x_{i},y_{j+1},t_{k+1}) \leq z^{2,\ell}_{i,j,k},\\
z^{2,\ell}_{i,j,k}, & \text{if }\rho_{i,j+1,k}\cinh(x_{i},y_{j+1},t_{k+1}) > z^{2,\ell}_{i,j,k}.
\end{dcases}
\end{equation}

\item If $z^{2,\ell}_{i,j,k} \leq 0$, we have $m^{2,\ell+1}_{i,j,k} = \max\{z^{2,\ell}_{i,j,k}, -\rho_{i,j,k}\cinh(x_i, y_j, t_{k+1})\}$, and hence $m^{2,\ell+1}_{i,j,k}$ only depends on $\rho_{i,j,k}$ in the following form (denoting this function by $ G^{2,-}_{i,j,k}$)
    \begin{equation}
m^{2,\ell+1}_{i,j,k} = G^{2,-}_{i,j,k}(\rho_{i,j,k}) := \begin{dcases}
-\rho_{i,j,k}\cinh(x_i, y_j, t_{k+1}), & \text{if }-\rho_{i,j,k}\cinh(x_i, y_j, t_{k+1}) \geq z^{2,\ell}_{i,j,k},\\
z^{2,\ell}_{i,j,k}, & \text{if } -\rho_{i,j,k}\cinh(x_i, y_j, t_{k+1}) < z^{2,\ell}_{i,j,k}.
\end{dcases}
\end{equation}
\end{itemize}
Therefore, depending on the sign of $z^{1,\ell}_{i,j,k}$, the minimizer $m^{1,\ell+1}_{i,j,k}$ is a piecewise linear function depending on either $\rho_{i,j,k}$ or $\rho_{i+1,j,k}$. Similarly, depending on the sign of $z^{2,\ell}_{i,j,k}$, the minimizer $m^{2,\ell+1}_{i,j,k}$ is a piecewise linear function depending on either $\rho_{i,j,k}$ or $\rho_{i,j+1,k}$.

Then, we plug the minimizers $m^1$ and $m^2$ into the original updating scheme to get a minimization problem of $\rho$ only:
\begin{equation}
{\scriptsize
\begin{split}
\rho^{\ell+1}_{i,j,k}
&= \argmin_{y\geq 0} (y-\alpha_{i,j,k}^\ell)^2 
  + \chara_{\{z^{1,\ell}_{i,j,k}\leq 0\}}(G^{1,-}_{i,j,k}(y) - z^{1,\ell}_{i,j,k})^2+ \chara_{\{z^{1,\ell}_{i-1,j,k} > 0\}}(G^{1,+}_{i-1,j,k}(y) - z^{1,\ell}_{i-1,j,k})^2
  \\
  &\quad\quad\quad\quad
  + \chara_{\{z^{2,\ell}_{i,j,k}\leq 0\}}(G^{2,-}_{i,j,k}(y) - z^{2,\ell}_{i,j,k})^2+ \chara_{\{z^{2,\ell}_{i,j-1,k} > 0\}}(G^{2,+}_{i,j-1,k}(y) - z^{2,\ell}_{i,j-1,k})^2\\
& = \argmin_{y\geq 0}
(y-\alpha_{i,j,k}^\ell)^2 + \chara_{\{z^{1,\ell}_{i,j,k}\leq 0, -y\cinh(x_i, y_j,t_{k+1}) \geq z^{1,\ell}_{i,j,k}\}}(y\cinh(x_i, y_j,t_{k+1}) + z^{1,\ell}_{i,j,k})^2\\
&\quad\quad\quad\quad\quad\quad\quad\quad + \chara_{\{z^{1,\ell}_{i-1,j,k} > 0, y\cinh(x_{i-1}, y_j, t_{k+1}) \leq z^{1,\ell}_{i-1,j,k}\}}(y\cinh(x_{i-1}, y_j, t_{k+1}) - z^{1,\ell}_{i-1,j,k})^2\\
&\quad\quad\quad\quad\quad\quad\quad\quad + \chara_{\{z^{2,\ell}_{i,j,k}\leq 0, -y\cinh(x_i, y_j, t_{k+1}) \geq z^{2,\ell}_{i,j,k}\}}(y\cinh(x_i, y_j, t_{k+1}) + z^{2,\ell}_{i,j,k})^2\\
&\quad\quad\quad\quad\quad\quad\quad\quad + \chara_{\{z^{2,\ell}_{i,j-1,k} > 0, y\cinh(x_i, y_{j-1}, t_{k+1}) \leq z^{2,\ell}_{i,j-1,k}\}}(y\cinh(x_i, y_{j-1}, t_{k+1}) - z^{2,\ell}_{i,j-1,k})^2.
\end{split}
}
\end{equation}
Given that each $\delta$-function can adopt values of either $0$ or $1$, there exist a total of $16$ distinct combinations, each corresponding to a unique minimizer. We integrate the positive parts of these minimizers into our pool of candidates, applying truncation due to the non-negativity of $\rho$. These $16$ candidates comprise the following ensemble:
\begin{equation*}
{\scriptsize
\begin{split}
 \mathcal{A}_1 = \Bigg\{\left(\frac{\alpha_{i,j,k}^\ell - \sum_{(v_0, v_1)\in \mathcal{C}} v_0v_1}{1 + \sum_{(v_0, v_1)\in \mathcal{C}} v_0^2}\right)_+ \colon \mathcal{C}\subseteq \Big\{(\cinh(x_i, y_j, t_{k+1}), z^{1,\ell}_{i,j,k}), (\cinh(x_{i-1}, y_j, t_{k+1}), -z^{1,\ell}_{i-1,j,k}), \\
 (\cinh(x_i, y_j, t_{k+1}), z^{2,\ell}_{i,j,k}), (\cinh(x_i, y_{j-1}, t_{k+1}), -z^{2,\ell}_{i,j-1,k})\Big\}\Bigg\}.
\end{split}
}
\end{equation*}
We also include the boundary points as in Section~\ref{sec:update_rho_m_L1_1d} to get the following set of all candidates
\begin{equation}\label{eqt:candidates_rho_L1_2d}
{\scriptsize
\begin{split}
\mathcal{A}_{i,j,k}^\ell = \Bigg\{0, \left(-\frac{z^{1,\ell}_{i,j,k}}{\cinh(x_i, y_j,t_{k+1})}\right)_+, \left(\frac{z_{i-1,j,k}^{1,\ell}}{\cinh(x_{i-1},y_j, t_{k+1})}\right)_+, \left(-\frac{z^{2,\ell}_{i,j,k}}{\cinh(x_i, y_j,t_{k+1})}\right)_+, \left(\frac{z_{i,j-1,k}^{2,\ell}}{\cinh(x_{i},y_{j-1}, t_{k+1})}\right)_+\Bigg\}\cup \mathcal{A}_1.
\end{split}
}
\end{equation}
Then, the minimizer $\rho_{i,j,k}^{\ell+1}$ is chosen from these candidates by comparing their objective function values.

\end{document}